\documentclass[a4paper]{article}
\usepackage{RR}
\usepackage{hyperref}
\usepackage{amssymb}
\usepackage{amsmath}
\RRdate{September 2007}
\RRauthor{% les auteurs
Houssem Haddar%\thanks{}%
  \and
Patrick Joly%\thanks{Footnote for second author}%
  \and
Hoai Minh Nguyen%\thanks{Footnote for second author}%
}
\authorhead{Haddar \& Joly \& Nguyen}
\RRtitle{Construction and analysis of approximate models for electromagnetic scattering from imperfectly conducting scatterers}
\RRetitle{Construction and analysis of approximate models for electromagnetic scattering from imperfectly conducting scatterers}
\titlehead{GIBC for imperfectly conducting scatterers}
\RRresume{Ce rapport est d\'edi\'e \`a la construction et l'analyse des
  conditions d'imp\'edances g\'en\'eralis\'ees (GIBCs) utilis\'ees en
  \'electromagn\'etisme comme
  approximations d'ordre sup\'erieur \`a une simple condition de conducteur
  parfait.  Nous nous int\'eressons ici au probl\`eme 3-D mod\'elis\'e par les
  \'equations de Maxwell en r\'egime harmonique. L'obtention des GIBCs se base sur un d\'eveloppement
  asymptotique par rapport \`a la profondeur de peau \`a l'int\'erieur de l'obstacle
  fortement conducteur. Ce d\'eveloppement est justifi\'e \`a tout ordre par
  des estimations d'erreurs sur les s\'eries tronqu\'ees, et
  nous donnons des formules explicites de ses termes jusqu'\`a l'ordre 3.
  Ces formules sont ensuite utilis\'ees pour obtenir les expressions des GIBCs. Les probl\`emes
  aux limites associ\'es aux GIBCs sont analys\'es et nous \'etablirons des
  estimations d'erreurs relatives \`a ces mod\`eles approch\'ees en terme de la conductivit\'e.
}
\RRabstract{This report is dedicated to the construction and analysis of
  so-called Generalized Impedance Boundary Conditions (GIBCs) used in electromagnetic
  scattering problems from imperfect conductors as higher order
  approximations of a perfect conductor condition. We consider here the 3-D case
  with Maxwell equations in a harmonic regime. The construction of GIBCs is
  based 
  on a scaled asymptotic expansion with respect to the skin depth. The
  asymptotic expansion is theoretically justified at any order and we give
  explicit expressions till the third order. These expressions are used to
  derive the GIBCs.   The associated boundary value problem is analyzed and error
  estimates are obtained in terms of the skin depth.}
\RRmotcle{GIBC, conditions aux limites effectives, forte conductivit\'e, diffraction electromagnetique, equations de Maxwell, couches
  limites, d\'eveloppement asymptotique, mod\`eles approch\'es }
\RRkeyword{GIBC, effective boundary conditions, highly conducting scatterers, electromagnetic scattering, Maxwell's equations, boundary
  layers, scaled asymptotic expansions, approximate models}
\RRprojets{POEMS et DeFI}
\RRtheme{\THNum} % cas de 5 themes
 \URRocq % pour ceux qui sont au centre de la France
\begin{document}
\makeRR   % cas d'un rapport de recherche
%% \makeRT % cas d'un rapport technique.
%% a partir d'ici, chacun fait comme il le souhaite
%\documentclass[11pt]{article}
%\usepackage{RR}
%\usepackage{a4wide}
%\usepackage{hyperref}
%\usepackage{amssymb}
%\usepackage{amsmath}
%\usepackage[dvips]{color, graphicx}
%\usepackage{subeqnar}

%\begin{document}
\newcommand{\eqn}{eqnarray}
\newcommand{\nc}{\newcommand}
\newcommand{\dsp}{\displaystyle}
\newcommand{\Oint}{\Omega_i}
\newcommand{\Oext}{\Omega_e}
\newcommand{\ptre}{\delta}
\newcommand{\eps}{\varepsilon}
\newcommand{\ui}{u_i}
\newcommand{\ue}{u_e}
\newcommand{\uue}{u^\ptre}
\newcommand{\uie}{u^\ptre_i}
\newcommand{\uee}{u^\ptre_e}
\newcommand{\fu}{\varphi}
\newcommand{\fue}{\fu^\ptre}
\newcommand{\ubke}{{\bar u_k^\ptre}}
\newcommand{\vke}{{v_k^\ptre}}
\newcommand{\lri}{\longrightarrow}
\newcommand{\bE}{{\bf E}}
\newcommand{\bH}{{\bf H}}
\newcommand{\wH}{{\widetilde H}}
\newcommand{\wu}{{\widetilde u}}
\newcommand{\hu}{{\widehat u}}
%%%%%%%%%%%%%%%%%%%%%%%%%%%%%%%%%%%%%%%%%%%%%%
\nc{\R}{\mathbb{R}}
\nc{\N}{\mathbb{N}}
\nc{\C}{\mathbb{C}}
\nc{\Gamt}{{\scriptscriptstyle{\Gamma}}}
\nc{\rt}{\mathrm{curl}} 
\nc{\dv}{\mathrm{div}} 
\nc{\divg}{\mbox{\rm div}_\Gamt \,} 
\nc{\grag}{\nabla_\Gamt \,}
\nc{\rotg}{\rt_\Gamt \,} 
\nc{\Rotg}{\vec{\rt_\Gamt} \,} 
\nc{\rotgm}{\rt_\Gamt^{\scriptscriptstyle M} \,} 
\nc{\Rotgm}{\vec\rt_\Gamt\,\!\!\!\!^{\scriptscriptstyle M} \,} 
\nc{\Lapg}{\Delta_\Gamt \,} 
\nc{\Deltag}{\vec\Delta_\Gamt \,} 
\nc{\xgg}{{x}_\Gamt}
\nc{\Ig}{{I}_\Gamt}
\nc{\grad}{\mbox{ grad }} 
\nc{\lE}{{\bf E}}
\nc{\lH}{{\bf H}} 
\nc{\mE}{{\mathbb E}}
\nc{\eg}{E}
\nc{\hg}{H}
\nc{\egie}{\eg^\ptre_i}
\nc{\egee}{\eg^\ptre_e} 
\nc{\hgie}{\hg^\ptre_i} 
\nc{\hgee}{\hg^\ptre_e}
\nc{\nEi}{{\widehat \mE}_i}
\nc{\egek}{\eg^k_e} 
\nc{\hgek}{\hg^k_e} 
\nc{\egik}{\eg^k_i}
\nc{\hgik}{\hg^k_i} 
\nc{\bfH}{{\bf H}} 
\nc{\oV}{\overline{V}}
\nc{\upE}{\tilde{\bE}} 
\nc{\upH}{\tilde{\bH}} 
\nc{\Ree}{{\cal R}e \;}
\nc{\Ime}{{\cal I}m \;}
\nc{\Ri}{\Rightarrow}
\nc{\tE}{\tilde{E}} 
\nc{\oE}{\overline{E}} 
\nc{\ege}{\eg_e}
\nc{\hge}{\hg_e} 
\nc{\egi}{\eg_i} 
\nc{\hgi}{\hg_i}
\nc{\oupE}{\overline{\upE}} 
\nc{\bfT}{{\bf T}}
\nc{\ttT}{\texttt{T}} 
\nc{\Etke}{\eg^{\ptre, k}}
\nc{\Htke}{\hg^{\ptre, k}} 
\nc{\Eke}{\eg^{\ptre, k}}
\nc{\Hke}{\hg^{\ptre, k}} 
\nc{\egeek}{\eg_e^{\ptre, k}}
\nc{\hgeek}{\hg_e^{\ptre, k}} 
\nc{\egerho}{\eg^{\ptre,\rho}}
\nc{\hgerho}{\hg^{\ptre,\rho}}
\nc{\oEerho}{\overline{\eg^{\ptre,\rho}}}
\nc{\oHerho}{\overline{\hg^{\ptre,\rho}}} 
\nc{\tu}{\tilde{u}}
\nc{\dEke}{{\cal E}^{\ptre,k}} 
\nc{\dHke}{{\cal H}^{\ptre,k}}
\nc{\bfdEk}{{\cal E}^{\ptre,k}} 
\nc{\wEk}{\widehat{\mE}^{\ptre,k}}
\nc{\tdEke}{\widetilde{\cal E}_e^{\ptre,k}} 
\nc{\tdHke}{\widetilde{\cal H}_e^{\ptre,k}} 
\nc{\bfdEke}{{\cal E}_e^{\ptre,k}} 
\nc{\bfdHke}{{\cal H}_e^{\ptre,k}} 
\nc{\bfdEki}{{\cal E}_i^{\ptre,k}} 
\nc{\bfdHki}{{\cal H}_i^{\ptre,k}} 
\nc{\bfdHk}{{\cal H}^{\ptre,k}} 
\nc{\dEkee}{\delta \eg_e^{\ptre, k}} 
\nc{\dHkee}{\delta \hg_e^{\ptre, k}} 
\nc{\dEkei}{\delta \eg_i^{\ptre, k}}
\nc{\dHkei}{\delta \hg_i^{\ptre,k}} 
\nc{\dEe}{{\cal E}_{3}^\ptre} 
\nc{\dHe}{{\cal H}_{3}^\ptre}
\nc{\olambda}{\overline{\lambda}} 
\nc{\ophi}{\overline{\phi}}
\nc{\oU}{\overline U} 
\nc{\supp}{\mbox{ supp }}
\nc{\dist}{\mbox{dist}}
\nc{\sqrti}{{\sqrt{i}}}
\nc{\hcurl}[1]{H(\rt,{#1})}
\nc{\hcurlt}[1]{\widetilde{H}(\rt,{#1})}
\nc{\hcurltz}[1]{\widetilde{H}_0(\rt,{#1})}
\nc{\Ee}{E_e}
\nc{\He}{H_e} \nc{\Ei}{E_i} \nc{\Hi}{H_i}
\nc{\Eie}{E^\ptre_i}
\nc{\Eee}{E^\ptre_e} \nc{\Hie}{H^\ptre_i} \nc{\Hee}{H^\ptre_e}
\nc{\Eek}{E^k_e} \nc{\Hek}{H^k_e} \nc{\Eik}{E^k_i}

%%%%%%%%%%%%%%%%%%%%%%%%%%%%%%%%%%%%%%%%%%%%%

\nc{\wt}[1]{\widetilde{#1}}
\newcommand{\nod}[2]{\left\| #1 \right\|_{#2}}
\newcommand{\nodt}[2]{\| #1 \|_{#2}}

\newtheorem{remark}{Remark}[section]
\newtheorem{definition}{Definition}[section]
\newtheorem{theorem}{Theorem}[section]
\newtheorem{lemma}{Lemma}[section]
\newtheorem{corollary}{Corollary}[section]

\newcommand{\proofend}{\hfill $\Box$ \vspace{2mm}}
\newcommand{\proof}{{\it Proof.~}}

%%%%%%%%%%%%% Fin des commandes de Houssem

%\title{Construction and analysis of approximate models for electromagnetic scattering form imperfectly conducting scatterers}
%\author{H. Haddar, H.M. Nguyen, P. Joly}
%\date{10/09/2005}
%\maketitle

\section{Introduction}

Generalized Impedance Boundary Conditions ($GIBC$) have become
a rather classical notion in the mathematical modeling of wave
propagation phenomena (see for instance, \cite{H-R-95} and, \cite{S-V-95}).
%. As an example, the book by Senior and Volakis, 
%``Approximate Boundary Conditions in Electromagnetics'' \cite{},
%is entirely devoted to this subject.  
They are used in electromagnetism for time harmonic scattering
problems from obstacles that are partially
or totally penetrable. The general idea is to replace the use of an ``exact
model'' inside (the penetrable part of) the obstacle by approximate
boundary conditions (also called equivalent or effective conditions). This
idea is pertinent if the boundary 
condition can be easily handled numerically, for instance when it is local. 
% The same type of idea led to the construction of local absorbing
% boundary conditions for the wave equation, \cite{E-M-77} 
% \cite{C-92} or more
% recently to the construction of On Surface Radiation Conditions, \cite{A-B-B-99}, \cite{A-B-01}.
% \\[12pt]
The  diffraction problem of electromagnetic waves by perfectly
conducting obstacles coated with a thin layer of dielectric material
is well suited for the notion of impedance conditions: 
due to the small (typically with respect to the wavelength) thickness
of the coating, the effect of the layer on the exterior medium is, as
a first approximation, local (see for instance, \cite{S-V-95}, \cite{H-R-95}, \cite{ENGUNED}, \cite{BENDLEMR}, 
\cite{A-C-91}).
\\[12pt]
The application we consider here is the diffraction of
 waves by highly conducting
 materials in electromagnetism. In such a case, it is the well-known  skin
 effect that creates a ``thin layer'' phenomenon. The high conductivity
 limitates the penetration of the wave  to a  boundary layer whose depth is
 inversely proportional to the square root of its magnitude. Then, here
 again, the effect of the obstacle is, as
a first approximation, local. 
% The numerical results presented in Figures \ref{fig1}-\ref{fig4} of
% section \ref{sect-num} illustrate this skin
% effect phenomena.
\\[12pt]
The first effective boundary
conditions for highly absorbing obstacles was proposed by Leontovich. This
condition ``sees''  only locally  the tangent plane to the
frontier. Later, Rytov, \cite{Rytov}, 
\cite{S-V-95} proposed an extension of the
Leontovitch condition, and his analysis was already based on the principle of
asymptotic expansions with respect to the small parameter in the problem: the
skin depth $\delta$. More recently, Antoine-Barucq-Vernhet \cite{A-B-V-01}
proposed a new derivation of such conditions based on the technique  of
pseudo-differential operator expansions. However, in all these works, the
rigorous mathematical justification of the resulting impedance conditions was
not treated.\\[12pt]
This paper is the continuation of the work in \cite{Haddar-Joly-Nguyen1}, in which we considered
the case of the scalar wave equation. Our objective is to extend the results of
to the case of 3D Maxwell's equations by constructing and analyzing GIBC's of
order 1, 2 and 3 (with respect to the skin depth, the small parameter of the
problem). These conditions are of impedance type (or $H-to-E$ nature): they
relates the tangential traces of the electric and magnetic fields via a 
local impedance operator.\\[12pt]
As in \cite{Haddar-Joly-Nguyen1}, the construction of the approximate conditions relies on
an asymptotic expansion of the exact solution, based on a scaling technique
and a boundary layer expansion in the neighborhood of the boundary of the
scatterer. If the organization of this paper contents is similar to 
\cite{Haddar-Joly-Nguyen1}, its technicality is much higher. Moving from the scalar wave
equation to the Maxwell system increases considerably the complexity of the
problem at two levels.
\begin{itemize}
\item The first one is linked with the algebra involved in the formal construction of the asymptotic 
  expansion of the exact solution (see Section \ref{derivation}). This
  is  essentially due to  the vectorial  nature of the unknowns and the
  expression of the curl operator in a parametric coordinates system (see
  Section \ref{preliminary}). The latter is  based on the
  formulas proposed in \cite{Haddar-Joly1}  with some simplifications.
\item The second one is related to the mathematical analysis on the GIBCs. This
  is not only due to the fact that we have to deal with usual functional
  analysis difficulties linked to Maxwell equations (in particular trace
  operators and compact
  embedding properties - see Sections \ref{errortrunk} and \ref{Sec-Analysis}
  and Appendix \ref{AppA}) but also 
  because we have to 
  face some new difficulties in the case of the third order condition. The
  tangential differential operators that would naturally appear  in the
  construction of the third order condition have not the good ``sign
  properties'' to be able to guarantee the existence of the approximate
  solution and the convergence (at optimal order) to the true solution. This
  leads us to apply various regularization procedures to construct the modified
  third order conditions (see Section \ref{modif}).
\end{itemize}
Our objectives in this work are essentially theoretical. The numerical
pertinence of obtained conditions have already been demonstrated in \cite{DHJ06-S}
where, in particular, the interest of using a third order condition rather
than a first or a second order condition is clearly shown. 
\\[12pt]
The outline of the article is as follows. Section
\ref{Sec_Model} contains a description of the physical and mathematical diffraction problem at study with
some basic stability properties of the solutions and asymptotic estimates with respect to
the conductivity.  We state the
main results of our paper in Section \ref{main}: the GIBCs
are presented in Sections \ref{Sec-GIBC1} and \ref{modif} while 
the corresponding error estimates are given in Section \ref{Sec_Error}. 
The formal construction of the asymptotic expansion is given in Section
\ref{derivation}.  This construction is
rigorously justified in Section
\ref{errortrunk} by proving optimal error estimates at each order. The last
section is dedicated to the study of the boundary value problems associated
with the GIBCs as well as the proof of optimal error estimates between these
solutions and truncated asymptotic expansions. The main result of our paper is
obtained as a combination of the results of Section~\ref{Sec-Analysis} and
Section \ref{errortrunk}. Some non standard technical results related to the $H(\rt)$ space 
(appropriate trace inequalities and special compact embedding properties) 
that may have their own interest have been gathered in Appendix A.
\section{Description of the physical model} \label{Sec_Model}
Let $\Oint$ be an open bounded domain in $\R^3$ with connected complement,
occupied by a homogeneous conducting medium. We denote by $\Gamma$ the
boundary of  $\Oint$ and assume that this boundary is a $C^{\infty}$ manifold.
We are interested in computing the  electromagnetic diffracted wave when the
conductivity of the medium, denoted by  $\sigma^\ptre$, is sufficiently high
($\ptre$ denotes a small parameter). More precisely we assume that $\sigma^\ptre \rightarrow
\infty$ as $\ptre \rightarrow 0$ and would like to study the  asymptotic
behavior of the diffracted electromagnetic field as $\ptre \rightarrow 0$ in
order to derive efficient
approximate models to compute the diffracted waves.
\\[12pt]
We assume that the exterior domain is homogeneous and the time and space
scales are chosen such  that the wave speed is 1 in this medium. The electromagnetic wave propagation is therefore governed by the following Maxwell's equations:
\begin{eqnarray*}
\left\{ \begin{array}{ll} \eps(x) \dsp \frac{\partial \lE^\ptre}{\partial
t} + \sigma^{\ptre}(x) \lE^\ptre - \rt\, \lH^\ptre
= F, & \mbox{ in } \Omega,  \\[6pt]
\dsp \frac{\partial \lH^\ptre }{\partial t} + \rt\, \lE^\ptre = 0, &
\mbox{ in } \Omega,
\end{array} \right.
\end{eqnarray*}
where $\Omega \subset \R^3$ denotes the propagative medium that we shall
assume to be regular and simply connected with connected boundary (for
instance an open ball), the functions $\sigma^{\ptre}(x)$ and $\eps(x)$ are
defined by:
\begin{equation*}
(\eps, \sigma^\ptre)(x) = \left\{ \begin{array}{ll}
(1, 0), \; & \mbox{  in  } \; \Oext, \\[6pt]
(\eps_r, \sigma^\ptre), \; & \mbox{  in  } \; \Oint \; 
\end{array} \right.
\end{equation*}
where $\Oext = \Omega \setminus \overline \Oint$ and where  $\eps_r > 0$
denotes the relative electric permittivity of the conducting medium. The
right-hand side  $F$ denotes some source term that we shall assume to be
harmonic in time : $F(x,t) = \mbox{Re}\left\{f(x)
~\mbox{exp} (i \omega t)\right\} $, where $\omega > 0$ denotes a given
frequency, and where $\mbox{Re}(z)$ denotes the real part of $z$. Hence, the
solutions are also time  harmonic:
\begin{eqnarray*}
\begin{array}{ll}
\lE^{\ptre}(x,t) = \mbox{Re}  \left\{ E^{\ptre}(x) \exp (i \omega
t)
\right\}, &
\lH^{\ptre}(x,t) = \mbox{Re}  \left\{ H^{\ptre}(x) \exp (i \omega
t) \right\},
\end{array} 
\end{eqnarray*}
where  the
electromagnetic field $(E^\ptre, H^\ptre)$ is solution to the harmonic Maxwell system:
\begin{eqnarray}\label{equation:1}
\left\{ \begin{array}{llll} (i) & (i \eps \omega + \sigma^\ptre) E^\ptre - \rt\, H^\ptre  = f, ~~ &
\mbox{ in } \; \Omega,\\[6pt]
(ii) & i \omega H^\ptre + \rt\, E^\ptre = 0,~~ &\mbox{ in } \;
\Omega.
\end{array} \right.
\end{eqnarray}
We assume that the support of the source term $f$ does not touch $\Oint$. The system of equations (\ref{equation:1})
has to be complemented with a boundary condition on the exterior
boundary $\partial \Omega$, for instance we work with the
following absorbing boundary condition %(see Remarks~\ref{rem1} and~\ref{rem2})
\begin{equation} \label{impedance}
E^\ptre_T - H^\ptre \times n = g,\quad \mbox{ on } \;
\partial \Omega,
\end{equation}
where $E_T := n \times (E \times n) $, $n$ is a normal
vector to  $\partial \Omega$ directed to the exterior of $\Omega$ and $g$
denoting some possible source term.
\begin{remark} \label{rem1}
According to (\ref{impedance}), the boundary $\partial \Omega$ can
be seen as a physical absorbing boundary where a standard
impedance condition is applied. The problem (\ref{equation:1},
\ref{impedance}) can also be seen as an approximation in a bounded
domain (namely $\Omega$) of the scattering problem in $\R^3
\setminus \Omega_i$. In such a case, the boundary condition on
$\partial \Omega$ has to be understood as an (low order)
approximation of the outgoing radiation condition at infinity and $g$ a source
term linked to the incident field.
\end{remark}
\begin{remark}\label{rem2}
In this paper, we could have treated as well the scattering
problem in $\R^3 \setminus \Omega_i$. The reader will easily be
convinced that the obtained results can be extended to this case
without any major difficulty. The only difference would lie in the
reduction to a bounded domain. This additional difficulty is
purely technical and not essential in the context of this paper
whose main purpose is the treatment of the ``interior boundary''
$\Gamma$.
\end{remark}
As mentioned above we are interested in describing the asymptotic behavior of the solution for
large $\sigma^{\ptre}$. As suggested by the expression of the analytic solution where $\Oint$ is the half space, the appropriate small length parameter can be defined as: 
\begin{equation*} %\label{defeps}
\ptre := 1/\sqrt{\omega \sigma^\ptre} \quad \Longleftrightarrow \quad \sigma^\ptre =1/(\omega \ptre^2).
\end{equation*}
This small parameter defines the so-called skin depth: the ``width''  of
the penetrable region inside the conducting medium is proportional to $\ptre$. 

\bigskip

\noindent For the construction of approximate models in the exterior
domain $\Oext$  it is useful to rewrite the
problem (\ref{equation:1}-\ref{impedance}) as a transmission
problem between $(\egie, \hgie) := (E^\ptre, H^\ptre)_{|\Oint}$ and  $(\egee,
\hgee) = (E^\ptre, H^\ptre)_{|\Oext}$ as follows:
\begin{eqnarray} \label{trans_eq_1}
&&\left\{ \begin{array}{ll}
 i \omega \egee - \rt\, \hgee = f, \;\; & \mbox{ in } \; \Oext,\\[6pt]
 i \omega \hgee + \rt\, \egee = 0, \;\; & \mbox{ in } \; \Oext, \\[6pt]
  {\eg}^\ptre_{e,\,T} - \hgee \times n = g, \;\; & \mbox{ on } \;\partial \Omega\\[6pt]
\egee \times n = \egie \times n , \;\;  & \mbox{ on } \; \Gamma,
\end{array} \right.
\\[6pt]
&& \label{trans_eq_2}
\left\{ \begin{array}{ll}
(i \eps_r\omega + \frac{1}{\omega \delta^2}) \egie - \rt\, \hgie  = 0, \;\; & \mbox{ in } \; \Oint,\\[6pt]
i \omega \hgie + \rt\, \egie = 0, \;\; & \mbox{ in } \; \Oint, \\[6pt]
\hgie \times n = \hgee \times n , \;\;  & \mbox{ on }\;\Gamma. 
\end{array} \right.
\end{eqnarray}
We have chosen to split the two transmission conditions (namely the continuity
of the tangential electric and magnetic fields) in such a way that the first
one appears as a boundary condition in (\ref{trans_eq_1}) for the interior
field while  the second one appears as a boundary condition in
(\ref{trans_eq_1}) for the interior field.
Roughly speaking, the approximate models are then obtained from replacing in
system (\ref{trans_eq_2})
the exact boundary condition on $\Gamma$ by an approximate one, whose
expression is derived from seeking appropriate asymptotic expansion of the
solution in the boundary layer inside $\Oint$.

\subsection{Existence-Uniqueness-Stability}
With $\hcurl{{\cal O}}$ denoting the space of functions $V \in L^2({\cal O})^3$ such that
$\rt\, V \in L^2(U)^3$, where ${\cal O}$ is an open domain of $\R^3$,  we define
\begin{equation}\label{defhcurlt}
\hcurlt{\cal O} = \{ V \in \hcurl{\cal O} \;\; ; \;\; V_T \in L^2_t(\partial {\cal O}) \}
\end{equation}
where $V_T$ is the tangential trace of $V$ (cf. Section \ref{preliminary} for
more details), 
$L^2_t(\partial {\cal O})$ denotes the space of functions $V \in 
L^2(\partial {\cal O})^3 $ such that $V \cdot n =0 $ on $\partial {\cal O}$,
where $n$ denotes a normal to  $\partial {\cal O}$. We recall that $\hcurlt{\cal O}$
is a Hilbert space with scalar product
$$
(U,V)_{\hcurlt{\cal O}} = (U, V)_{L^2({\cal O})}+ (\rt\, U, \rt\, V)_{L^2({\cal O})} + (U_T, V_T)_{L^2_t(\partial{\cal O})}.
$$
\begin{theorem} \label{exist-stab}
For given $f \in L^2(\Omega)^3$ and $g \in L^2_t(\partial \Omega)$ there exists
an unique solution $(E^\ptre, H^\ptre) \in \hcurlt{\Omega} \times
\hcurlt{\Omega}$ satisfying (\ref{equation:1}-\ref{impedance}). Moreover, there
exists a positive constant $C$ independent of $\ptre$ such that 
\begin{equation}\label{estim-ex}
\nodt{E^\ptre}{\hcurlt{\Omega}} + \frac{1}{\ptre} \, \nodt{E^\ptre}{L^2({\Oint})}
\le C\left(\nod{f}{L^2(\Omega)}+ \nod{g}{L^2_t(\partial\Omega)}\right).
\end{equation}
\end{theorem}

\proof~ The proof of existence and uniqueness can be found in \cite{Monk}
(Theorem 4.17). The solution $E^\ptre$ is constructed (using the Helmholtz
decomposition) as 
$$E^\ptre = E^\ptre_0 + \nabla p^\ptre$$
where $p^\ptre \in H^1_0(\Omega)$ and $E^\ptre_0 \in \hcurltz{\Omega}$ where 
$$
\hcurlt{\Omega} = \{ V \in \hcurlt{\Omega}\; ; \;\; \dv
((i\sigma^\ptre-\eps \omega^2) V) = 0\},
$$
will be equipped with the norm of $\hcurlt{\Omega}$, of which it is a closed
subspace.\\[12pt] 
Taking the divergence of (\ref{equation:1}-i) one easily see that the function
$p^\ptre$  is solution to
\begin{equation}\label{variational-1}
-\omega^2 (\eps \nabla p^\ptre, \nabla \phi)_{L^2(\Omega)} + \frac{i}{\delta^2}
 ( \nabla p^\ptre, \nabla \phi)_{L^2(\Oint)} =  i \omega (f, \nabla
 \phi)_{L^2(\Omega)} 
\end{equation}
for all $\phi \in H^1_0(\Omega)$. Choosing $\phi =- p^\ptre$ then successively considering
 the
real and imaginary parts  of the resulting equality, one
easily deduces that
\begin{equation} \label{estim-ex-0}
\nodt{\nabla p^\ptre}{L^2(\Omega)} + \frac{1}{\ptre} \, \nodt{\nabla p^\ptre}{L^2({\Oint})}
\le C\nod{f}{L^2(\Omega)},
\end{equation}
for some positive constant $C$ independent of $\ptre$, which proves estimate
(\ref{estim-ex}) for the gradient part of the solution. It remains to get an
estimate for $E^\ptre_0$. The proof is based on the compact
embedding of $\hcurltz{\Omega}$ into $L^2(\Omega)^3$, whose proof can be found
in \cite{Monk}, Theorem 4.7 (which is an adaptation of the proofs in
\cite{Weber} and \cite{Costabel}). The function $E^\ptre_0$ satisfies the
variational formulation
\begin{equation}\label{variational-2}
\begin{array}{l}
 ( \rt\,E^\ptre_0,  \rt\,  \Phi)_{L^2(\Omega)} - \omega^2 (\eps\,E^\ptre_0, \Phi)_{L^2(\Omega)}
  + \frac{i}{\ptre^2} (E^\ptre_0, \Phi)_{L^2(\Oint)}
+ i\omega (E^\ptre_{0,\, T}, \Phi_T)_{L^2_t(\partial \Omega)}
\\[6pt]
=i\omega \left\{  (f, \Phi)_{L^2(\Oext)}
-  (g, \Phi_T)_{L^2_t(\partial \Omega)} \right\} + \omega^2 (\eps\,\nabla p^\ptre, \Phi)_{L^2(\Omega)}
  - \frac{i}{\ptre^2} (\nabla p^\ptre, \Phi)_{L^2(\Oint)}
\end{array}
\end{equation}
for all $\Phi \in \hcurltz{\Omega}$,  where $p^\ptre$ is the solution to (\ref{variational-1}). 
We first prove (by contradiction) that there exists a positive constant $C$
independent of $\ptre$ such that  
\begin{equation}\label{estim-ex-1}
\nodt{E^\ptre_0}{L^2({\Omega})}
\le C\left(\nod{f}{L^2(\Omega)}+ \nod{g}{L^2_t(\partial\Omega)}\right).
\end{equation}
If not, then there would exist a sequence of data $(f^\ptre)$  and $(g^\ptre)$ such
that $$\nodt{f^\ptre}{L^2(\Omega)}+ \nodt{g^\ptre}{L^2_t(\partial\Omega)} = 1$$
and such that the corresponding solution $E^\ptre_0$ satisfies
$$\nodt{E^\ptre_0}{L^2({\Omega})} \rightarrow \infty \mbox{ as } \ptre
\rightarrow 0. $$
 With an added tilda denoting a division by
$\nodt{E^\ptre_0}{L^2({\Omega})}$, one observes that 
\begin{equation}\label{norm1}
\nodt{\tilde E^\ptre_0}{L^2({\Omega})} = 1,
\end{equation}
and satisfies, for all $\Phi \in \hcurltz{\Omega}$,
\begin{equation}\label{variational-3}
\begin{array}{l}
 ( \rt\,\tilde E^\ptre_0,  \rt\,  \Phi)_{L^2(\Omega)} - \omega^2 (\eps\,\tilde E^\ptre_0, \Phi)_{L^2(\Omega)}
  + \frac{i}{\ptre^2} (\tilde E^\ptre_0, \Phi)_{L^2(\Oint)}
+ i \omega(\tilde E^\ptre_{0,\, T}, \Phi_T)_{L^2_t(\partial \Omega)}
\\[6pt]
=i\omega \left\{  (\tilde f^\ptre, \Phi)_{L^2(\Oext)}
-  (\tilde g^\ptre, \Phi_T)_{L^2_t(\partial \Omega)} \right\} + \omega^2 (\eps\,\nabla \tilde p^\ptre, \Phi)_{L^2(\Omega)}
  - \frac{i}{\ptre^2} (\nabla \tilde p^\ptre, \Phi)_{L^2(\Oint)}.
\end{array}
\end{equation}
where $p^\ptre$ is the solution to (\ref{variational-1}) with $f$ replaced by $\tilde f^\ptre$. For $ \Phi = \tilde E^\ptre_0$, one gets
\begin{equation}\label{variational-3bis}
\begin{array}{l}
 \nodt{\rt\,\tilde E^\ptre_0}{L^2(\Omega)}^2 - \omega^2 (\eps\,\tilde E^\ptre_0, \tilde E^\ptre_0)_{L^2(\Omega)}
  + \frac{i}{\ptre^2} \nodt{\tilde E^\ptre_0}{L^2(\Oint)}^2
+ i \omega\nodt{\tilde E^\ptre_{0,\, T}}{L^2_t(\partial \Omega)}^2
\\[6pt]
=i\omega \left\{  (\tilde f^\ptre, \tilde E^\ptre_0)_{L^2(\Oext)}
-  (\tilde g^\ptre, \tilde E^\ptre_{0,\, T})_{L^2_t(\partial \Omega)} \right\} + \omega^2 (\eps\,\nabla \tilde p^\ptre, \tilde E^\ptre_0)_{L^2(\Omega)}
  - \frac{i}{\ptre^2} (\nabla \tilde p^\ptre, \tilde E^\ptre_0)_{L^2(\Oint)}.
\end{array}
\end{equation}
Taking first the imaginary part (\ref{variational-3bis}) and using estimate (\ref{estim-ex-0}) and (\ref{norm1}), one proves that
\begin{equation} \label{toto-bof}
\frac{1}{\ptre} \nodt{\tilde E^\ptre_0}{L^2(\Oint)}
+  \omega\nodt{\tilde E^\ptre_{0,\, T}}{L^2_t(\partial \Omega)} \le C\left( \nodt{\tilde f^\ptre}{L^2(\Omega)}+ \nodt{\tilde g^\ptre}{L^2_t(\partial\Omega)} \right)
\end{equation}
for some positive constant $C$ independent of $\ptre$. One deduces after
taking second the real part and using again estimate (\ref{estim-ex-0}) and
(\ref{norm1}) that also  
\begin{equation} \label{toto-bof2}
\nodt{\rt\,\tilde E^\ptre_0}{L^2({\Omega})} \le C_1  +  C_2 \left( \nodt{\tilde f^\ptre}{L^2(\Omega)}+ \nodt{\tilde g^\ptre}{L^2_t(\partial\Omega)} \right)
\end{equation}
for some different positive constants $C_1$ and $C_2$ independent of
$\ptre$. It is then observed that the sequence $(\tilde E^\ptre_0)$ is bounded
in $\hcurltz{\Omega}$. 
\\[10pt]
Applying (the trace) Lemma~\ref{trace-lemma} (proved in the Appendix) to $\tilde E^\ptre_0|_{\Oint}$ one
deduces from (\ref{toto-bof}) and (\ref{toto-bof2}) that 
$$
\|\tilde E^\ptre_0\times n\|_{H^{-\frac{1}{2}}(\Gamma)} \rightarrow 0 \qquad \mbox{as }
\ptre \rightarrow 0.
$$
Also from compactness theorems in Sobolev spaces and (\ref{toto-bof2}) one
can assume (up to an extracted subsequence)  that $E^\ptre_0\times n$ is
convergent  in
${H^{-\frac{1}{2}}(\partial \Omega)}$. 
Hence, (the compactness) Lemma~\ref{NewCompact-lemma} applied to the sequence $\tilde
E^\ptre_0 |_{\Oext}$, shows that  one can  extract a subsequence (also denoted
$\tilde E^\ptre_0|_{\Oext}$) that converges to some $\tilde E$ weakly in
$\hcurltz{\Oext}$ and strongly in $L^2(\Oext)^3$. Moreover
\begin{equation}\label{eq-uni-2}
\tilde E \times n = 0\;\;  \mbox{on}  \; \Gamma.
\end{equation}
Taking $\Phi$ with support inside $\Oext$ in
(\ref{variational-3}) and taking the limit as $\delta \rightarrow 0$,  show
that  $\tilde E$ satisfies (both $\tilde f^\ptre$ and $\tilde g^\ptre$ tend to
$0$) 
\begin{eqnarray} \label{eq-uni-1}
\rt\, \rt\, \tilde E - \omega^2 \tilde E = 0 \;\; & \mbox{in}  \; \Oext
\\[6pt]
i \omega \tilde E_T +  \rt \,\tilde E \times n = 0 &  \mbox{on}  \; \partial \Omega.
\end{eqnarray}
The uniqueness of the solution in $\hcurltz{\Omega}$ to
(\ref{eq-uni-1}-\ref{eq-uni-2}) proves that $\tilde E = 0 \;\; \mbox{in}  \;
\Oext. $ Since estimate (\ref{toto-bof}) also proves  that $\nodt{\tilde
E^\ptre_0}{L^2({\Oint})}$ $\rightarrow 0$, we get that
$$1 = \|\tilde E^\ptre_0\|_{L^2(\Omega)} \rightarrow \|\tilde
 E\|_{L^2(\Oext)}$$
which contradicts  $\tilde E=0 \;\; \mbox{in}  \;
\Oext$ and proves (\ref{estim-ex-1}). 
\\[10pt]
Now take $\Phi = E^\ptre_0$ in (\ref{variational-2}) and use the imaginary part then the real part as previously done to deduce
 \begin{equation}\label{estim-ex-2}
\nodt{E^\ptre_0}{\hcurlt{\Omega}} + \frac{1}{\ptre} \, \nodt{E^\ptre_0}{L^2({\Oint})}
\le C\left(\nod{f}{L^2(\Omega)}+ \nod{g}{L^2_t(\partial\Omega)}\right).
\end{equation}
Estimate (\ref{estim-ex}) is a straightforward consequence of (\ref{estim-ex-2}) and (\ref{estim-ex-0}).\proofend

\subsection{Exponential interior decay of the solution}
It is shown that the norm of the solution in a domain strictly interior
to $\Oint$ goes to $0$ faster than any power of $\ptre$. This is
a first way to express how the interior solution
concentrates near the boundary $\Gamma$. Let us indicate that this result will
also be a consequence of the asymptotic analysis performed in next sections,
however the methodology is more complex. The proof given here is a direct one
and is independent of the subsequent analysis.  The precise result is the
following: 
\begin{theorem}\label{apriori-estimate-2}
For any $\bar\nu> 0$ small enough so that $\Oint^{\bar\nu} := \{ x \in \Oint ; B(x,\bar\nu)
  \subset \Oint\}$ is a non-empty set, where $B(x, {\bar\nu})$ denotes the closed ball
  of center $x$ 
  and radius ${\bar\nu}$, there exist two positive constants
  $C_{\bar\nu} $ 
  and $c_{\bar\nu} $ independent of $\ptre$ such that
$$\nodt{\egie}{H(\rt,\Oint^{\bar\nu})} + \nodt{\hgie}{H(\rt,\Oint^{\bar\nu})}\le C_{\bar\nu} \exp (-c_{\bar\nu}/\ptre)(
\nod{f}{L^2(\Omega)}+ \nod{g}{L^2_t(\partial\Omega)}).$$
\end{theorem}

\proof The proof of this result follows the same lines as the scalar case
treated in
\cite{Haddar-Joly-Nguyen1}. For the reader convenience, we
hereafter give the basic ideas. 
\\[12pt]
We introduce a cut-off function $\phi_{{\bar\nu}} \in C^\infty(\Omega)$ such that
$$
\phi_{{\bar\nu}}(x) = 0  \; \mbox{ in } \; \Oext\; ,\; \quad \phi_{{\bar\nu}}(x) = \beta^{\bar\nu} \; \mbox{ in }\; \Oint^{\bar\nu},
$$
where the constant $\beta^{\bar\nu} > 0$ is chosen such that
\begin{equation} \label{condition_phi}
\|\nabla \phi_{{\bar\nu}}\|_\infty  < \frac{1}{4}.
\end{equation}
We set $\upE = \exp\big(\phi_{{\bar\nu}}(x)/\ptre\big) \eg^\ptre$. Straightforward
calculations show that
$$
\rt \, \eg^\ptre = \exp\big(-\phi_{{\bar\nu}}(x)/\ptre\big) \{ \rt\, \upE -
\frac{1}{\ptre} \nabla \phi_{\bar\nu} \times \upE \}.
$$
Hence,
\begin{eqnarray*}
\rt\, \rt\, \eg^\ptre &=& \exp\big(-\phi_{{\bar\nu}}(x)/\ptre\big) \; \left\{ \rt\, \rt\, \upE  -
\frac{1}{\ptre} \nabla \phi_{\bar\nu} \times \rt\,
\upE \right.\\ 
&&  \left. \quad - \frac{1}{\ptre} \; \rt\, (\nabla
\phi_{\bar\nu} \times \upE) + \frac{1}{\ptre^2}
 \nabla \phi_{\bar\nu} \times  (\nabla \phi_{\bar\nu} \times \upE) \right\}
\end{eqnarray*}
Therefore  $\upE$ satisfies, after multiplying the above
equality by $\exp(\phi_{\bar\nu} / \ptre)$
\begin{eqnarray} \label{system-vexp}
&& \rt\, \rt\, \upE  - \frac{1}{\ptre} \left(\nabla \phi_{\bar\nu} \times \rt
  \,\upE +  \rt\, (\nabla \phi_{\bar\nu} \times \upE) \right) \nonumber \\
&&  + \left ( (- \eps \omega^2  + i \omega \sigma^\ptre)\upE
 + \frac{1}{\ptre^2} \nabla
\phi_{\bar\nu} \times  (\nabla \phi_{\bar\nu} \times \upE) \right)= i \omega f  \;\; \mbox{ in
  } \; \Omega
\end{eqnarray}
Taking the $L^2(\Omega)^3$ scalar product of this equation (\ref{system-vexp}) by
$\upE$  and using Stokes formulas  yields,
(recall that $\upE = \eg^\ptre$ and $\phi_{\bar\nu} = 0$ in
$\Oext$)
\begin{equation} \label{identite-vexp-2}
\begin{array}{lll}
&& \nodt{\rt\, \upE}{L^2(\Oint)}^2  - \frac{1}{\ptre} \{
( \nabla \phi_{\bar\nu} \times \rt \upE, \upE)_{L^2(\Oint)} +  (\nabla \phi_{\bar\nu} \times
 \upE, \, \rt\, \upE)_{L^2(\Oint)} \}  \\[12pt]
&&  ( - \eps_r\omega^2 + \frac{i}{\ptre^2}) \nodt{\upE}{L^2(\Oint)}^2  + \frac{1}{\ptre^2} (\nabla
\phi_{\bar\nu} \times  (\nabla \phi_{\bar\nu} \times \upE),  \upE)_{L^2(\Oint)}^2 
 \\ [12pt]
&& \dsp =  i \omega \left\{  (f, \eg^\ptre)_{L^2(\Oext)}
-  (g, \eg^\ptre)_{L^2_t(\partial \Omega)} \right\} - \nodt{\rt\,
\eg^\ptre}{L^2(\Oext)}^2  \\ [12pt]
&& \quad +\omega^2 \; \nodt{\eg^\ptre}{L^2(\Oext)}^2  -i \omega
  \nodt{\eg^\ptre_T}{L^2_t(\partial \Omega)}^2
\end{array}
\end{equation}
Let us denote by  $L_\ptre$ the right hand side of the previous
equality. According to Theorem~\ref{exist-stab}, there
exists a constant $C$ independent
of $\ptre$ such that
$$
|L_\ptre| \le C \, (\nodt{f}{L^2(\Omega)}^2+\nodt{g}{L^2_t(\partial \Omega)}^2).
$$
On the other hand, thanks to inequality (\ref{condition_phi}), and making use
of the inequality 
$|b - a| \geq |b| - |a|$, we have the lower bound
$$
\begin{array}{l}
\dsp  \left|( - \eps_r\omega^2 + \frac{i}{\ptre^2}) \nodt{\upE}{L^2(\Oint)}^2  + \frac{1}{\ptre^2} (\nabla
\phi_{\bar\nu} \times  (\nabla \phi_{\bar\nu} \times \upE),  \upE)_{L^2(\Oint)}^2 
\right| 
\\[12pt]
\dsp \qquad\qquad \qquad\qquad \qquad  \ge \left(\frac{1 - \|\nabla \phi_{{\bar\nu}}\|_\infty^2}{\ptre^2} - \eps_r\omega^2 \right)\nodt{\upE}{L^2(\Oint)}^2 
\ge \frac{3}{4\ptre^2} \; \nodt{\upE}{L^2(\Oint)}^2,
\end{array} 
$$
for $\ptre$ is sufficiently small. One therefore deduces from
(\ref{identite-vexp-2}) 
$$
  \nodt{\rt\, \upE}{L^2(\Oint)}^2 + \frac{3}{4\ptre^2}\; \nodt{\upE}{L^2(\Oint)}^2
\le \frac{2}{\ptre}\; \|\nabla
\phi_{{\bar\nu}}\|_\infty\nodt{\upE}{L^2(\Oint)}\nodt{\rt\, \upE}{L^2(\Oint)} +
C \, \big(\nodt{f}{L^2(\Omega)}^2+\nodt{g}{L^2_t(\partial \Omega)}^2 \big). $$
Then, using again  (\ref{condition_phi}) and inequality $ 2|a||b| \le a^2 +
b^2$ one gets
  $$
  \nodt{\rt\, \upE}{L^2(\Oint)}^2 + \frac{3}{4\ptre^2}\nodt{\upE}{L^2(\Oint)}^2
\le 2
C \, \big(\nodt{f}{L^2(\Omega)}^2+\nodt{g}{L^2_t(\partial \Omega)}^2 \big). $$
The final estimate can now be easily deduced by noticing that $\eg^\ptre =
\exp\big(-\beta^{\bar \nu}/\ptre\big)\upE $ in $\Oint^{\bar\nu}$,  and by
using the Maxwell equations to get estimates on $\hg^\ptre$ from those on $\eg^\ptre$.  \proofend

% $$
% \begin{array}{lll}
% \nodt{\rt \upE}{L^2(\Oint)}^2 & - & \dsp \frac{2}{\ptre}\nodt{\nabla
%   \phi_{{\bar\nu}}}{L^\infty} \nodt{\rt
%   \upE}{L^2(\Oint)}\nodt{\upE}{L^2(\Oint)}\\[12pt]
% & + & \dsp  \frac{1}{2 \ptre^2} \;
% \nodt{\upE}{L^2(\Oint)}^2
%   \le C \, (\nodt{f}{L^2(\Omega)}^2+\nodt{g}{L^2_t(\partial \Omega)}^2).
% \end{array}
% $$
% Due to the inequality,
% $$\frac{2}{\ptre}\nodt{\nabla
%   \phi_{{\bar\nu}}}{L^\infty} \nodt{\rt \,\upE }{L^2(\Oint)}
% \nodt{\upE}{L^2(\Oint)}
% \leq \frac{1}{2} \nodt{\rt \,\upE\upE}{L^2(\Oint)}^2
% + \frac{2}{\ptre^2} \nodt{\nabla
%   \phi_{{\bar\nu}}}{L^\infty}^2 \nodt{\upE}{L^2(\Oint)}^2,
% $$
% and inequality (\ref{condition_phi}),
% \begin{equation} \label{identite-vexp-3}
% \frac{1}{2}\nodt{\rt \upE}{L^2(\Oint)}^2 +  \frac{1}{4
%   \ptre^2} \nodt{\upE}{L^2(\Oint)}^2
%   \le C \,(\nodt{f}{L^2(\Omega)}^2+\nodt{g}{L^2_t(\partial \Omega)}^2),
% \end{equation}
% for $\ptre$ small enough. To be continued !!! \proofend

\section{Statement of the main results}\label{main} 
We shall denote by $(\egeek,\hgeek)$, the approximate solutions in the
exterior domain $\Oext$, the presence of the integer $k$ meaning that these
fields will provide an approximation of order $O(\ptre^{k+1})$ of the exact 
exterior electromagnetic field $(\egee,\hgee)$, in a sense that will be made
precise by the error estimates (see Theorem~\ref{theoremerror}). 
They are obtained by solving  the standard Maxwell equations in the exterior domain $\Omega_e$ 
\begin{equation} \label{eqinter}
\left\{ \begin{array}{lcl} i \omega \egeek- \rt\, \hgeek   = f  \;\; &
\mbox{ in } \; \Oext,
\\[6pt]
i \omega \hgeek + \rt\, \egeek   = 0  \;\; & \mbox{ in } \; \Oext,
\\[6pt]
\dsp  \Eke_{e,\,T} - \hgeek \times n  = g\;\;\;& \mbox{ on } \;
\partial \Omega,  \end{array}\right.
\end{equation}
where $n$ denotes the normal to $\partial \Omega$ directed to the exterior of
$\Omega$, coupled with an appropriate GIBC on the interior boundary $\Gamma$ of the form
\begin{equation}\label{GIBCcond}
\egeek \times n + \omega \;{\cal D}^{\ptre, k}(\Hke_{e,\,T}) = 0,
\end{equation}
where $n$ denotes the normal to $\Gamma$ directed to the exterior of
$\Oext$, $\Hke_{e,\,T}$ is the tangential trace of $\Hke_{e}$, and where
${\cal D}^{\ptre, k}$ is an adequate local approximation of the {\it 
  H-to-E} map for the Maxwell equations inside $\Oint$, namely the operator:
$$
{\cal D}^{\ptre} : H^{-\frac{1}{2}}(\rotg,\Gamma) \quad \longrightarrow \quad
H^{-\frac{1}{2}}(\divg,\Gamma) 
$$
defined by
$$
{\cal D}^{\ptre} \varphi = - \frac{1}{\omega} E_i^{\ptre} \times n|_{\Gamma}
$$
where $\big(E_i^{\ptre}(\varphi),
H_i^{\ptre}(\varphi) \big)$ is the solution of the {\it interior}
  boundary value problem 
$$%\begin{equation} \label{interiorpb}
\left\{ \begin{array}{ll}
(i \eps_r\omega + \frac{1}{\omega \delta^2}) \egie(\varphi) - \rt\, \hgie(\varphi)  = 0, \;\; & \mbox{ in } \; \Oint,\\[12pt]
i \omega \hgie(\varphi) + \rt\, \egie(\varphi) = 0, \;\; & \mbox{ in } \; \Oint, \\[12pt]
H_{i,T}^{\ptre}(\varphi) = \varphi, \;\;  & \mbox{ on }\;\Gamma. 
\end{array} \right.
$$%\end{equation}
\subsection{The ``natural'' GIBCs for $k=0,1,2$.}  \label{Sec-GIBC1}
The approach that we shall use in Section~\ref{derivation}
for the formal derivation of the GIBCs leads to the following
expressions of ${\cal D}^{\ptre, k}$  (for $k=0, 1, 2, 3$),
\begin{equation}\label{Dk}
\left\{\begin{array}{ll}
{\cal D}^{\ptre, 0} =  0, \\[6pt]
{\cal D}^{\ptre, 1} = \ptre\sqrti  , \\[6pt]
{\cal D}^{\ptre, 2} = \ptre \sqrti  +  \ptre^2 ({\cal H}-{\cal C}) , 
\end{array} \right.
\end{equation}
 where, $\sqrti := \frac{\sqrt{2}}{2} + i \frac{\sqrt{2}}{2}$
denotes the complex square root of $i$ with positive real part, and 
${\cal C}$ and ${\cal H}$ are the  curvature and mean curvature tensors of
$\Gamma$ (we refer to Section~\ref{preliminary} for more details).
Note that the condition of order $0$ simply expresses the fact that the limit
exterior problem when $\delta$ goes to 0 corresponds to the perfectly
conducting boundary condition.
\subsection{The modified third order GIBC} \label{modif}
The same approach extended to $k=3$ would suggest to take:
\begin{equation} \label{GIBC3}
{\cal D}^{\ptre, 3} = {\cal D}_0^{\ptre, 3}
\end{equation}
where by definition (we refer to Section \ref{preliminary} for the
definition of the surface operators $\grag$, $\divg$, $\rotg$ and $\Rotg$)
\begin{equation}\label{D3_0}
{\cal D}_0^{\ptre, 3} : = \ptre \sqrti  +  \ptre^2 ({\cal H}-{\cal C}) 
 + \frac{\ptre^3}{2 \sqrti} \left({\cal C}^2 -  {\cal H}^2 +\eps_r \omega^2  + \grag \divg
 + \Rotg \rotg \right).
\end{equation}
However, we did not succeed in proving that such a choice was mathematically
sound due to the presence of the
second order surface operator $\grag \divg + \Rotg \rotg$. As a self-adjoint
operator in $L^2_t(\Gamma)$, this operator (more precisely the associated
quadratic form) has no fix
sign. This induces difficulties in the study of the forward problem via
variational techniques and, as a consequence, the well-posedness of the 
corresponding boundary value problem is not clear: this is a new difficulty
with respect to the scalar wave equation. \\[12pt]
This is why we propose hereafter another third order condition, that
(formally) gives the same order of accuracy as the one in (\ref{Dk}) but
admits good mathematical properties with respect to stability and error
estimates. The reader will easily notice that the proposed modifications are not the
only possible ones (see for instance Remarks \ref{remmodif1} and
\ref{remmodif2}), we exhibit only one particular choice. We shall hereafter
present the intuitive reasons that led us to introduce these
modifications, postponing the rational justification to the error analysis of Section
\ref{errorGIBC-section}. 
\\[12pt]
% The idea to get this new
% condition is to modify the operator ${\cal D}_0^{\ptre, 3}$ in order to
% circumvent the difficulties encountered when trying to prove
% existence and uniqueness of solutions when the former expression of ${\cal
%   D}^{\ptre, 3}$ is used.
% The basic ingredient for the approximation consists in using two times  the
% formal substitution 
% $$
% 1 - \ptre^2 P_\Gamt = (1+ \ptre^2 P_\Gamt)^{-1} + O(\ptre^3)
% $$
% where $P_\Gamt$ denotes an appropriate second order positive surface operator.
% \\[12pt]
% We explain below how we construct the modified operator ${\cal D}^{\ptre,
%   3}$. 
The first desirable (and probably necessary) property is the absorption
property: 
$$
{\cal R}e \int_{\Gamma} {\cal D}^{\ptre, 3} \varphi \cdot \bar{\varphi} \; d
\sigma \geq 0,
$$
for any smooth tangential vector field $\varphi$ on $\Gamma$. Such a property
is satisfied by the exact DtN operator and expresses the absorbing nature of
the conductive medium:
$$
{\cal R}e \int_{\Gamma} {\cal D}^{\ptre, 3} \varphi \cdot \bar{\varphi} \; d
\sigma = \frac{1}{\omega \delta^2} \; \int_{\Oint} |\egie(\varphi)|^2 \; dx.
$$
It will play an essential role in proving the {\it uniqueness}
of solutions. One can observe that this condition is satisfied by
${\cal D}^{\ptre, 1}$ and ${\cal D}^{\ptre, 2}$. For ${\cal D}^{\ptre, 3}$, we
see that 
$$
{\cal R}e \; {\cal D}_0^{\ptre, 3} = \ptre \frac{\sqrt{2}}{2} +  \ptre^2 ({\cal
  H}-{\cal C})  
 + \frac{\ptre^3}{2 \sqrt{2}} \left({\cal C}^2 -  {\cal H}^2 +\eps_r \omega^2\right)  +
   \frac{\ptre^3}{2 \sqrt{2}} \; \grag \divg 
 + \frac{\ptre^3}{2 \sqrt{2}} \;\Rotg \rotg.
$$
The problem comes from the operator $\grag \divg$ which is negative in the
$L^2$ sense. However, we can write formally
\begin{equation} \label{splitting}
\left| \begin{array}{lll}
\dsp \frac{{\ptre}\sqrt{2}}{2} + \frac{\ptre^3}{2 \sqrt{2}} \; \grag \divg
& = & \dsp \frac{{\ptre}}{2\sqrt{2}} + \frac{{\ptre}}{2\sqrt{2}} + \frac{\ptre^3}{2
  \sqrt{2}} \; \grag \divg \\[12pt]
& = & \dsp \frac{\ptre}{2\sqrt{2}} \; (1 - {\ptre^2} \;
\grag \divg)^{-1} + O(\ptre^5),
\end{array} \right.
\end{equation}
which suggests to define the real part of ${\cal D}^{\ptre, 3}$ as 
\begin{equation} \label{reelle}
\left| \begin{array}{lcl}
{\cal R}e \; {\cal D}^{\ptre, 3} & = &  \dsp\frac{ \ptre}{2\sqrt{2}}   +  \ptre^2 ({\cal 
  H}-{\cal C})  + \frac{\ptre^3}{2 \sqrt{2}} \left({\cal C}^2 -  {\cal H}^2
  +\eps_r \omega^2\right)
 + \frac{\ptre^3}{2 \sqrt{2}} \;\Rotg \rotg
\\[12pt]
&+ & \dsp\frac{\ptre}{2\sqrt{2}} \; (1 - {\ptre^2} \; \grag \divg)^{-1}
\end{array} \right.
\end{equation}
\begin{remark} \label{remmodif1}
The approximation process (\ref{splitting}) is analogous to the process used
in the construction of absorbing boundary conditions for the wave equations,
see \cite{E-M-77, C-92} for instance, where the Pad\'e approximations are preferred to Taylor
approximations in order to enforce the stability of the resulting approximate problem.\\[12pt]
In (\ref{splitting}), the splitting $\dsp\frac{ \ptre\sqrt{2}}{2} = \dsp\frac{
  \ptre}{2\sqrt{2}} + \dsp\frac{ \ptre}{2\sqrt{2}}$ is somewhat
  arbitrary and could be changed into 
$$\dsp\frac{ \ptre\sqrt{2}}{2} = (1 - \alpha) \; \frac{ \ptre\sqrt{2}}{2} + 
\alpha \; \frac{ \ptre\sqrt{2}}{2}$$
for any $\alpha \in \; ]0,1[$. Our choice corresponds to $\alpha = 1/2$.
\end{remark}
The second modification was guided by the {\it existence} proof for the
boundary value problem associated to the boundary condition
(\ref{GIBCcond}). We realized that it was useful that the imaginary
part of  ${\cal D}^{\ptre, 3}$ satisfies a ``Garding type'' inequality, namely
that the principal part of this operator be positive in the $L^2$
sense.  This property is not satisfied by the imaginary part of ${\cal
  D}_0^{\ptre,3}$:~\\[4pt]  
$$
{\cal I}m \; {\cal D}_0^{\ptre, 3} = \ptre \frac{\sqrt{2}}{2} -   
\frac{\ptre^3}{2 \sqrt{2}} \left({\cal C}^2 -  {\cal H}^2 +\eps_r
  \omega^2\right)  - \frac{\ptre^3}{2 \sqrt{2}} \grag \divg
 - \frac{\ptre^3}{2 \sqrt{2}} \; \Rotg \rotg.
$$
This time, the problem is due to the negative operator $- \Rotg
\rotg$. The same manipulation as for the real part of ${\cal D}^{\ptre, 3}$
% suggests to 
% $$
% \frac{{\ptre}}{2\sqrt{2}} - \frac{\ptre^3}{2 \sqrt{2}} \;\Rotg
% \rotg =  
%  \frac{\ptre}{2\sqrt{2}} \; (1 + {\ptre^2} \; \Rotg
% \rotg)^{-1} +
% O(\ptre^5),$$
suggests to define the imaginary part of ${\cal D}_r^{\ptre, 3}$ as 
\begin{equation} \label{imagi}
\left|\begin{array}{lcl}
{\cal I}m \; {\cal D}_r^{\ptre, 3} &=& \dsp\frac{ \ptre}{2\sqrt{2}} +   
\frac{\ptre^3}{2 \sqrt{2}} \left({\cal C}^2 -  {\cal H}^2 +\eps_r
  \omega^2\right)  - \frac{\ptre^3}{2 \sqrt{2}} \grag \divg
\\[12pt]
& +&  \dsp \frac{\ptre}{2\sqrt{2}} 
\; (1 + {\ptre^2} \; \Rotg \rotg)^{-1}.
\end{array} \right.
\end{equation}
Modifications (\ref{reelle}) and (\ref{imagi}) lead us to introduce the
operator  
\begin{equation} \label{D3modified0}
\left|\begin{array}{lcl}
\dsp \tilde {\cal D}^{\ptre, 3}& =& \dsp \ptre \frac{\sqrt{i}}{2} + \ptre^2 ({\cal H}-{\cal C}) 
 + \frac{\ptre^3}{2 \sqrti} \left({\cal C}^2 -  {\cal H}^2 +\eps_r
 \omega^2\right) 
\\[12pt] 
&&\dsp + \; \frac{\sqrt{2}}{4} \; \ptre \; \big(\; (1   - {\ptre^2} \; \grag
 \divg)^{-1} + {\ptre^2} \; \Rotg \rotg
 \big) 
\\[12pt] 
&&\dsp+ \; i \; \frac{\sqrt{2}}{4}\;\ptre \; \big( \; (1   + {\ptre^2} \; 
 \Rotg \rotg)^{-1} -  {\ptre^2}\; \grag \divg\big).
\end{array} \right.
\end{equation}
~\\[10pt]
which formally satisfies $\tilde {\cal D}^{\ptre, 3} = {\cal D}^{\ptre, 3}_0 +
O(\ptre^5)$. \\[12pt]
It turns out that even if this condition is suitable for
variational study of existence and uniqueness of the resulting boundary value
problem, it did not enable us to have a direct proof of optimal error
estimates (although we think it can be achieved by constructing the full
asymptotic expansion associated with the associated boundary value problem). We realized that the difficulties encountered in the analysis are
related to the fact that the operator $\tilde {\cal D}^{\ptre, 3}$ is a
pseudo-differential operator of order 2, whiles the exact impedance operator
which maps continuously $H^{-1/2}(\rotg, \Gamma)$ into $H^{-1/2}(\divg,
\Gamma)$ is more something between an operator of order $-1$
and an operator of order $1$. This gave us the idea to force our approximate
operator to be of order $0$ by applying a {\it regularization} process (the
Yosida regularization) to the
operators $\Rotg \rotg$ and $\grag \divg$
\begin{equation} \label{regYos}
\left| \begin{array}{ll}
\Rotg \rotg \simeq \Rotg \rotg \, (1   + {\ptre^2} \; 
 \Rotg \rotg)^{-1}  &\quad \mbox{in $O(\ptre^2)$}, \\[12pt]
\grag \divg \simeq\grag \divg \, (1   - {\ptre^2} \; \grag
 \divg)^{-1}\big)&\quad \mbox{in $O(\ptre^2)$}.
\end{array} \right.
\end{equation}
Such an approximation is consistent with the $O(\ptre^5)$ accuracy provided by
$\tilde {\cal D}^{\ptre, 3}$ since $\Rotg \rotg$ and $\grag \divg$ are
multiplied by $\ptre^3$. Moreover, it does not affect the good sign properties
of the real and imaginary parts of the operator since we ``divide'' by positive
operators. Therefore, we propose for the third order condition the following expression
\begin{equation} \label{D3modified}
\left|\begin{array}{lcl}
\dsp  {\cal D}^{\ptre, 3}& :=& \dsp \ptre \frac{\sqrt{i}}{2} + \ptre^2 ({\cal H}-{\cal C}) 
 + \frac{\ptre^3}{2 \sqrti} \left({\cal C}^2 -  {\cal H}^2 +\eps_r
 \omega^2\right) 
\\[12pt] 
&&\dsp + \; \frac{\sqrt{2}}{4} \; \ptre \; \big(\; (1   - {\ptre^2} \; \grag
 \divg)^{-1} + {\ptre^2} \; \Rotg \rotg \, (1   + {\ptre^2} \; 
 \Rotg \rotg)^{-1}
 \big) 
\\[12pt] 
&&\dsp+ \; i \; \frac{\sqrt{2}}{4}\;\ptre \; \big( \; (1   + {\ptre^2} \; 
 \Rotg \rotg)^{-1} -  {\ptre^2}\; \grag \divg\, (1   - {\ptre^2} \; \grag
 \divg)^{-1}\big).
\end{array} \right.
\end{equation}
\begin{remark}\label{remmodif2}
The regularization process (\ref{regYos}) can also be seen as an analogous to
the {\it stabilization} process used in numerical methods such as stabilized
finite elements or discontinuous Galerkin methods, in order to ensure optimal
error estimates.
Here again, in (\ref{regYos}), we chosed arbitrarily equal to 1 the
regularization constant in the term in factor of $\ptre^2$. We could have
chosen for instance
$$
\Rotg \rotg \simeq \Rotg \rotg \, (1   + \beta \; {\ptre^2} \; 
 \Rotg \rotg)^{-1}
$$
where $\beta$ is a positive constant. However, the interest
of choosing $\beta = 1$, is that we make appear the same operators $(1   +
{\ptre^2} \; \Rotg \rotg)^{-1}$ and $(1   - {\ptre^2} \; \grag
 \divg)^{-1})$ that are already present in the expression of $\tilde {\cal
   D}^{\ptre, 3}$), which is of interest from the numerical point of view.
\end{remark}
% which also formally satisfies ${\cal D}^{\ptre, 3} = {\cal D}^{\ptre, 3}_0 +
% O(\ptre^5)$. We remark that 
% ${\cal D}^{\ptre, 3}$ is formally a  ``pseudo differential'' operator of order
% $0$ while  $\tilde {\cal D}^{\ptre, 3}$ is of order $2$.
The operator happens to have the good {\it consistency}, {\it coercivity} and
{continuity} properties to lead to optimal error estimates. More precisely:
\begin{itemize}
\item One can check that 
\begin{equation} \label{errorD}
{\cal D}^{\ptre, 3} = {\cal D}^{\ptre, 3}_0 + \ptre^5 \;{\cal R}^{\ptre, 3}
\end{equation}
where  the operator ${\cal R}^{\ptre, 3}$, given by 
$$
{\cal R}^{\ptre, 3} = \frac{\sqrt{2}}{4} (1+i) \; \left[
\big(1   - {\ptre^2} \; \grag \divg\big)^{-1}\; (\grag \divg)^2
+ \big(1   + {\ptre^2} \; \Rotg \rotg\big)^{-1}\; (\Rotg \rotg)^2
\right],
$$
maps continuously $H^{s+4}_t(\Gamma)$ in $H^{s}_t(\Gamma)$ and
satisfies the uniform bound
\begin{equation} \label{boundR}
\|{\cal R}^{\ptre, 3}\|_{{\cal L}\big(H^{s+4}_t(\Gamma);H^{s}_t(\Gamma)\big)}
\; \leq \; 1. 
\end{equation}
\item One can prove (see Lemma \ref{lembuy}) that ${\cal D}^{\ptre, 3}$ has
  the following fundamental properties (obviously satisfied by ${\cal
    D}^{\ptre, 1}$ and  
${\cal D}^{\ptre, 2}$)
$$ 
\forall \varphi \in L^2_t(\Gamma), \quad 
\|{\cal D}^{\ptre, k} \varphi\|_\Gamma \le C_1 \; \ptre \; \|\varphi\|_\Gamma,
\quad
\Ree ({\cal D}^{\ptre, k} \varphi, \varphi)_\Gamma \ge C_2 \; \ptre \; 
\|\varphi\|_\Gamma^2.
$$
with $C_1$ and $C_2$ strictly positive constants.
These appear to be sufficient properties to transform the consistency
properties ${\cal D}^{\ptre, 3}$ into optimal error estimates (see the
proof of Lemma~\ref{lemerreor}).
\end{itemize}
% \begin{remark} \label{rigoureux}
% To be more precise concerning the formal consistence of the derived
% expressions, one can check that 
% $$
% \tilde{\cal D}^{\ptre, 3} = {\cal D}^{\ptre, 3}_0 + \ptre^5 \; \tilde{\cal R}^{\ptre, 3}
% $$
% where  the operator $\tilde{\cal R}^{\ptre, 3}$ is given by 
% $$
% \tilde {\cal R}^{\ptre, 3} = \frac{\sqrt{2}}{4} \; \left[
% \big(1   - {\ptre^2} \; \grag \divg\big)^{-1}\; (\grag \divg)^2
% + i \big(1   + {\ptre^2} \; \Rotg \rotg\big)^{-1}\; (\Rotg \rotg)^2
% \right],
% $$
% maps continuously $H^{s+4}_t(\Gamma)$ in $H^{s}_t(\Gamma)$ and
% satisfies the uniform bound
% \begin{equation} \label{boundR0}
% \|\tilde{\cal R}^{\ptre, 3}\|_{{\cal L}\big(H^{s+4}_t(\Gamma);H^{s}_t(\Gamma)\big)}
% \; \leq \; 1/2. 
% \end{equation}
% On the other hand, 
% $$
% {\cal D}^{\ptre, 3} = {\cal D}^{\ptre, 3}_0 + \ptre^5 \;{\cal R}^{\ptre, 3}
% $$
% where  the operator ${\cal R}^{\ptre, 3}$ is given by 
% $$
% {\cal R}^{\ptre, 3} = \frac{\sqrt{2}}{4} (1+i) \; \left[
% \big(1   - {\ptre^2} \; \grag \divg\big)^{-1}\; (\grag \divg)^2
% + \big(1   + {\ptre^2} \; \Rotg \rotg\big)^{-1}\; (\Rotg \rotg)^2
% \right],
% $$
% maps continuously $H^{s+4}_t(\Gamma)$ in $H^{s}_t(\Gamma)$ and
% satisfies the uniform bound
% \begin{equation} \label{boundR}
% \|{\cal R}^{\ptre, 3}\|_{{\cal L}\big(H^{s+4}_t(\Gamma);H^{s}_t(\Gamma)\big)}
% \; \leq \; 1. 
% \end{equation}
% \end{remark}
% Our analysis and error estimate apply to the third order
% condition where the operator ${\cal D}^{\ptre, 3}$ is replaced by ${\cal
%   D}^{\ptre, 3}_m$.

\subsection{Existence, uniqueness and error estimates} \label{Sec_Error}
The natural functional spaces for the solutions of the approximate problems
vary according to the regularity of their traces on $\Gamma$. We shall
distinguish the case $k=0$ for which we set 
$$
{\cal V}_H^0 = \{ H \in H(\rt,\Oext) \; ; \; 
(H \times n)|_{\partial \Omega}  \in L^2_t(\partial \Omega) \}, \;\;
 {\cal V}_E^0 = \{ E \in {\cal V}_H^0 \;;\; (\eg \times n)|_\Gamma=0\},
$$
from the case $k=1$, $2$ or $3$ for which we set
$$
{\cal V}_H^k = {\cal V}_E^k = \hcurlt{\Oext}
$$
(see (\ref{defhcurlt}) for the definition of $\hcurlt{\Oext}$).
 Then we have the following central
theorem, that uses and combines the partial results of Sections 4-6.
\begin{theorem} \label{theoremerror}
For $k=0$, 1, 2 or 3,  there exists $\ptre_k$ such that for $\ptre
\le \ptre_k$, the boundary value problem ((\ref{eqinter}),
(\ref{GIBCcond})) has a unique solution $(\egeek, \hgeek) \in
{\cal V}_E^k \times {\cal V}_H^k$. Moreover, if $B$ is a bounded 
domain containing $\Omega$, then there
exists a constant $C_k$, independent of $\ptre$, such that
\begin{equation*}% \label{error}
\| E^\ptre_e - \egeek \|_{H(\rt, \Oext)} \le C_k \; \ptre^{k + 1}.
\end{equation*}
\end{theorem}
\begin{remark}
For $k=0, 1$, the above theorem holds for all $\ptre$.
\end{remark}

\section{Formal derivation of the GIBCs} \label{derivation}
\subsection{Preliminary material}
\label{preliminary}
We recall in this section some well known facts
about differential geometry and differential operators. 
\paragraph{Local coordinates.} Let $n$ be the normal field
  defined on $\Gamma$ and directed to the interior of $\Oint$. For a
  sufficiently small 
  positive constant ${\bar\nu}$ (see condition (\ref{condnubar}) below)  we define
$$\Oint^{\bar\nu} = \{ x \in \Oint \; ; \; \mathrm{dist}(x, \partial \Oint) < {\bar\nu}\}.$$
To any $x \in \Oint^{\bar\nu}$ we uniquely associate the local parametric
coordinates 
$(\xgg, \nu) \in \Gamma \times (0, {\bar\nu})$  through
\begin{equation}
\label{param}
x = \xgg + \nu \, n, \;\;\; x \in \Oint^{\bar\nu}.
\end{equation}
\paragraph{Tangential (or surface) differential operators.} 
In what follows we deal with various fields defined on $\Gamma$: scalar fields
$\varphi$ (with values in $\C$), vector fields $V$ (with values in $\C^3$) and
matrix (or tensor) fields ${\bf A}$ (with values in ${\cal L}(\C^3$)). By definition:
\begin{itemize}
\item A vector field $V$ is tangential if and only if $V\cdot n = 0$ (as a scalar
  field along $\Gamma$).
\item  A matrix field ${\bf A}$ is tangential if and only if ${\bf A} \; n = 0$
  (as a vector field along $\Gamma$).
\end{itemize}
For simplicity, we assume that these fields have at least $C^1$
regularity, but this can be removed by interpreting the derivatives in the sense of
distributions.\\[12pt] 
We recall that the surface gradient operator $\grag$ is defined by:
$$
 \quad \grag \varphi(x_{\Gamma}) = \nabla
\hat{\varphi} (x_{\Gamma}), \quad \forall \varphi : \Gamma \rightarrow \R,
$$
where $\hat{\varphi}$ is the 3-D vector field defined locally in
$\Oint^{\bar\nu}$ by $\hat{\varphi} (x_{\Gamma} + \nu \, n) =
\varphi(x_{\Gamma})$.  Note that $\grag \varphi$ is 
a tangential vector field. We can define in the same way
the surface gradient of a vector field as a tangential matrix field whose
columns are the surface gradients of each component of the vector field.\\[12pt]
We denote by $-\divg$ the $L^2(\Gamma)$- adjoint
of $\grag$ : $-\divg$ maps a tangential vector field into a scalar field.
More generally, if ${\bf A}(\xgg)$ is a tangential matrix field
on $\Gamma$, we  define the operator  ${\bf A}\grag$  for a scalar field
$\varphi(\xgg)$ by $$({\bf A}\grag) u : = {\bf A} (\grag u).$$
In the same way, we define the operator $({\bf A}\grag) \cdot$ acting on a  
tangential vector field $V(\xgg)$ as:
$$
({\bf A}\grag) \cdot V := \sum_{i=1}^3 ({\bf A}\grag V_i )_i,
$$
where the subscript $i$ denotes the $i^{th}$ component of a vector in the
canonical basis of $\R^3$. 
% We remark that for ${\cal R} = \Ig$ we respectively retrieve
% the definition of the surface gradient of $u$ and surface divergence of $V$.
\\[12pt]
We then define the surface curl
of a tangential vector filed $V(\xgg)$  
and the surface vector curl of a scalar function $\varphi(\xgg)$ as
$$
\rotg V := \divg (V \times n) \;\; \mbox{ and } \;\; \Rotg \varphi : = (\grag \varphi)
\times n.   
$$ 
Of course, the various operators $\grag$, $\divg$, ${\cal R}\grag$, $\rotg$
and $\Rotg$ applies in principle to functions defined on
$\Gamma$. However, they can obviously be understood as (partial) differential
operators acting on fields defined in the (3-D) domain 
$\Gamma \times (0,{\bar\nu})$.
For instance, if $\phi \in C^1(\Gamma \times (0,{\bar\nu}))$, we define $\grag
\phi
\in C^0(\Gamma \times (0,{\bar\nu}))^3$ as follows (with obvious notation):
\begin{equation*} %\label{gradgamma3D}
\forall \; \nu \in (0,{\bar\nu}), \quad \Big[\grag \phi\Big](\cdot, \nu)] :=
\grag \Big[\phi(\cdot, \nu)\Big].
\end{equation*}
We apply similar rules to $\divg$, ${\cal R}\grag$, $\rotg$
and $\Rotg$. The extension of these definition in the sense of
distributions is also elementary (as soon as $\Gamma$ is $C^{\infty}$).
\paragraph{Geometrical tools.} In what follows, and for the sake of the notation
conciseness, we shall most of time not explicitly indicate the
dependence on $\xgg$ of the functions, except when we feel it necessary. We shall be more precise in mentioning the possible dependence
with respect to the normal coordinate $\nu$.\\[12pt]
A particularly fundamental tensor field  is the
curvature tensor ${\cal C}$, defined by  ${\cal C}:= \nabla_\Gamt n$.  We
recall that ${\cal C}$ is  
symmetric and ${\cal C} \, n = 0$. We denote $c_1$, $c_2$ the other two
eigenvalues of ${\cal C}$ (namely the {\it principal curvatures}) associated
with tangential eigenvectors $\tau_1$, $\tau_2$ ($\tau_1 \cdot n = \tau_2
\cdot n = 0$). We also introduce
\begin{equation} \label{courbures}
g := c_1 c_2 \quad \mbox{ and } \quad {h} := \frac{1}{2}(c_1 + c_2)
\end{equation}
which are respectively the {\it Gaussian} and {\it mean curvatures} of
$\Gamma$, and also introduce the
associated matrix fields:
\begin{equation} \label{courbures2}
{\cal H} = h \, \Ig \;\; \mbox{ and } \;\; {\cal G} = g \, \Ig,
\end{equation}
where $\Ig(\xgg)$ denotes the projection operator on the tangent
plane to $\Gamma$ at $\xgg$.\\[12pt]
Let us introduce (this is the Jacobian of the transformation
$(\xgg, \nu) \rightarrow x$ - see (\ref{param}))
\begin{equation*}% \label{defJ}
J(\nu) \; \big(= J(\nu, \xgg) \big) \;:= \mathrm{det} ({ I + \nu \, {\cal C}})
= 1 + 2 \nu h  + \nu^2 g,
\end{equation*}
and we choose $\bar{\nu}$ sufficiently small in such a way that 
\begin{equation} \label{condnubar}
\forall \; \nu < \bar{\nu}, \quad \forall \; \xgg \in \Gamma, \quad J(\nu,
\xgg)) = 1 + 2 \nu h(\xgg) + \nu^2 g(\xgg) > 0.
\end{equation}
Thus, for each $\nu < \bar{\nu}$, there exists a tangential
matrix field $\xgg \rightarrow {\cal R}_\nu(\xgg)$ such that
$$(I + \nu \, {\cal C}(\xgg)) \; {\cal R}_\nu(\xgg)  = \Ig(\xgg). $$
More precisely, there exists a tangential matrix field on $\Gamma$, ${\cal
  M}(\xgg)$ , such that:
\begin{equation*}% \label{defM}
\Ig + \nu {\cal M} := J(\nu) \; {\cal R}_\nu, \quad \forall \; \xgg \in \Gamma,
\quad \forall \; \nu < \bar{\nu}.
\end{equation*}
One easily sees (using for instance the eigenbasis ($\tau_1$,
$\tau_2$, n)  of ${\cal C}$) that 
\begin{equation*}% \label{propM}
{\cal M} = 2 {\cal H} - {\cal C} \quad \mbox{ and } \quad {\cal M} \, {\cal C} = {\cal G}. 
\end{equation*}
\paragraph{The curl operator in local coordinates.} The basic step of our
forthcoming calculations will be to rewrite the Maxwell equations in the domain $
\Oint^{\bar\nu}$,by  using the local
coordinates.
For this, we need the expression of the curl operator in the variables $(\xgg,
\nu)$.  It is shown in \cite{Haddar-Joly1} that the curl of a 3-D vector field $V:
\Omega_i^{\bar\nu} \rightarrow \R^3$  is given in parametric coordinates by :
\begin{equation*}% \label{identite_1}
\rt\, V = \left[ ({\cal R}_\nu \grag) \cdot (\widehat V \times n)\right] n + \left[
  {\cal R}_\nu \grag (\widehat V \cdot n) \right] \times n -  ({\cal R}_\nu {\cal C}
  \widehat V) \times n - \partial_\nu (\widehat V \times n),
\end{equation*}
where  $V$ and $\widehat V$ (defined on $\Gamma \times (0, \bar \nu)$) are
related by 
$$
\widehat V(\xgg, \nu) = V(\xgg+\nu \,n).
$$
This formula can be written in a more convenient form (for us), after
multiplication by $J(\nu)$:
\begin{equation*} %\label{curlpara}
\left| \begin{array}{lll}
J(\nu) \; \rt\, V &=&  \dsp \left[ \big( (I + \nu {\cal M})
    \grag\big) \cdot (\widehat V \times   n)\right] n  + 
\left[ (I + \nu {\cal M})  \grag (\widehat V \cdot n)
\right] \times n 
\\[12pt]
\dsp &- & \dsp   \left[({\cal C} + \nu {\cal G}) \widehat V
\right] \times n - J(\nu) \;\partial_\nu (\widehat V \times n),
\end{array} \right.
\end{equation*}
or, in an equivalent form, 
\begin{equation} \label{identite_2}
J(\nu) \; \rt\, V = \left( C_{\Gamma} + \nu \; C_{\Gamma}^M \right) \widehat V
- J(\nu) \; \partial_\nu (\widehat  V \times n)
\end{equation}
where we have introduced the notation
\begin{equation}\label{notationC}
\left\{ \begin{array}{l}
C_{\Gamma} \widehat V = \big( \rotg \widehat V \big)n + \Rotg(\widehat V \cdot
n) - {\cal C}\widehat V \times n\\[12pt]
C_{\Gamma}^M\widehat V = \big(\rotgm \widehat V \big)n + \Rotgm(\widehat V
\cdot n) - {\cal 
  G} \widehat V \times n\\[12pt]
\Rotgm u :=({\cal M} \grag u) \times n \quad \mbox{and} \quad \rotgm \widehat  V  = ({\cal M} \grag) \cdot (\widehat V \times n).
\end{array} \right.
\end{equation}
This  expression is convenient for the asymptotic matching procedure,
described hereafter,
because we made explicit the (polynomial)
dependence of the operators with respect to $\nu$.
\paragraph{Functional spaces on $\Gamma$ and trace spaces.}~\\[12pt]
We assume that the definition of $H^s(\Gamma)$ for any real $s$ is well known.  We shall denote~by
$$
( \cdot , \cdot)_{\Gamma} \quad \mbox{ and } \left< \cdot , \cdot\right>_{\Gamma},
$$
respectively the inner product in $L^2(\Gamma)^3$ and the duality bracket
${\cal D}'(\Gamma)^3 - {\cal D}(\Gamma)^3$.\\[12pt]
Next, we introduce some notation for spaces of tangent vector fields along
$\Gamma$. For any $s \leq 0$, we set:
\begin{equation*}% \label{defHst}
\left\{ \begin{array}{l}
H^s_t(\Gamma) = \{ V \in H^s(\Gamma)^3 \; /\; V \cdot n = 0 \mbox{
  on } \Gamma \} \qquad (H^0_t(\Gamma) = L^2_t(\Gamma) \\[12pt]
H^{-s}_t(\Gamma) = \{ V \in H^{-s}(\Gamma)^3 \; /\; \left< V, \,\varphi \,n
\right>_{\Gamma} = 0, \; \forall \; \varphi \in H^s(\Gamma) \} \quad (\equiv
  \big(H^s_t(\Gamma)\big)' 
\end{array} \right.
\end{equation*}
as well as 
\begin{equation*}% \label{defHstrotdiv}
\left\{ \begin{array}{l}
H^s(\divg, \Gamma) = \{ V \in H^s_t(\Gamma)^3 \; / \;\divg V \in H^s(\Gamma) \}
\\[12pt] 
H^s(\rotg, \Gamma) = \{ V \in H^s_t(\Gamma)^3 \; / \;\rotg V \in H^s(\Gamma) \}
\end{array} \right.
\end{equation*}
equipped with their natural graph norms (we notice that $H^0(\divg,\Gamma)$ and $H^0(\rotg,\Gamma)$ are often denoted
by respectively $H(\divg,\Gamma)$ and 
 $H(\rotg,\Gamma)$). 
Finally, we recall the well known trace theorems:
\begin{theorem}
The two trace mappings
\begin{equation*} %\label{traces}
\left\{ \begin{array}{llll}
u \in C^{\infty}(\Oext)^3 & \mapsto & u\times n|_{\Gamma} &  \\[12pt]
u \in C^{\infty}(\Oext)^3 & \mapsto & u_T = u - (u\cdot n) n \quad \big(\equiv
n \times (u \times n) \big) & 
\end{array} \right.
\end{equation*}
can be extended as continuous and surjective linear applications
from $H(\rt, \Oext)$ onto $H^{-
  \frac{1}{2}}(\divg, \Gamma)$ and  $H^{- \frac{1}{2}}(\rotg, \Gamma)$ respectively.
Moreover, $H^{-\frac{1}{2}}(\divg, \Gamma)$ is the dual of $H^{-
  \frac{1}{2}}(\rotg, \Gamma)$ and one has the Green's formula:
\begin{equation*}% \label{green_rot}
\int_{\Omega} \big( \rt \,u \cdot v - u \cdot \rt \, v \big) \; dx =
\left< u \times n, v_T \right>_{\Gamma} = - \left< v \times n, u_T
\right>_{\Gamma} 
\end{equation*}
\end{theorem}
$\forall \; (u,v) \in H(\rt, \Oext)^2$.
\subsection{The asymptotic ansatz.} \label{asymptotic_ansatz}
To formulate our ansatz, it is
useful to 
introduce a cutoff function $\chi \in C^\infty(\Oint)$ such that
$\chi = 1$ in $\Oint^{{\bar\nu}}$ and $\chi = 0 $ in $\Oint
\setminus \Oint^{2{\bar\nu}}$ for a sufficiently small $\bar\nu >0$. For this ansatz we are not interested in the part
of the solution inside the support of 
$(1 - \chi)$, since we already know that the norm of the solution in this part
exponentially decay to $0$
as $\ptre$ goes to $0$ (this is Theorem \ref{apriori-estimate-2}). For the
remaining part of the solution, we postulate the following
expansions:
\begin{equation} \label{ansatz_1}
\left| \begin{array}{ll}
\egee(x) = \ege^0(x) + \ptre \; \ege^1(x) + \ptre^2 \; \ege^2(x) + \cdots \;\;\; \mbox{
for } x \in \Oext,\\[6pt]
\hgee(x) = \hge^0(x) + \ptre \; \hge^1(x) + \; \ptre^2 \hge^2(x) + \cdots \;\;\; \mbox{
for } x \in \Oext,\\
\end{array} \right.
\end{equation}
where $\ege^\ell$, $\hge^\ell$, $\ell=0,\, 1,\, \cdots$ are functions defined on $\Oext$ and
\begin{equation} \label{ansatz_2}
\left| \begin{array}{ll}
\chi(x) \egie(x) = \egi^0(\xgg, \nu/\ptre) + \ptre \; \egi^1(\xgg, \nu/\ptre) +
  \ptre^2 \; \egi^2(\xgg, \nu/\ptre) + \cdots   \;\;\; \mbox{
  for } x \in \Oint^{\bar\nu} \\[6pt]
\chi(x) \hgie(x) = \hgi^0(\xgg, \nu/\ptre) + \ptre \; \hgi^1(\xgg, \nu/\ptre) +
  \ptre^2 \; \hgi^2(\xgg, \nu/\ptre) + \cdots   \;\;\; \mbox{
  for } x \in \Oint^{\bar\nu}
\end{array} \right.
\end{equation}
where $x$, $\xgg$ and $\nu$ are as in (\ref{param}) and where
$\egi^\ell (\xgg, \eta)$, $\hgi^\ell (\xgg, \eta)$ : $\Gamma \times
\R^+ \mapsto {\mathbb C}$ satisfy
\begin{equation} \label{cond1_EHik}
\left| \begin{array}{ll}
\dsp \mbox{For a.e. } \xgg \in \Gamma, \quad \int_0^{+\infty} |\egi^\ell(\xgg, \eta)|^2 \; d \eta < + \infty = 0, \\[12pt]
\dsp \mbox{For a.e. } \xgg \in \Gamma,
\quad \int_0^{+\infty} |\hgi^\ell(\xgg, \eta)|^2 \; d \eta < + \infty.
\end{array} \right.
\end{equation}
\begin{remark}
The condition (\ref{cond1_EHik}) will imply that   $\egi^\ell$ and $\hgi^\ell$ are
exponentially decreasing with respect~to~$\eta$. 
\begin{equation*}% \label{cond1_EHik-2}
\left| \begin{array}{ll}
\dsp \mbox{For a.e. } \xgg \in \Gamma, \quad  \lim_{\eta \rightarrow \infty}
\egi^\ell(\xgg, \eta) = 0, \quad \mbox{(exponentially fast)} \\[6pt]
\dsp \mbox{For a.e. } \xgg \in \Gamma, \quad \lim_{\eta \rightarrow \infty} \hgi^\ell(\xgg, \eta) = 0, \quad \mbox{(exponentially fast)}
\end{array} \right.
\end{equation*}
which is coherent with the existence of a boundary layer suggested by
Theorem~\ref{apriori-estimate-2}. 
\end{remark}
\begin{remark}
Expansion (\ref{ansatz_2}) makes sense since the local
coordinates $(\xgg, \nu)$ can be used inside the support of $\chi$.
\end{remark}
In the next step, we shall identify the set of equations satisfied by
$(\ege^\ell, \hge^\ell)$ and $(\egi^\ell, \hgi^\ell)$,  $\ell \ge 0$ by
writing, formally, that we want to solve the transmission problem
(\ref{trans_eq_1})-(\ref{trans_eq_2}). \\[12pt]
In the sequel, it is useful to introduce the notation
\begin{equation} \label{tuik}
\left| \begin{array}{ll}
\widetilde E_i^\ptre(\xgg, \eta) := \egi^0(\xgg, \eta) + \ptre\; 
\egi^1(\xgg, \eta) + \ptre^2 \; \egi^2(\xgg, \eta) + \cdots \;\;\; (\xgg,
\eta) \in \Gamma \times \R^+,\\[6pt]
\widetilde H_i^\ptre(\xgg, \eta) := \hgi^0(\xgg, \eta) + \ptre\; 
\hgi^1(\xgg, \eta) + \ptre^2 \; \hgi^2(\xgg, \eta) + \cdots \;\;\; (\xgg,
\eta) \in \Gamma \times \R^+.
\end{array} \right.
\end{equation}
so that ansatz (\ref{ansatz_2}) has to be understood as
\begin{equation} \label{ansatz_3}
\left| \begin{array}{ll} \chi(x) \egie(x) =  \widetilde
E_i^\ptre(\xgg, \nu/\ptre) + O(\ptre^{\infty}) \;\;\; \mbox{ for }
x \in \Oint^{\bar\nu}, \\ [6pt] \chi(x) \hgie(x) =  \widetilde
H_i^\ptre(\xgg, \nu/\ptre) + O(\ptre^{\infty}) \;\;\; \mbox{ for }
x \in \Oint^{\bar\nu}.
\end{array} \right.
\end{equation}
\subsection{The equations for the exterior fields.}
This is the easy part of the job. The equations are directly derived from
(\ref{trans_eq_1})  and we obtain that
$(\egek, \hgek)$ satisfy
\begin{equation} \label{eqinterk}
\left\{ \begin{array}{ll}
 i \omega \egek - \rt\, \hgek = f_k, \;\; & \mbox{ in } \; \Oext,\\[6pt]
 i \omega \hgek + \rt\, \egek = 0, \;\; & \mbox{ in } \; \Oext, \\[6pt]
  {\egek}_{\,|T} - \hgek \times n = g_k, \;\; & \mbox{ on } \;\partial \Omega
\end{array}\right.
\end{equation}
where we have set $f_0 = f$, $g_0 = g$ and $f_k = 0$, $g_k =0$, for $ k \ge
1$, (\ref{eqinterk}) being complemented with the interface condition 
\begin{equation} \label{interf-0}
\egek|_\Gamma(\xgg) \times n = \egik(\xgg, 0) \times n , \;\; \mbox{ for } \; \xgg \in \Gamma,
\end{equation}
which completely defines $(\egek, \hgek)$ if $\egik(\xgg, 0) \times n$ is
known. % During the formal construction of the GIBCs we shall consider $(\egek,
%\hgek)$ as being a given data.

\newcommand{\hgiet}{\widetilde H_i^\ptre}
\newcommand{\egiet}{\widetilde E_i^\ptre}
\newcommand{\Opun}{{\bf C}_{\Gamt}\,}
\newcommand{\Opde}{{\bf C}_{\Gamt}^M\,}

%\begin{remark}
%In dimension two, if $s$ denotes the curvilinear abscissa on $\Gamma$ and $c$
%the curvature then identity (\ref{identite_1}) becomes
%\begin{equation}
%{\bar\nu} = \left(\frac{1}{1+\nu \, c}\right) \, \partial_s\, \left(
%  \frac{1}{1+\nu \, c} \right) \, \partial_s +
%\left(\frac{1}{1+\nu \, c}\right) \,\partial_\nu \,(1+\nu \, c) \,\partial_\nu
%\end{equation}
%\end{remark}

\subsection{The equations for the interior fields.}
\label{Asymptotic formal matching}
\noindent As indicated above, we need to compute the interior fields
$\egik$. The 
principle consists into expressing this field in terms of the tangential boundary
values of $(\hge^\ell)$, $\ell \le k$ by solving the interior equations. More precisely, we now substitute the expansion (\ref{tuik},
\ref{ansatz_3}) into  the system (\ref{trans_eq_2}) and assume that
the quantity: 
\begin{equation*}% \label{Hbord}
\hgee \times n = \sum_{k=0}^{+ \infty} \ptre^k \; \hge^k \times n 
\end{equation*}
is known on $\Gamma$.
We need of course to rewrite the equations of (\ref{trans_eq_2}) in the local 
``scaled'' coordinates 
$$(\xgg, \eta = \nu/\delta).$$
Using formula (\ref{identite_2}) with $\nu = \ptre \eta$ we obtain:
\begin{equation} \label{asymtotic_equation-1}
\hspace*{-0.6cm}\left\{ \begin{array}{ll}
\dsp J(\ptre \eta) \, \Big(i \eps_r\omega + \frac{1}{\omega {\ptre}^2}\Big)
\egiet - \big(\Opun + {\ptre \eta} \; \Opde\big) \hgiet +
\frac{J(\ptre\eta)}{\ptre} \; \partial_\eta \hgiet \times n  = O(\ptre^{\infty}), &
\hspace*{-0.1cm} \mbox{ in } \; 
\Gamma \times [0,\frac{\bar{\nu}}{\ptre}),\hspace*{-0.1cm}\\[12pt]  
\dsp iJ(\ptre\eta)\; \omega \, \hgiet  + (\Opun +
{\ptre \eta} \; \Opde) \egiet - \frac{J(\ptre\eta)}{\ptre} \, \partial_\eta
\egiet \times n = O(\ptre^{\infty}), & \hspace*{-0.1cm} \mbox{ in } \; \Gamma \times [0,
\frac{\bar{\nu}}{\ptre}), \hspace*{-0.1cm} 
\end{array} \right.
\end{equation}
%where we have dropped the $O(\ptre^{\infty})$ terms. 
These equations are complemented by the boundary condition
\begin{equation} \label{boundary_condition}
\hgiet(\xgg, 0) \times n + O(\ptre^{\infty}) = \hgee \times n \;(\xgg), \quad \xgg \in \Gamma.
\end{equation}
% completed by the boundary condition
% \begin{equation} \label{asymtotic_equation-1bord}
% \hgiet(\xgg, 0) \times n = \fue(\xgg) \times n + O(\ptre^\infty), \;\;   \xgg
% \in \Gamma. 
% \end{equation}
The substitution of (\ref{tuik}, \ref{ansatz_3}) into
(\ref{asymtotic_equation-1}, \ref{boundary_condition}) leads to a sequence of
problems that enable us
to  inductively determine the fields $(\egi^k, \hgi^k)$. The computations are
relatively delicate but straightforward. The most difficult task is to
explain the recurrence properly, which is the aim of this section. In 
Section~\ref{explicit}, we shall compute explicitly the first terms of the
expansions.\\[12pt]
It turns out to be very useful to
make a change of unknown concerning the electric field. This is motivated by
the observation that
\begin{equation}\label{premier terme}
\egi^0 = 0.
\end{equation}
This fact can be explained along the following lines: indeed from (\ref{tuik})
and (\ref{ansatz_3}) one deduces (at least formally) that:
$$
\|\egie\|^2_{L^2(\Oint)} \sim \; \delta \; \int_{\Gamma} \int_0^{+\infty} |\egi^0(\xgg, \eta)|^2
\; d \eta \; d \sigma.
$$
Therefore, the a priori estimate (\ref{estim-ex-1}), which says that
$\|\egie\|^2_{L^2(\Oint)} = O(\ptre^2)$, implies $\egi^0 = 0.$
\\[12pt]
The expansion for the electric field therefore starts with $
\ptre\;  \egi^1$ while for the magnetic field $\hgi^0\neq 0$. In some sense there is a natural shift of one power of $\ptre$ between
the expansions of the electric and magnetic fields. This is why we introduce
the ``normalized'' electric field: 
\begin{equation} \label{scaledE}
\nEi^{\ptre} = \frac{1}{\ptre} \; \egiet ,
\end{equation} 
and we seek an expansion of the form
\begin{equation} \label{new-expansion_E}
\nEi^{\ptre}(\xgg, \eta) := \nEi^0(\xgg, \eta) + \ptre\; 
\nEi^1(\xgg, \eta) + \ptre^2 \; \nEi^2(\xgg, \eta) + \cdots \;\;\; (\xgg,
\eta) \in \Gamma \times \R^+,
\end{equation} 
with the correspondence 
\begin{equation*} %\label{correspondance}
\egi^{k+1} = \nEi^k, \quad k \geq 1.
\end{equation*} 
We then rewrite (\ref{asymtotic_equation-1}) as a system of equations for
$(\nEi^{\ptre},\hgiet) $ (we have multiplied the first  equation by
$\ptre$) 
$$
\left\{ \begin{array}{ll}
\dsp J(\ptre \eta) \, \Big(i \eps_r \ptre \omega + \frac{1}{\omega }\Big)
\; \nEi^{\ptre} - \big(\ptre \; \Opun + 
{\ptre^2 \eta} \; \Opde\big) \hgiet + J(\ptre\eta) \; \partial_\eta
\hgiet \times n  = 0, \;\; & \mbox{ in } \; \Gamma \times [0,\frac{\bar{\nu}}{\ptre}),\\[12pt]
\dsp iJ(\ptre\eta)\; \omega \, \hgiet  + (\ptre \, \Opun +
{\ptre^2 \, \eta} \; \Opde) \; \nEi^{\ptre}  - {J(\ptre\eta)} \, \partial_\eta
\nEi^{\ptre} \times n = 0, \;\; & \mbox{ in } \; \Gamma \times [0,\frac{\bar{\nu}}{\ptre}),
\end{array} \right.
$$
~\\that we can rewrite by separating the ``$\ptre$-independent'' part, that we
keep in the left hand side, from the remaining terms, that we put in the
right hand side, as follows:
\begin{equation} \label{asymtotic_equation-new}
\left\{ \begin{array}{llll}
\dsp \; \; \partial_\eta \hgiet \times n + \frac{1}{\omega} \; \nEi^{\ptre}
& = & \dsp \sum_{\ell=1}^{4} \; \ptre^{\ell} A_H^{(\ell)} (\nEi^{\ptre},\hgiet)
 & \mbox{ in } 
\; \Gamma \times \R^+,  
\\[18pt]
\dsp - \; \partial_\eta \nEi^{\ptre} \times n + i \omega \, \hgiet & = &
\dsp\sum_{\ell=1}^{2} \; \ptre^{\ell} A_{E}^{(\ell)} (\nEi^{\ptre},\hgiet), &
\mbox{ in   } \; \Gamma \times \R^+. \end{array} \right.
\end{equation}
% \begin{equation} \label{asymtotic_equation-new-old}
% \left\{ \begin{array}{llll}
% \dsp \; \; \partial_\eta \hgiet \times n + \frac{1}{\omega} \; \nEi^{\ptre}
% & = & \dsp \sum_{\ell=1}^{3} \; \ptre^{\ell} A_{HH}^{(\ell)} \; \hgiet +
% \sum_{\ell=1}^{3} \; \ptre^{\ell} A_{HE}^{(\ell)} \; \nEi^{\ptre}, & \mbox{ in
%   } 
% \; \Gamma \times \R^+,  
% \\[12pt]
% \dsp - \; \partial_\eta \nEi^{\ptre} \times n + i \omega \, \hgiet & = &
% \dsp\sum_{\ell=1}^{3} \; \ptre^{\ell} A_{EH}^{(\ell)} \; \hgiet +
% \sum_{\ell=1}^{3} \; \ptre^{\ell} A_{EE}^{(\ell)} \; \nEi^{\ptre}, & \mbox{ in
%   } \; \Gamma \times \R^+. \end{array} \right.
% \end{equation}
The linear operators $\big \{ A_{H}^{(\ell)},  \; \ell = 1,2,3 \; \big \}$ are given
by:
\begin{equation*} %\label{lesoperateursA_H}
\left\{ \begin{array}{lll}
\dsp A_{H}^{(1)} \; (u,v) = \Opun \; v - 2 h \eta \; \big( \partial_\eta v
\times n + \frac{1}{\omega} \; u\big),\\[12pt] 
\dsp A_{H}^{(2)} \; (u,v) = - i \eps_r \omega \; u + \eta \; \Opde \; v - g
\eta^2 \; \big( \partial_\eta v 
\times n + \frac{1}{\omega} \; u\big),\\[12pt]
\dsp A_{H}^{(3)} \; (u,v) = - 2 \eta h \; i \eps_r \omega \; u,\\[12pt]
\dsp A_{H}^{(4)} \; (u,v) = - 2 \eta^2 g \; i \eps_r \omega \; u,\\[12pt]
\end{array}\right.
\end{equation*}
and the linear operators $\big \{ A_{E}^{(\ell)}, \; \ell = 1,2 \; \big \}$
are given by: 
\begin{equation*}% \label{lesoperateursA_E}
\left\{ \begin{array}{lll}
\dsp A_{E}^{(1)} \; (u,v) = - \Opun \; u + 2 h \eta \; \big( \partial_\eta u
\times n - i {\omega} \; v\big) ,\\[12pt]
\dsp A_{E}^{(2)} \; (u,v) = - \eta \; \Opde \; u + g \eta^2 \; \big(
\partial_\eta u \times n - i {\omega} \; v\big).
\end{array}\right.
\end{equation*}
% The linear operators $\big \{ \; (A_{HH}^{(\ell)}, A_{HE}^{(\ell)},
% A_{EH}^{(\ell)}, A_{EE}^{(\ell)} \ ), \; \ell = 1,2,3 \; \big \}$ are defined
% by:\\ 
% \begin{equation} \label{lesoperateursA-old}
% \begin{array}{ll}
% \left\{ \begin{array}{lll}
% \dsp A_{HH}^{(1)} \; u = \Opun  u - 2 \eta h \; \partial_\eta u \times n,\\[12pt]
% \dsp A_{HH}^{(2)} \; u = \eta \; \Opde  u - \eta^2 g \; \partial_\eta u \times
% n,\\[12pt] 
% \dsp A_{HH}^{(3)}= 0,
% \end{array} \right.
% &
% \left\{ \begin{array}{lll}
% \dsp A_{HE}^{(1)} \; u = - \, \big( i \eps_r \omega + \frac{1}{\omega}) \; u,\\[12pt]
% \dsp A_{HE}^{(2)} \; u = - \, (2 \eta  h \; i \eps_r \omega + \eta^2
% \frac{g}{\omega})\; u,\\[12pt] 
% \dsp A_{HE}^{(3)}= - \, \eta^2 g \; i  \eps_r \omega \; u,
% \end{array} \right. 
% \vspace*{1cm}
% \\
% \left\{ \begin{array}{lll}
% \dsp A_{EH}^{(1)} \; u = - \, 2 \; \eta  h \; i \omega \; u,\\[12pt]
% \dsp A_{EH}^{(2)} \; u = - \, \eta^2 g \; i \omega \; u,\\[12pt]
% \dsp A_{EH}^{(3)}= 0,
% \end{array} \right.
% &
% \left\{ \begin{array}{lll}
% \dsp A_{EE}^{(1)} \; u = - \, \Opun  u + 2 \eta  h \;
%     \partial_\eta u \times n,\\[12pt]
% \dsp A_{EE}^{(2)} \; u = - \, \eta \; \Opde  u + 2 \eta^2 g \; \partial_\eta u
% \times 
% n,\\[12pt] 
% \dsp A_{EE}^{(3)}= 0.
% \end{array} \right.
% \end{array}
% \end{equation}
~\\Substituting (\ref{tuik}) and (\ref{new-expansion_E})
into (\ref{asymtotic_equation-new}) then equating the same powers of $\ptre$ leads to the
following systems:
\begin{equation} \label{induction-ehh}
\left\{ \begin{array}{llll}
\dsp \; \; \partial_\eta \hgi^k \times n + \frac{1}{\omega} \; \nEi^{k}
& = & \dsp \sum_{\ell=1}^{4} \;  A_H^{(\ell)}
(\nEi^{k - \ell},\hgi^{k - \ell})
, & \mbox{
  in } \; \Gamma \times \R^+,  
\\[18pt]
\dsp - \; \partial_\eta \nEi^{k} \times n + i \omega \, \hgi^k & = &
\dsp \sum_{\ell=1}^{2} \; A_E^{(\ell)}
(\nEi^{k - \ell},\hgi^{k - \ell}), & \mbox{ in
  } \; \Gamma \times \R^+, \end{array} \right.
\end{equation}
for $k = 0, 1, 2, \cdots$, with the convention $\nEi^{\ell} = \hgi^{\ell}= 0$
for $\ell < 0$.\\[12pt]
Of course, these equations have to be complemented with the conditions (see
(\ref{boundary_condition}) and (\ref{cond1_EHik}))
\begin{equation} \label{cond_infini_limites_EH}
\left\{ \begin{array}{l}
\dsp  \hgi^k(\xgg, 0) \times n = \hge^k(\xgg, 0),\\[12pt] 
\dsp  \int_0^{+ \infty} |\hgi^k(\xgg,
\eta)|^2 \; d \eta < + \infty, \quad  \mbox{ and } \quad \int_0^{+ \infty} |\nEi^{k}(\xgg, 
\eta)|^2 \; d \eta  < + \infty, 
\end{array} \right.
\end{equation}
$\forall \; \xgg \in \Gamma.$ 

The reader can already notice how  the roles of the
variables 
$\eta$ and $\xgg$ have been separated. The variable $\xgg$ appears as parameter for
determining $(\nEi^{k}, \hgi^{k})$ from the previous $(\nEi^{\ell}
,\hgi^{\ell})$'s since, for each $\xgg$, one simply has to solve an ordinary differential
system in the variable $\eta$. The solutions to this inductive system of
equations can be expressed in a general way using the result of the following
technical lemma. For that purpose it is useful to introduce 
\begin{equation*}% \label{defPk}
{\bf P}_k(\Gamma, \R^+; \C^3) := \Big\{ u(\xgg, \eta) = \sum_{j=1}^k a_j(\xgg)
\; \eta^j, \; a_j \in C^{\infty}(\Gamma; \C^3) \; \Big\}.
\end{equation*}
\begin{lemma} \label{technique}
Let $(f,g) \in {\bf P}_k(\Gamma, \R^+; \C^3)^2$ and $\varphi \in
C^{\infty}(\Gamma; \R^3)$, Then  the problem,
\begin{equation} \label{pb_technique}
\begin{array}{l}
\mbox{Find } (u,v) \in C^{\infty}(\Gamma; C^{\infty}(\R^+))^2 \mbox{ such
  that, }\\[6pt]
\left\{ \begin{array}{llll}
\dsp \; \; \; \partial_\eta v \times n + \frac{1}{\omega} \; u
 =  \dsp e^{- \sqrti \, \eta} \; f(\eta, \cdot),  & \mbox{ in
  } \; \Gamma \times \R^+,
\\[6pt]
\dsp - \; \partial_\eta u \times n + i \omega \, v   = \dsp e^{-
 \sqrti \, \eta} \; g(\eta, \cdot),  & \mbox{ in
  } \; \Gamma \times \R^+. \end{array} \right.
\end{array}
\end{equation}
with the conditions: 
\begin{equation} \label{cond_infini_limites}
\left\{ \begin{array}{l}
\dsp \forall \; \xgg \in \Gamma, \quad u(\xgg, 0) \times n = \varphi(\xgg),
\\[12pt] 
\dsp \forall \; \xgg \in \Gamma, \quad \int_0^{+ \infty} |u(\xgg,
\eta)|^2 \; d \eta < + \infty, \quad  \int_0^{+ \infty} |v(\xgg, 
\eta)|^2 \; d \eta  < + \infty, 
\end{array} \right.
\end{equation}
has a unique solution, which is of the form
\begin{equation} \label{forme_solution}
u(\xgg, \eta) = e^{- \sqrti \, \eta} \; p(\xgg, \eta) \quad \mbox{ and }
\quad v(\xgg, \eta) =
e^{- \sqrti \, \eta}
\; q(\xgg, \eta)
\end{equation}
with $(p,q) \in {\bf P}_{k+1}(\Gamma, \R^+; \C^3)^2$ and with the square root
definition $\sqrti := \frac{\sqrt{2}}{2} \; (1 + i)$.
\end{lemma}

\proof
This is a simple exercise on ordinary differential equations. 
For the
uniqueness, assuming $\varphi = 0$ and $f=g=0$, we can eliminate $v$ between
the two equations of (\ref{pb_technique}) and obtain, with $u_T = n \times(u \times n) = u - (u \cdot n) \; n$:
$$
\partial_{\eta \eta} u_T - i \; u_T = 0.
$$
Since we also have $u_T(\xgg,0) = 0$, the only solution satisfying the second condition of (\ref{cond_infini_limites})
 is $u_T = 0$. One deduces from the 
second equation of (\ref{pb_technique}) that  $v=0$ which in turn implies $u=0$
using the first  equation of (\ref{pb_technique}). \vspace{2mm}\\
We proceed in the same way to prove the existence of solutions of the
form (\ref{forme_solution}). First, the projection on $n$ of the two equations
 directly gives
$$
\left\{ \begin{array}{l}
(u \cdot n) = \omega \; e^{- \sqrti \, \eta} \; (f \cdot n), \\[12pt]
(v \cdot n) = \omega \; e^{- \sqrti \, \eta} \; (f \cdot n).
\end{array} \right.
$$
For the tangential components, one easily gets
$$
\left\{ \begin{array}{l}
\partial_{\eta \eta} u_T - i \; u_T = \widetilde{f}_T \; e^{- \sqrti \, \eta}, \\[12pt]
\partial_{\eta \eta} v_T - i \; v_T = \widetilde{g}_T \; e^{- \sqrti \, \eta},
\end{array} \right.
$$
where 
$$
\left\{ \begin{array}{l}
\dsp \widetilde{f}_T := \big(\partial_{\eta} g - \sqrti \; \eta g \big) \times n + i
\omega f_T, \quad  \in {\bf P}_{k}(\Gamma, \R^+; \C^3),\\[12pt]
\dsp \widetilde{g}_T := n \times \big(\partial_{\eta} f - \sqrti \; \eta f \big) -
\frac{1}{\omega} g_T, \quad  \in {\bf P}_{k}(\Gamma, \R^+; \C^3).
\end{array} \right.
$$
Thus, the key point is the resolution of the scalar 
differential  equation:
\begin{equation} \label{pb_scalaire}
\partial_{\eta \eta} \psi - i \psi = s(\eta) \; e^{- \sqrti \, \eta}, \quad \mbox{in } \R^+, 
\end{equation}
with $\psi(0) = \psi_0$. The unknown $\psi$ is sought in $L^2(\R^+)$ and
$s(\eta)$
is given in $P_m$, the space of polynomials in $\eta$ of degree less than $m$.
We prove that the solution of this equation has the form:
$$
\psi(\eta) = e^{- \sqrti \, \eta} \; p(\eta) \quad \mbox{with } p \in P_{m+1}.
$$
Indeed, after subtracting  $$\psi_0 \; e^{- \sqrti \, \eta},$$
we can assume that $\psi_0 = 0$. Introducing the space:
$$
\widetilde{P}_{m+1} = \left\{ p \in P_m \; /\; p(0) = 0 \right \} ,
$$
Observing that $\left[ \partial_{\eta \eta} - i \right] \big\{ e^{- \sqrti \,
  \eta} \; x^{\ell } \big\}= e^{- \sqrti \, \eta} \; \big\{ \ell \; x^{\ell-2}
\left[ 2 
\sqrti \; x + \ell-1 \right]\big\}$, we deduce that
$$
\left( \partial_{\eta \eta} - i \right) \in {\cal L}(e^{- \sqrti \, \eta} \;\widetilde{P}_{m+1}, e^{- \sqrti \, \eta} \;P_m).
$$
From the uniqueness result for $L^2$ solutions of (\ref{pb_scalaire}) with
Dirichlet condition at $\eta=0$, we deduce that the
operator $\left( \partial_{\eta \eta} - i \right)$ is injective in the space
$e^{- \sqrti \, \eta} \; \widetilde{P}_{m+1}$. Since
$e^{- \sqrti \, \eta} \;P_m$ and $e^{- \sqrti \, \eta} \;\widetilde{P}_{m+1}$ have the same dimension, this operator is an
isomorphism from $e^{- \sqrti \, \eta} \;P_m$ into $e^{- \sqrti \, \eta} \;\widetilde{P}_{m+1}$, which concludes the proof.
 \proofend
 
\noindent As an application of this lemma we obtain the following result.
\begin{theorem} \label{lesEkHk}
The fields $ H_e^k \times n \in C^{\infty}(\Gamma; \C^3)$ being given, there exists a
unique sequence $$\left\{ \;(\nEi^{k},\hgi^k) \in {C}^{\infty} (\Gamma; \C^3)^2,
\; k = 0, 1, 2, \cdots \right\}$$ satisfying the sequence of
problems (\ref{induction-ehh})-(\ref{cond_infini_limites_EH}). Moreover, 
\begin{equation} \label{formEH}
(e^{\sqrti \, \eta} \;\nEi^{k}, \; e^{\sqrti \, \eta} \;\hgi^k) \; \in \; 
{\bf P}_k(\Gamma, \R^+; \C^3)^2.
\end{equation}
\end{theorem} 

\proof
We prove this theorem using an induction on $k$. Suppose that the existence and uniqueness of
$(\nEi^{\ell},\hgi^{\ell})$ and the property
(\ref{formEH}) have been guaranteed up to $k-1$ and
 \begin{equation} \label{recEH1}
\big( e^{\sqrti \, \eta} \; \nEi^{m}, \; e^{\sqrti \, \eta} \hgi^{m} \big) \in
{\bf P}_{m}(\Gamma, \R^+; \R^3)^2, \quad m = 0, \cdots, k-1.
\end{equation}
We also include in the recurrence the assumption (by convention ${\bf
  P}_{-1}(\Gamma, \R^+; \R^3) = 0$):
\begin{equation} \label{recEH2}
\left\{ \begin{array}{l}
e^{\sqrti \, \eta} \; \big(- \; \partial_\eta \nEi^{m} \times n + i \omega
\hgi^{m}\big)\in  {\bf P}_{m-1}(\Gamma, \R^+;\C^3)^2 \quad m = 0, \cdots, k-1,
\\[12pt] 
\dsp e^{\sqrti \, \eta} \; \big(\partial_\eta \hgi^{m} \times n +
\frac{1}{\omega} \; \nEi^{m}\big) \in  {\bf P}_{m-1}(\Gamma, \R^+;
\C^3)^2, \quad m = 0, \cdots, k-1.
\end{array} \right.
\end{equation}
According to Lemma \ref{technique}, to prove (\ref{recEH1}) and (\ref{recEH2})
for $m=k$, it suffices to show that
the two right hand sides of (\ref{induction-ehh}) are of the form:
$$
e^{-\sqrti \, \eta} \; p \quad \mbox{with} \quad p \in 
{\bf P}_{k-1}(\Gamma, \R^+; \C^3)^2.
$$
To verify this, let us consider $(u,v)$ satisfying:
\begin{equation}\label{recEH2-b}
\left\{ \begin{array}{l} 
\big( \; e^{\sqrti \, \eta} \; u, \; e^{\sqrti \, \eta} \; v \big) \in  {\bf
  P}_{m}(\Gamma, \R^+; \C^3)^2,
\\[12pt]
\big( \; e^{\sqrti \, \eta} \; (- \; \partial_\eta u \times n + i \omega v), \;
e^{\sqrti \, \eta} \; (\partial_\eta v \times n + 
\frac{1}{\omega} \; u ) \; \big) \in  {\bf P}_{m-1}(\Gamma, \R^+;\C^3)^2 .
\end{array}\right.
\end{equation}
The fact that the operators $A_{H}^{(\ell)}$, $\ell = 3,\,4$, are
polynomials of degree $2$ in $\eta$, means that 
\begin{equation} \label{prop1}
e^{\sqrti \, \eta} \; A_{H}^{(\ell)} \; (u,v) \in {\bf P}_{m + 2}(\Gamma,
\R^+;\C^3), \mbox{for } \ell = 3,\,4.
\end{equation}
Moreover, using assumption (\ref{recEH2-b}) and the explicit form of
$A_{H}^{(\ell)}$ and $A_{H}^{(\ell)}$ for $\ell = 1, 2$,  we observe that
\begin{equation} \label{prop2}
\left \{ \begin{array}{l}
(\; e^{\sqrti \, \eta} \; A_{H}^{(1)} \; (u,v), \; e^{\sqrti \, \eta} \;
A_{E}^{(1)} \; (u,v) \:) \in {\bf P}_{m}(\Gamma,\R^+;\C^3)^2,\\[12pt]
(\; e^{\sqrti \, \eta} \; A_{H}^{(2)} \; (u,v), \; e^{\sqrti \, \eta} \;
A_{E}^{(2)} \; (u,v) \:) \in {\bf P}_{m+1}(\Gamma,\R^+;\C^3)^2.
\end{array} \right.
\end{equation}
Therefore, since  $e^{\sqrti \, \eta}
\;(\nEi^{k-\ell},\hgi^{k-\ell}) \in {\bf P}_{k-\ell}(\Gamma, \R^+;\R^3)^2$
(according to (\ref{recEH1})),
we deduce from (\ref{prop1}) that 
\begin{equation*}% \label{part1}
e^{\sqrti \, \eta} \; \sum_{l=3}^4 A_{E}^{(\ell)} \;
(\nEi^{k-\ell},\hgi^{k-\ell}) \in {\bf P}_{k-1}(\Gamma, \R^+;\C^3)^2 , 
\end{equation*}
while, thanks to (\ref{recEH2}) and by noticing that $k
-\ell + 2 \leq k-1$ for $ \ell \geq 3$, we deduce from (\ref{prop2}) that
\begin{equation*} %\label{part2}
\left\{ \begin{array}{l}
\dsp e^{\sqrti \, \eta} \; \sum_{l=1}^2 \; A_{E}^{(\ell)} \;
(\nEi^{k-\ell},\hgi^{k-\ell}) \in {\bf P}_{k-1}(\Gamma, \R^+;\C^3)^2,\\[12pt]
\dsp e^{\sqrti \, \eta} \; \sum_{l=1}^2 \; A_{H}^{(\ell)} \;
(\nEi^{k-\ell},\hgi^{k-\ell}) \in {\bf P}_{k-1}(\Gamma, \R^+;\C^3)^2 .
\end{array} \right.
\end{equation*}
This  leads to the desired property.

\vskip2mm

\noindent To end the proof, it suffices to check the result for $k=0$ and $k=1$,
which will be done in Section~\ref{explicit} (see Remarks~\ref{rec1}-\ref{rec2}). 
\proofend

\subsection{Explicit computation of the interior fields for $k =1,2,3$}
\label{explicit} 
This section is devoted to the presentation of the technical  details related
to the computation of the asymptotic terms up to the order $k=3$. 
In the sequel, we shall  systematically use the following formulas, 
deduced from (\ref{notationC}),
\begin{equation} \label{util}
(\Opun V) \cdot n = \rotg V_T, \quad \Opun V \times n = \Rotg (V \cdot n) \times n - \big( {\cal C} V \times n \big) \times n,
\end{equation}
\begin{equation} \label{utilM}
(\Opde V) \cdot n = \rotgm V_T, \quad \Opde V \times n = \Rotgm (V \cdot n) \times n - \big( {\cal G} V \times n \big) \times n.
\end{equation}
{\bf Computation of $(\nEi^0 \equiv E_i^1, \hgi^0)$}.
For $k=0$, (\ref{induction-ehh}) gives
\begin{equation} \label{eq_k=0}
\left\{ \begin{array}{llll} 
\dsp \; \; \partial_\eta \hgi^0 \times n + \frac{1}{\omega} \; \nEi^{0}
& = & 0, & \mbox{
  in } \; \Gamma \times \R^+,   
\\[12pt] 
\dsp - \; \partial_\eta \nEi^{0} \times n + i \omega \, \hgi^0, & = & 
0 & \mbox{ in
  } \; \Gamma \times \R^+. \end{array} \right. 
\end{equation}
whose unique $L^2$ solution satisfying $\hg_{i,T}^0(\xgg, \eta) =
H^0_{e,T}(\xgg)$ is given by:
\begin{equation} \label{sol_k=0}
\left\{ \begin{array}{l}
\nEi^{0}(\xgg, \eta) \equiv E_i^1(\xgg, \eta)= \sqrti \;  \omega \;
\big(H^0_{e} \times n\big) (\xgg) \; e^{-\sqrti \, \eta},
\\[12pt]
\hg_{i}^0(\xgg, \eta)= H^0_{e,T}(\xgg) \; e^{-\sqrti \, \eta},
\end{array} \right. 
\end{equation}
from which we deduce the useful information for the construction of the
GIBCs, namely:
\begin{equation} \label{traceET_k=1}
E_i^1 \times n \; (\xgg, 0) = - \sqrti \; \omega \;  H^0_{e,T}  (\xgg).
\end{equation}
\begin{remark} \label{rec1}  Notice that (\ref{sol_k=0}) proves in particular
  Theorem~\ref{lesEkHk} for $k=0$.
\end{remark}
{\bf Computation of $(\nEi^1 \equiv E_i^2, \hgi^1)$}.
For $k=1$, (\ref{induction-ehh}) gives, using (\ref{eq_k=0})
\begin{equation} \label{eq_k=1}
\left\{ \begin{array}{llll} 
\dsp \; \; \partial_\eta \hgi^1 \times n + \frac{1}{\omega} \; \nEi^{1}
& = & \Opun \hgi^{0}, & \mbox{
  in } \; \Gamma \times \R^+,   
\\[12pt] 
\dsp - \; \partial_\eta \nEi^{1} \times n + i \omega \, \hgi^1 & = & 
- \Opun \nEi^{0},& \mbox{ in
  } \; \Gamma \times \R^+. \end{array} \right. 
\end{equation}
We project (\ref{eq_k=1}) on $n$, use (\ref{util}) 
and the expressions (\ref{sol_k=0}) for $\nEi^{0}$ and $\hgi^{0}$, to obtain: 
\begin{equation} \label{solnorm_k=1}
\left\{ \begin{array}{ll} 
\dsp  \nEi^{1} \cdot n = \omega \; \Opun \hgi^{0}\cdot n = \omega \big[ \rotg
H^0_{e,T}\big] (\xgg) \; e^{-\sqrti \, \eta} , &  \mbox{
  in } \; \Gamma \times \R^+,   
\\[12pt] 
\hgi^1 \cdot n  =  
\dsp\frac{i}{\omega} \; \Opun \nEi^{0} \cdot n =  - \frac{1}{\sqrti} \;\big[
\rotg \big(H^0_{e} \times n\big)\big] (\xgg) \; e^{-\sqrti \, \eta} & \mbox{ in
  } \; \Gamma \times \R^+. \end{array} \right. 
\end{equation}
Next, we eliminate $\nEi^{1}$ in (\ref{eq_k=1}) and get the following equation in $H^1_{i,T}$ 
$$
\left(\partial^2_{\eta\eta} - i \right)H^1_{i,T}= n \times \partial_{\eta}
\big[\Opun H_i^0 \big]- \frac{1}{\omega} \; n \times \big(\Opun \nEi^{0}
\times n \big)
$$
We use again (\ref{util}) and (\ref{sol_k=0}) to transform the right hand
side. Using the following identity, that can easily be deduced from the
definitions (\ref{courbures}) and (\ref{courbures2})
\begin{equation*}% \label{id1}
\big({\cal C}(V \times n)\big)\times n - {\cal C}V = -2{\cal H} \; V \;\; \mbox{ for all } V \in \R^3,
\end{equation*}
we finally get after some easy manipulations
\begin{equation*} %\label{eqHT_k=1}
\left(\partial^2_{\eta\eta} - i \right)H^1_{i,T} = 2 \sqrti \; {\cal H}
\; H^0_{e,T}(\xgg) \; e^{-\sqrti \, \eta},
\end{equation*}
whose unique $L^2$ solution satisfying $\hgi^0(\xgg, \eta) =
H^0_{e,T}(\xgg)$ is given by:
\begin{equation} \label{HT_k=1}
H^1_{i,T}(\xgg, \eta) = \Big( H^1_{e,T}(\xgg) - \eta \; {\cal H} \;
H^0_{e,T}(\xgg) \Big) \; e^{-\sqrti \, \eta}.
\end{equation}
Coming back to the first equation of (\ref{eq_k=1}), we get
\begin{equation} \label{ET_k=1}
\mE^1_{i}\times n \; (\xgg, \eta) = \omega \; 
\Big( -\sqrti\, H^1_{e,T}(\xgg)+ ({\cal C}-{\cal H})\, H^0_{e,T}(\xgg)
+ \eta \; \sqrti \; {\cal H} \; H^0_{e,T}(\xgg) \Big)\;  e^{-\sqrti \, \eta}. 
\end{equation}
In particular 
\begin{equation} \label{traceET_k=2}
E_i^2 \times n \; (\xgg, 0) = \omega \; \Big( - \sqrti \;  H^1_{e,T}  (\xgg) + ({\cal C}-{\cal H}) \, H^0_{e,T}(\xgg) \Big) .
\end{equation}
\begin{remark} \label{rec2}  Notice that (\ref{solnorm_k=1}), (\ref{HT_k=1})
  and (\ref{ET_k=1}) prove in particular Theorem~\ref{lesEkHk} for $k=1$.
\end{remark}
{\bf Computation of $(\nEi^2 \equiv E_i^3, \hgi^2)$}. The calculations are much harder and tedious than for the two previous cases. That is why we shall restrict ourselves to the main steps. Also, for the sake of simplicity, we shall often omit to mention the dependence of the various quantities we manipulate with respect to $\xgg$.\\[12pt]
For $k=2$, (\ref{induction-ehh}) gives, using (\ref{eq_k=1})
\begin{equation} \label{eq_k=2}
\hspace*{-0.2cm}\left\{ \begin{array}{llll} 
\dsp \; \; \partial_\eta \hgi^2 \times n + \frac{1}{\omega} \; \nEi^{2}
\hspace*{-0.1cm} & = & \hspace*{-0.1cm} r_H^2, & \mbox{ in
  } \; \Gamma \times \R^+
\\[12pt] 
\dsp - \; \partial_\eta \nEi^{2} \times n + i \omega \, \hgi^2
\hspace*{-0.1cm} & = & \hspace*{-0.1cm} r_E^2, &
\mbox{   in } \; \Gamma \times \R^+,   
\end{array} \right. 
\end{equation}
where we have set
\begin{equation*}% \label{RHSeq_k=2}
\hspace*{-0.2cm}\left\{ \begin{array}{llll} 
r_H^2 \hspace*{-0.1cm} & = & \hspace*{-0.1cm} 
\Opun \hgi^{1} - i \; \eps_r \omega \, \mE^0_{i} + 
\eta \; \big(\Opde - 2 h \; \Opun \big) \hgi^{0},
\\[12pt] 
r_E^2 \hspace*{-0.1cm} & = & \hspace*{-0.1cm}  
- \Opun \nEi^1 - \eta \; \big(\Opde - 2 h \; \Opun \big) \mE^0_{i} .
\end{array} \right. 
\end{equation*}
We can go directly to the evaluation of $H_{i,T}^2$ which satisfies (apply $n \times \partial_{\eta}$ to the first equation of (\ref{eq_k=2}), divide the second equation by $\omega$ and add the two results)
\begin{equation} \label{eqHT_k=2}
\left(\partial^2_{\eta\eta} - i \right) H_{i,T}^2= n \times \partial_{\eta} \,
r_{H}^2 - \frac{1}{\omega} \; r_{E,T}^2. 
\end{equation}
The next step consists in expressing the right hand side of (\ref{eqHT_k=2})
 in terms of the previous $(\mE_i^{\ell}, H_i^{\ell})$'s. 
Using (\ref{util}), (\ref{utilM}) and the fact that $H_i^0 \cdot n = 0$, we first compute that
$$
  n \times r_H^2 = n \times \Rotg \big(\hgi^1 \cdot n \big) - {\cal C} H_{i,T}^1 - \eta \; \big( {\cal G} - 2  h {\cal C} \big) H_{i,T}^0 - i \; \eps_r \omega \, \big( n \times \mE^0_{i}  \big),
$$
Next, we use the expressions (\ref{sol_k=0}), (\ref{solnorm_k=1}) and
(\ref{HT_k=1}) and the identity
$$
n \times \Rotg \big( \rotg (V \times n) \big) = - \grag \, \big( \divg V \big)
$$
to obtain
$$
  n \times r_H^2  =  \left[ \frac{1}{\sqrti} \; \big( \grag \, \divg + \eps_r \omega^2 \big)  H^0_{e,T}-  {\cal C} H^1_{e,T} \right] \; e^{- \sqrti \eta} 
+ \eta \; \big(3 h {\cal C} - {\cal G}\big) H^0_{e,T} \; e^{- \sqrti \eta} \; .
$$
After differentiation, we get 
\begin{equation} \label{der_rH}
\left| \begin{array}{lll}
  n \times \partial_{\eta}r_H^2  &= & 
\left[ - \big( \grag \, \divg + \eps_r \omega^2 \big)  H^0_{e,T} 
+ \big(3 h {\cal C} - {\cal G}\big) H^0_{e,T}
 + \sqrti \; {\cal C} H^1_{e,T} \right] \; e^{- \sqrti \eta}
\\[12pt]
 & &  \quad -\eta \, \sqrti \; \big(3 h {\cal C} - {\cal G}\big) H^0_{e,T}\; e^{- \sqrti \eta} \; .
\end{array} \right.
 \end{equation}
In the same way, using again (\ref{util}), (\ref{utilM}) and the fact that $\mE_i^0 \cdot n =0$, we calculate 
$$
 r_{E,T}^2 = n \times ( r_{E}^2 \times n) = - \rotg \big(\mE_i^1 \cdot n \big) + 
\big( {\cal C} \mE_i^1 \big) \times n + \eta \; \big[ \big( {\cal G} - 2 h {\cal C}\big) \mE_i^0  \big] \times n.
$$
Next, we notice that $${\cal C} V  = - {\cal C} 
\big[ (V \times n)\times n \big]  \quad \mbox{(and the same with ${\cal G})$}$$
 and use the expressions 
(\ref{sol_k=0}), (\ref{solnorm_k=1}) and (\ref{ET_k=1}) respectively for 
$\mE_i^0 \times n$, $\mE_i^1 \cdot n$
 and $\mE_i^1 \times n$ to obtain
$$
\left| \begin{array}{lll}
\dsp - \frac{1}{\omega} \; r_{E,T}^2  &= & \left[\;
\Rotg \big( \rotg H^0_{e,T}\big) +  \eta \; \sqrti
 \; \Big( \big({\cal G} - 3 h \, {\cal C} \big) H_{e,T}^0 \times n\Big) \times n
\; \right] \; e^{- \sqrti \eta}
\\[12pt]
& + &
\left[\; {\cal C} \Big( \big( {\cal H} - {\cal C}) H_{e,T}^0 \times n\Big) \times n - 
\sqrti \; \Big( {\cal C} \big( H_{e,T}^1 \times n \big) \Big) \times n \; \right] 
\; e^{- \sqrti \eta}
\end{array} \right.
$$
This can be written in a simplified form, using the following identities that
hold for all $V \in \R^3$ and that are easily deduced from (\ref{courbures}) and (\ref{courbures2})
$$
\left\{ \begin{array}{l}
 \left\{{\cal C} ( ({\cal H} - {\cal C}) V) \times n) \right\} \times n = (3 h{\cal C} - {\cal C}^2 - 2 {\cal H}^2)V, 
\\[12pt]
({\cal C}(V \times n))\times n = ({\cal C} - 2 {\cal H}) V,
\\[12pt]
\left\{(3{\cal H} {\cal C}-{\cal G})(V\times n)\right\} \times n = (3{\cal H} {\cal C} + {\cal G}  - 6 {\cal H}^2) V.
\end{array} \right.
$$
We obtain
\begin{equation} \label{rET}
\left| \begin{array}{lll}
\dsp - \frac{1}{\omega} \; r_{E,T}^2  &= & 
\left[ \Rotg \big( \rotg H^0_{e,T}\big) + 
\big(3{\cal H} {\cal C} - {\cal C}^2 - 2 {\cal H}^2\big) H^0_{e,T} 
\right] \; e^{- \sqrti \eta}
\\[12pt]
& - &  \sqrti \; ({\cal C}- 2 {\cal H}) H^1_{e,T}  
\; e^{- \sqrti \eta} + \; \eta \; \sqrti \; \big(3{\cal H} {\cal C} + {\cal G}  - 
6 {\cal H}^2\big) \; H^0_{e,T} \; e^{- \sqrti \eta}. 
\end{array} \right.
\end{equation}
Substituting (\ref{der_rH}) and (\ref{rET}) in (\ref{eqHT_k=2}) leads to the following equation
 \begin{equation*}% \label{eh2}
\left|\begin{array}{lcll}
\dsp \left(\partial^2_{\eta\eta} - i \right)\hg^{2}_{i,T} &=& e^{-\sqrti \, \eta} &\left\{2 \; \sqrti \; {\cal H} \; H^1_{e,T} + ({\cal C}^2 + 2 {\cal H}^2  -{\cal G}) \, H^0_{e,T} \right.  \\[12pt]
&&& \left.\;  - (\Deltag + \eps_r \omega^2) \; H^0_{e,T} 
- \eta \; \sqrti \; \big(6{\cal H}^2 - 2 {\cal G}\big) \; H^0_{e,T}\right\},
\end{array} \right.
\end{equation*}
where $\Deltag : = \grag \divg -\Rotg \rotg$ is the vectorial Laplace Beltrami operator.\\[12pt]
Since the $L^2$ solution to 
$$
\dsp (\partial_{\eta\eta}^2 - i) u = (a + b \, \eta)\,e^{-\sqrti \, \eta} \quad \mbox{in } \R^+, 
$$
is given by
$$
u(\eta) = \left(u(0) + \left(\frac{a}{2\sqrti} - \frac{b}{4i}\right) \eta + \frac{b}{4\sqrti} \; \eta^2 \right)\, e^{-\sqrti \, \eta},
$$
we deduce that $\hg^{2}_{T}$ is given by the expression
\begin{equation*}% \label{expr-ht2}
\left| \begin{array}{ll}
\hg^{2}_{i,T}(\xgg, \eta) = e^{-\sqrti \, \eta} &\left\{ H^2_{e,T} - \eta \;
  {\cal H} H^1_{e,T}- \dsp \frac{\eta}{2\sqrti}\; \big ({\cal C}^2-{\cal H}^2\big)\; H^2_{e,T}
\right.
\\[12pt]
& \dsp \left.
+\frac{\eta}{2\sqrti}\; \big (\Deltag + \eps_r \omega^2 \big) \; H^2_{e,T}+ \frac{\eta^2}{2} \; \big(3{\cal H}^2 - {\cal G}\big) H^1_{e,T} \right\}. 
\end{array} \right.
\end{equation*}
Finally we go back to the first equation of to obtain, after lengthy
calculations that we do not detail here, 
\begin{equation*}% \label{expr-et3}
\left| \begin{array}{llll}
\mE^{2}_{i,T} \times n = \dsp \omega \; e^{-\sqrti \, \eta} \; \Big\{ \hspace*{-0.2cm} & \dsp
-\sqrti \; H^2_{e,T}+ ({\cal C} -{\cal H}) \; H^1_{e,T} -\frac{1}{2 \sqrti} \;
 \big({\cal C}^2 - {\cal H}^2 \big) \; H^0_{e,T}
\\[12pt]& \dsp-\frac{1}{2 \sqrti} \; \big( \,\eps_r \omega^2 + 
\grag \divg + \Rotg \rotg \big) \,
H^0_{e,T}
\\[12pt]
&\dsp + \; \eta \; \left(\sqrti{\cal H} H^1_{e,T} + \frac{1}{2} \big (5 {\cal H}^2 -6 {\cal H}{\cal C} + {\cal C}^2 - \Deltag -\eps_r \omega^2 \big)\, H^0_{e,T} \right)
\\[12pt]
& \dsp - \; \eta^2 \; \frac{\sqrti}{2} \; \big(3{\cal H}^2 - {\cal G}\big) H^0_{e,T} \; \Big \}. 
\end{array} \right.
\end{equation*}
In particular, for $\eta = 0$,
\begin{equation} \label{traceET_k=3}
\left| \begin{array}{llll}
E^{3}_{i,T} \times n = \dsp \omega \; e^{-\sqrti \, \eta} \; \Big\{ \hspace*{-0.2cm} & \dsp
-\sqrti \; H^2_{e,T}+ ({\cal C} -{\cal H}) \; H^1_{e,T} -\frac{1}{2 \sqrti} \;
 \big({\cal C}^2 - {\cal H}^2 \big) \; H^0_{e,T}
\\[12pt]& \dsp-\frac{1}{2 \sqrti} \; \big( \,\eps_r \omega^2 + 
\grag \divg + \Rotg \rotg \big) \,
H^0_{e,T}\; \Big \}.
\end{array} \right.
\end{equation}
\subsection{Construction of the GIBCs}
\label{construc}
The GIBCS  of order $k$  is obtained by considering the
truncated expansions in $\Oext$
\begin{equation*} %\label{defEetronc}
\eg^{\ptre}_{e,k} := \sum_{\ell = 0}^{k} \, \ptre^\ell \, \eg_e^\ell \quad \mbox{ and } \quad \hg^{\ptre}_{e,k} := \sum_{\ell = 0}^{k} \, \ptre^\ell \, \hg_e^\ell
\end{equation*}
as (formal) approximations of order $k+1$ of $\eg^{\ptre}_e$ and
$\hg^{\ptre}_e$ respectively (notice that $k$ appears here as a subscript
while it appears as an exponent in the notation of the solution of the
approximate problem (\ref{eqinter}, \ref{GIBCcond})).  Using the ``second'' interface condition,
namely (\ref{interf-0}), one has 
\begin{equation} \label{bdry-gibc-0}
\eg^{\ptre}_{e,k}|_{\Gamma}(\xgg)  \times n=  \sum_{\ell = 0}^{k} \, \ptre^\ell \, \eg_i^\ell(\xgg, 0) \times n\;\; \mbox{for } \xgg \in \Gamma.
\end{equation} 
Substituting into (\ref{bdry-gibc-0}) the computed expressions of $\eg_i^\ell(\xgg, 0)$
 for $\ell = 1,2$ and $3$, respectively given by 
 (\ref{traceET_k=1}), (\ref{traceET_k=2}) and (\ref{traceET_k=3}))  leads to an identity  
of the form 
\begin{equation} \label{bdry-gibc-1}
\eg^{\ptre}_{e,k} \times n + \omega \, {\cal D}^{\ptre, k}\big[(
\hg^{\ptre}_{e,k})_T\big] = \ptre^{k+1} \, \varphi^\ptre_k \;\; \mbox{on } \Gamma ,\quad
\mbox{ for } k = 0,1,2,
\end{equation}
and
\begin{equation} \label{bdry-gibc-1-3}
\eg^{\ptre}_{e,3} \times n + \omega \, {\cal D}_0^{\ptre, 3}\big[(
\hg^{\ptre}_{e,3})_T\big] = \ptre^{4} \, \varphi^\ptre_{3,0} \;\; \mbox{on } \Gamma, 
\end{equation}
where ${\cal D}^{\ptre, k}$, $k=0,1,2$ are given by (\ref{Dk}) and ${\cal
  D}_0^{\ptre, 3}$ is given by (\ref{D3_0}) and where ${\varphi^\ptre_k} \in
  C^{\infty}(\Gamma)^3$, $k=0,1,2$ are tangential vector fields given by
\begin{equation} \label{lesphi}
\left\{ \begin{array}{ll}
\varphi^\ptre_0 &= 0, \\[12pt]
 \varphi^\ptre_1 &=  \sqrti \; \omega H_{e,T}^1, \\[12pt]
 \varphi^\ptre_2 &= \sqrti \; \omega \; H_{e,T}^2 + \omega ({\cal C}-{\cal H}) 
\big(H_{e,T}^1 + \ptre \; H_{e,T}^2 \big),
\end{array} \right.
\end{equation}
and  obviously satisfy the estimates (for $\ptre$ small enough)
\begin{equation} \label{estimphi}
\|\varphi_k^\ptre \|_{H^s_t(\Gamma)} \leq C_k(s), \quad k=0,1,2
\end{equation}
where $ C_k(s)$ is independent of $\ptre$, while $\varphi^\ptre_{3,0} \in C^{\infty}(\Gamma)^3$ is given by\\[4pt]
\begin{equation} \label{phi30}
\left\{ \begin{array}{ll}
 \varphi^\ptre_{3,0} & = \sqrti \; \omega 
\; H_{e,T}^3 + \omega ({\cal C}-{\cal H}) 
\big(H_{e,T}^3 + \ptre \; H_{e,T}^2 \big)\\[12pt]
& \dsp + \; \frac{1}{2 \sqrti} \;
 \big({\cal C}^2 - {\cal H}^2 \big) \; \big( H^1_{e,T} + \ptre \;  H^2_{e,T}
+  \ptre^3 \;  H^3_{e,T} \big)\\[12pt]
& \dsp + \;\frac{1}{2 \sqrti} \; \big( \,\eps_r \omega^2 + 
\grag \divg + \Rotg \rotg \big) \; \big( H^1_{e,T} + \ptre \;  H^2_{e,T}
+  \ptre^3 \;  H^3_{e,T} \big).
\end{array} \right.
\end{equation}
The GIBC (\ref{GIBCcond}) is obtained for $k=0,1,2$ by neglecting the
right-hand side of (\ref{bdry-gibc-1}). For $k = 3$, 
the same process leads to the condition (\ref{GIBC3}) that is
modified according to the process explained in Section~\ref{modif}.  
Notice that according to that construction, we have
\begin{equation} \label{identite3-bis}
\eg^{\ptre}_{e,3} \times n + \omega \, {\cal D}^{\ptre, k}\big[(
\hg^{\ptre}_{e,3})_T\big] = \ptre^{4} \, \varphi^\ptre_3 \;\; \mbox{on } \Gamma, \quad
\mbox{where } \varphi^\ptre_3 = \varphi^\ptre_{3,0} + \ptre \; {\cal R}^{\ptre,3} \big[(
\hg^{\ptre}_{e,3})_T\big],
\end{equation}
and using the property (\ref{boundR}) of ${\cal R}^{\ptre,3}$, 
\begin{equation} \label{estimphi3}
\|\varphi_3^\ptre \|_{H^s_t(\Gamma)} \leq C_3(s) \; 
\end{equation}
where $C_3(s) $ is independent of $\ptre$.
\subsection{Towards the theoretical justification of the GIBCs} 
\label{erroranalysis}
Our goal in the next two sections is to justify the GIBCs (\ref{GIBCcond})
by  estimating the errors
$$%\begin{equation} \label{deferror}
\Ee^\ptre - E_e^{\ptre,k} \quad \mbox{and} \quad ~\He^\ptre - H_e^{\ptre,k},
$$%\end{equation}
where $(E_e^{\ptre,k}, H_e^{\ptre,k})$ is the solution of the approximate problem
((\ref{eqinter}), (\ref{GIBCcond})), whose well-posedness will be
shown in Section~\ref{well-posedness} (see Theorem~\ref{exist-unicite2*}). It appears non trivial to work directly
with the differences $E_e^\ptre - E_e^{\ptre,k}$ and $H_e^\ptre -
H_e^{\ptre,k}$, we shall use the truncated series $(\eg^{\ptre}_{e,k}, \hg^{\ptre}_{e,k})$ introduced in Section
\ref{construc} as intermediate quantities. Therefore, the error
analysis is split into two steps:
\begin{enumerate}
\item Estimate the differences $\Ee^\ptre - \eg^{\ptre}_{e,k}
$ and $\He^\ptre - \hg^{\ptre}_{e,k} $ ; this is done in
Section \ref{errortrunk}, and more precisely in Lemma~\ref{estime-troncature_theorem} and Corollary~\ref{estime-troncature_theorem2}.
\item Estimate the difference
$\eg^{\ptre}_{e,k} - E_e^{\ptre, k}$ and $ \hg^{\ptre}_{e,k} -
H_e^{\ptre, k}$ ; this is done in Section \ref{errorGIBC} and
more precisely in Theorem \ref{error2}.
\end{enumerate}
\begin{remark}
Notice that step 1 of the proof is completely independent on the
$GIBC$ and will be valid for any integer $k$. Also, for $k=0$, the
second step is useless since $\widetilde E^{\ptre,0} =
E^{\ptre,0}$.
\end{remark}
\section{Error estimates for  the truncated expansions} \label{errortrunk}
\subsection{Main results} \label{main_results}
Let us introduce the fields $\Etke_\chi(x), ~ \Htke_\chi(x) : \Omega
\mapsto {\mathbb C}^3$ such that
\begin{equation*}
\Etke_\chi(x) \;\; = \;\; \left\{\begin{array}{ll} \dsp \sum_{\ell
= 0}^k \; \ptre^\ell \Ee^\ell(x) = \eg^{\ptre}_{e,k} \; , & \; \mbox{ for } x
\in \Oext, 
\\[12pt]
\chi(x) \dsp \sum_{\ell = 0}^k \; \ptre^\ell \Ei^\ell(\xgg,
\nu/\ptre) & \; \mbox{ for } x \in \Oint,
\end{array}\right.
\end{equation*}
\begin{equation*}
\Htke_\chi(x) \;\; = \;\; \left\{\begin{array}{ll} \dsp \sum_{\ell
= 0}^k \; \ptre^\ell \He^\ell(x) = \hg^{\ptre}_{e,k} \; , & \; \mbox{ for } x \in
\Oext, 
\\[12pt]
\chi(x) \dsp \sum_{\ell = 0}^k \; \ptre^\ell \Hi^\ell(\xgg,
\nu/\ptre) & \; \mbox{ for } x \in \Oint,
\end{array}\right.
\end{equation*}
where the local coordinates $\dsp \xgg$ and $\dsp \nu$ are defined as in
Section~\ref{preliminary} and the cut-off function $\chi$ is defined as
in section~\ref{asymptotic_ansatz}. These fields are good candidates to be good
approximations of the exact fields $(E^{\ptre}, H^{\ptre})$.  The main result of this section
is:
\begin{lemma}\label{estime-troncature_theorem}
For any integer $k$, there exists a constant $C_k$ independent of
$\ptre$ such that
\begin{equation} \label{lesestimations}
\left\{\begin{array}{clcl} (i) & \| E^\ptre - \Etke_\chi
\|_{H(\rt,\Omega)} &\le & C_k \; \ptre^{k+\frac{1}{2}},
\\[8pt]
(ii) & \| E^\ptre - \Etke_\chi \|_{L^2(\Oint)} &\le& C_k \;
\ptre^{k + \frac{3}{2}},
\\[8pt]
(iii) & \| E^\ptre  \times n - \Etke_\chi \times n
\|_{H^{-\frac{1}{2}}(\Gamma)} & \le& C_k \; \ptre^{k + 1}.
\end{array}\right.
\end{equation}
\end{lemma}

The proof of Lemma \ref{estime-troncature_theorem}, postponed to Section
\ref{Theproof},  rely on a fundamental a priori estimates  that we shall
state and prove in Section~\ref{techlem1}. We first give a straightforward
corollary of Lemma~\ref{estime-troncature_theorem}.

\begin{corollary}
\label{estime-troncature_theorem2}
For any integer $k$, there exists a constant $\widetilde C_k$ independent of
$\ptre$ such that:
\begin{equation*}\left \{ \begin{array}{lll}
\| \Ee^\ptre - \eg^{\ptre}_{e,k}\|_{H(\rt, \Omega_e)}\le \; \widetilde
C_k \; \ptre^{k+1}, \\[12pt]
\| \He^\ptre - \hg^{\ptre}_{e,k}\|_{H(\rt, \Omega_e)}\le \; \widetilde
C_k \; \ptre^{k+1}.
\end{array} \right.
\end{equation*}
\end{corollary}

\proof~ Simply write
$$
\Ee^\ptre - \eg^{\ptre}_{e,k} = \Ee^\ptre - \eg^{\ptre}_{e,k+1} + \ptre^{k+1}
\Ee^{k+1}
$$
which yields, since $\eg^{\ptre}_{e,k} = E_\chi^{\ptre,k+1}$ in
$\Oext$,
$$
\| \Ee^\ptre - \eg^{\ptre}_{e,k}\|_{H(\rt, \Omega_e)} \leq \|\Ee^\ptre -
E^{\ptre,k+1}_{\chi}\|_{H(\rt, \Omega_e)} +
\ptre^{k+1} \; \| E_e^{k+1}\|_{H(\rt,\Omega_e)}.
$$
Using the estimate  (\ref{lesestimations}-i) of Lemma~\ref{estime-troncature_theorem},
we get
$$
\| E_e^\ptre - \eg^{\ptre}_{e,k}\|_{H(\rt, \Omega_e)} \leq C_{k} \;
\ptre^{k+\frac{3}{2}} + \; \ptre^{k+1} \; \| \Ee^{k+1}\|_{H(\rt,
\Omega_e)} \leq \widetilde C_{k} \; \ptre^{k+1}.
$$
The estimates for $\He^\ptre - \hg^{\ptre}_{e,k}$ is an immediate consequence of
\begin{equation*}
\left\{ \begin{array}{llll} - i \omega (\He^\ptre - \hg^{\ptre}_{e,k}) + \rt
(\Ee^\ptre - \eg^{\ptre}_{e,k}) = 0 &\, \, \mbox{
in } ~~~ \Omega_e, \\[12pt]
i \omega (\Ee^\ptre - \eg^{\ptre}_{e,k})  + \rt (\He^\ptre - \hg^{\ptre}_{e,k}) =  0 &\, \,
\mbox{ in } ~~~ \Omega_e.
\end{array}\right.
\end{equation*}
\proofend

\subsection{A fundamental a priori estimate} \label{techlem1}
The proof of lemma \ref{estime-troncature_theorem} relies of the following fundamental technical lemma.
\begin{lemma} \label{stab-eps}
Assume that  $\bE^\ptre \in H(\rt, \Omega)
$   satisfies 
\begin{equation} \label{maxhom}
\left\{
\begin{array}{llll}
& \rt \, \rt \,\bE^\ptre - \omega^2 \bE^\ptre = 0,  &\;\; \mbox{ in } \;
\Oext, \\[4pt]
& i \omega \bE^\ptre_T  - \rt \, \bE^\ptre  \times n=  0 , & \;\; \mbox{ on } \;
\partial \Omega, 
\end{array}
\right.
\end{equation}
together with the following inequality
\begin{equation} \label{estime-stab}
 \begin{array}{cc} \dsp \left| \int_\Omega \left( |\rt \,
\bE^\ptre|^2 - \omega^2 |\bE^\ptre|^2 \; \right) ~ dx  + i \omega
\left( \int_{\partial \Omega} |\bE^\ptre \times n|^2 \; ds +
\frac{1}{\ptre^2} \int_{\Oint} |\bE^\ptre|^2 \; dx \right) \right| \\[12pt]
\leq \dsp A \left( \ptre^{s+\frac{1}{2}} \; \|\bE^\ptre \times n
\|_{H^{-\frac{1}{2}}(\Gamma)} +  \ptre^{s} \;
\|\bE^\ptre\|_{L^2(\Oint)} \right),
\end{array}
\end{equation}
for some non-negative constants $A$ and $s$ independent of
$\ptre$. Then there exists a constant $C$ independent of $\ptre$ such
that
\begin{equation} \label{estim_fond}
\| \bE^\ptre \|_{H( \rt, \Omega)} \le  C \,\ptre^{s+1},
\quad
\| \bE^\ptre \|_{L^2(\Oint)} \le C \,\ptre^{s+2},\quad
\| \bE^\ptre \times n \|_{H^{-\frac{1}{2}}(\Gamma)} \le C
\,\ptre^{s+\frac{3}{2}}, 
\end{equation}
for sufficiently small $\ptre$.
\end{lemma}
\proof~For convenience, we shall denote by $C$ a positive constant
whose value may change from one line to another but remains
independent of $\ptre$. We divide the proof in two steps.\\[12pt]
{\bf Step 1.} We first prove by contradiction that $\| \bE^\ptre \|_{L^2(\Omega)}\le  C \,\ptre^{s+1}$. This is the main step of the proof which will use two
important technical lemmas \ref{trace-lemma} and \ref{NewCompact-lemma}, that are proven
in the Appendix. \\[12pt] 
Assume that the positive quantity
$$\lambda^\ptre := \ptre^{-(s+1)} \; \|\bE^\ptre\|_{L^2(\Omega)}$$ is
unbounded as $\ptre \rightarrow 0$. After extraction of a subsequence, still
denoted $\bE^\ptre$ with $\delta \rightarrow 0$, we can assume that
$\lambda^\ptre \rightarrow + \infty.$
Let $\upE^\ptre =
\bE^\ptre/\|\bE^\ptre\|_{L^2(\Omega)}$ (so that $\|\upE^\ptre\|_{L^2(\Omega)} =1$). \\[12pt]
Our goal is to show that, up to the
extraction of another subsequence, $\upE^\ptre$ converges strongly
in $L^2(\Oext)$ and to obtain a contradiction by looking at the limit field
$\upE$.\\[12pt]
To show this, we wish to apply to $\upE^\ptre$ the compactness result of Lemma
\ref{NewCompact-lemma} with ${\cal O} = \Oext$. Since $\dv \,\bE^\ptre = 0$
and since $\upE^\ptre$ is bounded in $L^2(\Oext)$, we only
need to show that:
\begin{equation} \label{Amontrer}
\left| \begin{array}{ll}
\mbox{\bf (i)} &  \rt \, \upE^\ptre \mbox{ is bounded in } L^2 (\Oext),\\[12pt]
\mbox{\bf (ii)} & \upE^\ptre \times n|_{\partial \Omega} \mbox{ converges in }
H^{-\frac{1}{2}}(\partial \Omega), \\[12pt]
\mbox{\bf (iii)}& \upE^\ptre \times n|_{\Gamma} \mbox{ converges in }
H^{-\frac{1}{2}}(\Gamma).
\end{array} \right.
\end{equation}
We first notice that after division by
$\|\bE^\ptre\|_{L^2(\Omega)}$, the inequality (\ref{estime-stab})
yields
\begin{equation} \label{estime-stab-1}
\left| \quad \begin{array}{rr} \dsp \left| \int_\Omega \left( |\rt \,
\upE^\ptre|^2 - \omega^2 |\upE^\ptre|^2 \; \right)  dx  + i
\omega\left( \int_{\partial \Omega} |\upE^\ptre \times n|^2 \; ds
+ \frac{1}{\ptre^2} \int_{\Oint} |\upE^\ptre|^2 \; dx \right)
\right|
\\[12pt]
\dsp \le \quad \frac{A}{\lambda^\ptre} \left( \ptre^{-1} \;
\|\upE^\ptre\|_{L^2(\Oint)} + \ptre^{-\frac{1}{2}}
\;  \|\upE^\ptre \times n\|_{H^{-\frac{1}{2}}(\Gamma)} \right).
\end{array} \right.
\end{equation}
We shall now establish estimates on the two terms in the right hand
side of (\ref{estime-stab-1}) in terms of $\| \rt \, \upE^\ptre \|_{L^2(\Omega)}$
(namely inequalities (\ref{ineq6}) and (\ref{ineqtrace})). 
\\[12pt]
Considering the imaginary part of the left hand side of
(\ref{estime-stab-1}), we observe that since 
$1/\lambda^\ptre$ is bounded, 
$$
\| \upE^\ptre \|_{L^2(\Oint)}^2 \le \; C \; \ptre^{\frac{3}{2}}\;
  \|\upE^\ptre \times n \|_{H^{-\frac{1}{2}}(\Gamma)} +
C \; \ptre \; \|\upE^\ptre\|_{L^2(\Oint)}.
$$
Next, we use the trace inequality (\ref{trace-ineq}) of
Lemma~\ref{trace-lemma} with ${\cal O} = \Oint$ to get
$$
\| \upE^\ptre \|_{L^2(\Oint)}^2 \le \; C \; \ptre^{\frac{3}{2}} \;
\|\upE^\ptre\|_{L^2(\Oint)}^{\frac{1}{2}} \;
\left(\|\upE^\ptre\|_{L^2(\Oint)}^{\frac{1}{2}} +
\|\rt \, \upE^\ptre\|_{L^2(\Oint)}^{\frac{1}{2}} \right)
+ C \; \ptre \; \|\upE^\ptre\|_{L^2(\Oint)},
$$
which yields, after division by
$\|\upE^\ptre\|_{L^2(\Oint)}^{\frac{1}{2}}$,
\begin{equation} \label{ineq1}
\| \upE^\ptre \|_{L^2(\Oint)}^{\frac{3}{2}} \le C_1 \;  \ptre \;
\|\upE^\ptre\|_{L^2(\Oint)}^{\frac{1}{2}} +
C_2 \; \ptre^{\frac{3}{2}} \; \|\rt \,
\upE^\ptre\|_{L^2(\Oint)}^{\frac{1}{2}} .
\end{equation}
Let $K$ be a positive
constant to be fixed later.
Using Young's inequality $ ab \leq 2/3 \; a^{3/2}
+ 1/3 \; b^{3}$ with $a = K^{-1} \; \ptre$ and $b = K \;
\|\upE^\ptre\|_{L^2(\Oint)}^{\frac{1}{2}}\; $, \\we get
\begin{equation} \label{ineq2}
\ptre \; \|\upE^\ptre\|_{L^2(\Oint)}^{\frac{1}{2}} \le
\frac{2}{3} \; K^{-\frac{3}{2}} \; \ptre^{\frac{3}{2}} + \frac{K^3}{3}
\; \|\upE^\ptre \|_{L^2(\Oint)}^{\frac{3}{2}}.
\end{equation}
Choosing $ C_1 K^3 = {3}/{2}$ and substituting (\ref{ineq1}) into
(\ref{ineq2}), one deduces
\begin{equation}
\| \upE^\ptre \|_{L^2(\Oint)}^{\frac{3}{2}} \le C \; \ptre^{\frac{3}{2}}
\; \left( 1 + \|\rt \,\upE^\ptre\|_{L^2(\Oint)}^{\frac{1}{2}} \right),
\end{equation}
which yields
\begin{equation} \label{ineq6}
\ptre^{-1} \; \| \upE^\ptre \|_{L^2(\Oint)} \le C \; \left( 1 + \|\rt\,
\upE^\ptre\|_{L^2(\Oint)}^{\frac{1}{3}} \right).
\end{equation}
Now considering the real part of the left hand side of
(\ref{estime-stab-1}) 
and using the fact that 
$\| \upE^\ptre \|_{L^2(\Omega)} = 1$, we observe that
\begin{equation} \label{ineq4}
\| \rt \,\upE^\ptre \|_{L^2(\Omega)}^2 \le C\; \left( 1+
  \ptre^{-\frac{1}{2}} \; \|\upE^\ptre \times n \|_{H^{-\frac{1}{2}}(\Gamma)} +
  \ptre^{-1} \; \|\upE^\ptre\|_{L^2(\Oint)}\right).
\end{equation}
On the other hand, after multiplication by $\ptre^{-\frac{1}{2}}$, the
trace inequality (\ref{trace-ineq}) applied to $\upE^\ptre$ is equivalent to
$$
\ptre^{-\frac{1}{2}} \; \|\upE^\ptre \times n
\|_{H^{-\frac{1}{2}}(\Gamma)} \leq C \; \ptre^{\frac{1}{2}} \;
\left\{ \ptre^{-1} \; \| \upE^\ptre
  \|_{L^2(\Oint)}\right\}
+ C \; \left\{ \ptre^{-1}\| \upE^\ptre
  \|_{L^2(\Oint)}\right\}^{\frac{1}{2}} \; \|\rt\,
\upE^\ptre\|_{L^2(\Oint)}^{\frac{1}{2}}.
$$
After applying the Cauchy-Schwarz inequality to the second term of the right
hand side of the above inequality, we easily get, since  $\ptre$ is bounded
\begin{equation} \label{ineqtrace}
\left| \begin{array}{lll}
\ptre^{-\frac{1}{2}} \; \|\upE^\ptre \times n\|_{H^{-\frac{1}{2}}(\Gamma)}
 & \leq  & C \;  \left\{\ptre^{-1}\| \upE^\ptre\|_{L^2(\Oint)} + \; \|\rt\,
\upE^\ptre\|_{L^2(\Oint)} \right\}, \\[12pt]
& \leq  & C \;  \left\{ 1+ \|\rt\,
\upE^\ptre\|_{L^2(\Oint)}^{\frac{1}{3}} + \|\rt\,
\upE^\ptre\|_{L^2(\Oint)}\right\}, 
\end{array} \right.
\end{equation}
where we used (\ref{ineq6}) for the second inequality. Substituting
(\ref{ineqtrace}) into (\ref{ineq4}) shows that
\begin{equation*} %\label{estrot}
 \| \rt\,\upE^\ptre \|_{L^2(\Omega)}^2 \le C \; \left ( 1+ \|\rt\,
\upE^\ptre\|_{L^2(\Oint)}^{\frac{1}{3}} + \|\rt\,
\upE^\ptre\|_{L^2(\Oint)} \right).
\end{equation*}
This proves (\ref{Amontrer}-{\bf(i)}). We also deduce thanks to (\ref{ineq6}) and (\ref{ineqtrace}) that 
\begin{equation} \label{estl2trace}
\ptre^{-1}
||\upE^\ptre||_{L^2(\Oint)} \; \mbox{ and } \; \ptre^{-\frac{1}{2}} \;
||\upE^\ptre \times n ||_{H^{-\frac{1}{2}}(\Gamma)} \quad \mbox{are bounded,}
\end{equation}
which proves in particular (\ref{Amontrer}-{\bf(ii)}) ($\upE^\ptre \times n$
converges to 0 in $H^{-\frac{1}{2}}(\Gamma)$). This also means that the right
hand 
side of (\ref{estime-stab-1}) remains bounded. Thus, going
back to (\ref{estime-stab-1}) shows that $
||\upE^\ptre \times
n||_{L^2(\partial \Omega)}$ is bounded, 
which proves (\ref{Amontrer}-{\bf(iii)})  by the compactness of the $L^2(\partial
\Omega)$ embedding into $H^{-\frac{1}{2}}(\partial \Omega)$.\\[12pt]
Now we shall conclude the proof of Step 1. From  (\ref{Amontrer}-{\bf(i)}) one deduces that  $\upE^\ptre
$ is a bounded sequence in $ H(\rt, \Omega)$, therefore, up to extracted
subsequence, we can assume that  $\upE^\ptre$ weakly converges in $ H(\rt,
\Omega)$ to some $\upE$. Considering the restriction to $\Oext$, thanks to
(\ref{Amontrer}) we can apply the compactness result of Lemma \ref{NewCompact-lemma} and
deduce that an extracted subsequence of $\upE^\ptre
$, denoted again by $\upE^\ptre$ for simplicity, strongly converges to  $\upE$ in
$L^2(\Oext)$. On the other hand, we observe from (\ref{estl2trace}) that  $\upE^\ptre$
strongly converges  to $0$ in
$L^2(\Omega_i)$, hence  $\upE = 0$ in $\Oint$, which implies in particular
\begin{equation} \label{proplim}
 \upE \times n = 0
\mbox{ on } \Gamma. 
\end{equation}
Passing to the weak limit in equations (\ref{maxhom}) one easily verify that
\begin{eqnarray}\label{equation:4}
\left\{ \begin{array}{lll} \rt \,\, \rt \,\upE - \omega^2 \upE =0, ~~ &
\mbox{ in } \; \Omega_e,
\\[4pt]
i \omega \upE_T - \rt \, \upE \times n = 0, ~~ & \mbox{ on } \; \partial
 \Omega,
\end{array} \right.
\end{eqnarray}
The uniqueness of solutions to  (\ref{equation:4})-(\ref{proplim}) in $ H(\rt,
\Oext)$
implies that also $\upE =0 $ in $\Oext$. We therefore obtain that $\dsp \upE^\ptre$ converges to 0 in $L^2(\Omega)$
which is contradiction with $\dsp ||\upE^\ptre||_{L^2(\Omega)} = 1$.
Consequently $\lambda^{\ptre}$ is bounded, that is to say 
\begin{equation} \label{est_1}
\|\bE^\ptre \|_{L^2(\Omega)} \le C \; \ptre^{s+1}.
\end{equation}
{\bf Step 2.} We shall now proceed with the proof of estimates
(\ref{lesestimations}). Considering the imaginary part of the left hand side of
estimate (\ref{estime-stab}) and applying Lemma~\ref{trace-lemma} (with ${\cal
O} = \Oint$) yields 
\begin{equation*}%\label{est_2}
\| \bE^\ptre \|_{L^2(\Oint)}^2 \le C \left( \ptre^{s+\frac{5}{2}}
\;  \|\rt \,\bE^\ptre \|_{L^2(\Oint)}^{\frac{1}{2}}\; 
\|\bE^\ptre \|_{L^2(\Oint)}^{\frac{1}{2}} + \ptre^{s+2}\; \|\bE^\ptre
\|_{L^2(\Oint)}\right).
\end{equation*}
Using two times the Young inequality $ab \leq 1/2 (a^2
+ b^2)$, the first time with 
$$a =  \ptre^{\frac{1}{2}} \; \|\rt\, \bE^\ptre \|_{L^2(\Oint)}^{\frac{1}{2}} \quad
\mbox{ and } \quad b = \|\bE^\ptre \|_{L^2(\Oint)}^{\frac{1}{2}} \; ,$$
and the second time with
$$a = \|\bE^\ptre \|_{L^2(\Oint)} \quad \mbox{ and } \quad b = \ptre^{s+2},$$
leads to (we also use $\|\rt\, \bE^\ptre \|_{L^2(\Oint)} \leq \|\rt\, \bE^\ptre
\|_{L^2(\Omega)}\; $)
\begin{equation}\label{est_2bis}
\| \bE^\ptre \|_{L^2(\Oint)}^2 \le C \; \left(\ptre^{2s+4} + \ptre^{s+2}
\left(\| \bE^\ptre \|_{L^2(\Oint)} + \ptre \| \rt\, \bE^\ptre
  \|_{L^2(\Omega)}\right)\right).
\end{equation}
On the other hand, considering this the real part of the left hand side of
estimate (\ref{estime-stab}) and using (\ref{est_1}), we get
\begin{equation*}%\label{est_3}
\| \rt \,\bE^\ptre \|_{L^2(\Omega)}^2 \le C \; \left(\ptre^{2s+2} +
\ptre^{s+\frac{1}{2}} \; \|\rt \,\bE^\ptre \|_{L^2(\Oint)}^{\frac{1}{2}}
\; \|\bE^\ptre \|_{L^2(\Oint)}^{\frac{1}{2}} + \ptre^{s} \; \|\bE^\ptre
\|_{L^2(\Oint)}
\right),
\end{equation*}
which gives, using Young's inequality once again, 
\begin{equation}\label{est_3bis}
\| \rt\, \bE^\ptre \|_{L^2(\Omega)}^2 \le  C \; \left(\ptre^{2s+2} + \ptre^{s}\;
\left(\| \bE^\ptre \|_{L^2(\Oint)} + \ptre \; \| \rt\, \bE^\ptre
  \|_{L^2(\Omega)}\right)\right).
\end{equation}
Combining (\ref{est_2bis}) and  (\ref{est_3bis}) leads to
$$
\| \bE^\ptre \|_{L^2(\Oint)}^2 + \ptre^2 \; \| \rt\, \bE^\ptre
  \|_{L^2(\Omega)}^2 \le C \left(\ptre^{2s+4} + \ptre^{s+2} \;
\left(\| \bE^\ptre \|_{L^2(\Oint)} + \ptre \; \| \rt \,\bE^\ptre
  \|_{L^2(\Omega)}\right)\right),
$$
which yields 
$$
\| \bE^\ptre \|_{L^2(\Oint)} + \ptre \; \| \rt \,\bE^\ptre
  \|_{L^2(\Omega)} \le C \; \ptre^{s+2}
$$
and in particular the second inequality of (\ref{estim_fond}). The third inequality of (\ref{estim_fond}) is a direct consequence of
the first two ones and the application of Lemma~\ref{trace-lemma} in
$\Oint$.\proofend 
\begin{remark} Notice that since we simply used in the first step of the proof
  the fact that
$1/\lambda^{\ptre}$ is bounded, we have proved in fact  that
$$
\lim_{\ptre \rightarrow 0} \; \ptre^{-(s+1)} \; \|E^\ptre \|_{L^2(\Omega)}= 0.
$$
\end{remark}
\subsection{The proof of Lemma \ref{estime-troncature_theorem}.} 
\label{Theproof}
Let us introduce, for each integer $k$, the {\it error fields}
\begin{equation}
\bfdEk= \Ee^\ptre - \Etke_{\chi}, \quad \dHke= \He^\ptre - \Htke_{\chi}.
\end{equation} 
The
idea of the proof is to show that $\dEke$ satisfies an a priori
estimate of the type (\ref{estime-stab}) and then to use the
stability lemma~\ref{stab-eps}. To prove such an estimate, we
shall use the equations satisfied by $(\dEke, \dHke)$,
respectively in $\Oint$ and $\Oext$. \\[12pt]
\noindent{\bf The equations in $\Oext$.}
It is straightforward to check that in the exterior domain $\Oext$, the
errors
$$%\begin{equation} \label{errors-ext}
(\bfdEke, \bfdHke) := (\bfdEk|_{\Oext},\bfdHk|_{\Oext})
$$%\end{equation}
satisfies the homogeneous equation:
\begin{eqnarray} \label{errorOextpre1}
\left\{ \begin{array}{lll}
(i)& \rt \, \bfdHke + i \omega \bfdEke =0, &\;\;\; \textrm{ in } \Omega_e, \\[12pt]
(ii)&  \rt \, \bfdEke - i \omega \bfdHke=0, &\;\;\; \textrm{ in }
\Omega_e,\\
\end{array} \right.
\end{eqnarray}
and
\begin{equation}\label{errorboundary}
({\bfdEke})_T  - \bfdHke \times n = 0 , \;\;\; \textrm{ on }
\partial \Omega.\\[12pt]
\end{equation}
\noindent Eliminating $\bfdHke$ in (\ref{errorOextpre1}), we get 
\begin{equation}\label{errorOext}
\left \{ \begin{array}{ll}
\rt \, \big( \rt \, \bfdEke \big) - \omega^2 \bfdEke = 0, & \;\;\; \textrm{ in }
\Omega_e.\\[12pt]
\rt \, \bfdEke \times n + i \omega \big(\bfdEke\big)_T = 0, & \;\;\; \textrm{ on }
\partial \Omega.
\end{array} \right. \end{equation}
{\bf The equations in $\Oint$.} Now consider the restrictions to $\Oint$ and
set 
$$%\begin{equation} \label{errors-ixt}
(\bfdEki, \bfdHki) := (\bfdEk|_{\Oint},\bfdHk|_{\Oint}).
$$%\end{equation}
It is also useful to introduce the fields 
\begin{equation*}% \label{notation2}
E^\ptre_{i,k} (\xgg, \nu) := \sum_{\ell = 0}^{k} \ptre^{\ell} \;
E_{i}^{\ell}(\xgg, \frac{\nu}{\ptre}), \quad
H^\ptre_{i,k}(\xgg, \eta) := \sum_{\ell = 0}^{k} \ptre^{\ell} \; H_{i}^{\ell}(\xgg, \frac{\nu}{\ptre}),
\end{equation*}
so that using the local coordinates, we can write
\begin{equation*}
\Etke_{\chi}(x) = \chi \; E^\ptre_{i,k}(\xgg, \eta), \quad 
\Htke_{\chi}(x) = \chi \;
H^\ptre_{i,k}(\xgg, \eta) \quad 
\mbox{in } \Oint.
\end{equation*}
Our goal is to show that $(\Etke_{\chi}, \Htke_{\chi})$ satisfy the 
``interior equations'' except that two {\it small} source terms appear at the
right hand side, respectively due to the cut-off function $\chi$ and the
truncation of the series at order $k$. We first compute that
\begin{equation}\label{formulecompose}
\left\{ \begin{array}{lll} \dsp  \rt
\Htke_{\chi} + i \omega \Etke_{\chi} - \frac{1}{\omega \, \ptre^2} \Etke_{\chi} &= \dsp \chi \left (
i \omega E^\ptre_{i,k} + \rt H^\ptre_{i,k} -
  \frac{1}{\omega \, \ptre^2} E^\ptre_{i,k} \right) + \nabla \chi \times H^\ptre_{i,k}, \\[12pt]
\rt \Etke_{\chi} - i \omega \Htke_{\chi}  &= \chi \left (\rt E^\ptre_{i,k}  - i
\omega H^\ptre_{i,k}\right) + \nabla \chi \times E^\ptre_{i,k}.
\end{array} \right.
\end{equation}
Thanks to the exponentially decaying nature of $E^\ptre_{i,k}(\xgg,\eta)$ 
and $H^\ptre_{i,k}(\xgg,\eta)$ with respect to $\eta$ (cf. Theorem~\ref{lesEkHk}), 
the terms in factor of $\nabla \chi$ are exponentially small in
$\ptre$.\\[12pt] 
It remains to compute the terms in factor of $\chi$. These calculations are
tedious, but the idea is simple and consists - in some sense - to do the same
calculations as in Section~\ref{Asymptotic formal matching} but in the reverse
sense. According to (\ref{scaledE}) we define 
\begin{equation*}% \label{Escale}
\mE^\ptre_{i,k} := \frac{E^\ptre_{i,k}}{\delta} = \sum_{p=0}^{k-1} \ptre^p \; \mE_i^p.
\end{equation*}
With the notation of Section~\ref{Asymptotic formal matching}, we have 
\begin{equation} \label{utile}
\rt E^\ptre_{i,k}  - i
\omega H^\ptre_{i,k} = r_{i}^{\delta,k}(\xgg , \nu /\ptre) 
\end{equation}
where the function $r_{i}^{\delta,k}(x ,\eta)$ is given by
$$
r_{i}^{\delta,k} = \partial_\eta \mE^\ptre_{i,k} \times n -  i \omega
\, H^\ptre_{i,k} + \sum_{\ell=1}^{2} \; \ptre^{\ell} A_{E}^{(\ell)}
(\mE^\ptre_{i,k},H^\ptre_{i,k}).
$$
Replacing $\mE^\ptre_{i,k}$ and $H^\ptre_{i,k}$ by their polynomial expansion in $\delta$, we get
$$
r_{i}^{\delta,k} =   \sum_{p=0}^{k-1} \ptre^{p} \; 
\big( \partial_{\eta} \mE_i^p \times n -  i \omega \,H_i^p \big)
-  i \omega \, \ptre^{k} \; H_i^k + 
\sum_{\ell=1}^{2} \; \ptre^{\ell} \; \sum_{p=0}^{k-1}
\ptre^{p} \; A_{E}^{(\ell)} (\mE_i^p, H_i^p).
$$
Using the equations (\ref{induction-ehh}) satisfied by the $\mE_i^p$'s
and $H_i^p$'s, we get
$$
r_{i}^{\delta,k} =  -  i \omega \, \ptre^{k} \; H_i^k + 
\sum_{\ell=1}^{2} \; \sum_{p=0}^{k-1}
\ptre^{p+\ell} \; A_{E}^{(\ell)} (\mE_i^p, H_i^p) - \sum_{\ell=1}^{2} \sum_{p=0}^{k-1} \ptre^{p}
\; A_{E}^{(\ell)} (\mE_i^{p-\ell}, H_i^{p-\ell}) .
$$
Applying the change of index $p+l \rightarrow p$ in the first sum, we get
$$
r_{i}^{\delta,k} =  -  i \omega \, \ptre^{k} \; H_i^k + 
\sum_{\ell=1}^{2} \; \sum_{p=0}^{k-1+\ell}
\ptre^{p} \; A_{E}^{(\ell)} (\mE_i^{p-\ell}, H_i^{p-\ell}) - \sum_{\ell=1}^{2} \sum_{p=0}^{k-1} \ptre^{p}
\; A_{E}^{(\ell)} (\mE_i^{p-\ell}, H_i^{p-\ell}),
$$
that is to say
$$
r_{i}^{\delta,k} =
-  i \omega \, \ptre^{k} \; H_i^k + 
\sum_{\ell=1}^{2} \; \sum_{p=k}^{k-1+\ell}
\ptre^{p} \; A_{E}^{(\ell)} (\mE_i^{p-\ell}, H_i^{p-\ell}) .
$$
Paying attention to the above expression and using the form of the functions
$\mE_i^p$ and $H_i^p$ (cf. Theorem \ref{lesEkHk}), we see that
$$%\begin{equation} \label{res2}
\rt \, \mE^\ptre_{i,k}  - i \omega H^\ptre_{i,k} =  \ptre^{k} \; \big( g_{k,0}^{\ptre} + \ptre
\; g_{k,1}^{\ptre}) \quad \mbox{in } \supp \chi,
$$%\end{equation}
where the functions $g_0^{\ptre}$ and $g_1^{\ptre}$ are of the form
\begin{equation} \label{formgq}
\dsp g_{k,q}^{\ptre}(x) = p_{k,q}(\xgg, \frac{\nu}{\delta}) \; e^{-\sqrt{i} \,
  \frac{\nu}{\delta}}, \quad  p_{k,q} \in {\bf P}_k(\Gamma, \R^+; \C^3), \quad
q=0,1 .
\end{equation} 
From (\ref{formgq}), we easily deduce that
\begin{equation} \label{propgq}
\|\chi \; g_{k,q}^{\ptre}\|_{L^2(\Oint)} \leq C_{k,q} \; \ptre^{\frac{1}{2}},
  \quad \|\chi \; \rt \, g_{k,q}^{\ptre}\|_{L^2(\Oint)} \leq C'_{k,q} \; \ptre^{\frac{1}{2}},
  \quad q = 0,1.
\end{equation} 
In the same way, using again local coordinates, we have 
$$
i \omega E^\ptre_{i,k} + \rt H^\ptre_{i,k} -
  \frac{1}{\ptre^2} E^\ptre_{i,k} = 
\frac{1}{\delta} \; s_{i}^{\delta,k}(x , \nu /\ptre), 
$$
with
$$
s_{i}^{\delta,k} = \partial_\eta H^\ptre_{i,k}  \times n 
+ \frac{1}{\omega} \; \mE^\ptre_{i,k}
- \dsp \sum_{\ell=1}^{4} \; \ptre^{\ell} \;
A_H^{(\ell)} (\mE^\ptre_{i,k}, H^\ptre_{i,k} ) 
$$
Replacing $\mE^\ptre_{i,k}$ and $H^\ptre_{i,k}$ by their polynomial expansion
in $\delta$ we get
$$
\left| \begin{array}{lll}
s_{i}^{\delta,k} & = & \dsp \sum_{p=0}^{k-1} \; \ptre^{p} \; 
\dsp \big( \partial_{\eta} H_i^p \times n + \frac{1}{\omega} \,\mE_i^p \big)
- \dsp \sum_{\ell=1}^{4} \; \ptre^{\ell} \; \sum_{p=0}^{k-1} \ptre^{p} \;
A_H^{(\ell)} (\mE_i^p, H_i^p ) \\[12pt]
& + & \dsp \ptre^{k} \; \partial_{\eta} H_i^p \times n - \sum_{\ell=1}^{4} \;
\ptre^{\ell + k}  \; A_H^{(\ell)} (0, H_i^k).
\end{array} \right.
$$
Using equations (\ref{induction-ehh}) satisfied by the $\mE_i^p$'s
and $H_i^p$'s, we get 
$$
\left| \begin{array}{lll}
s_{i}^{\delta,k} & = & \dsp \sum_{\ell=1}^{4} \; \sum_{p=0}^{k-1} \; \ptre^{p} \; \dsp A_H^{(\ell)} (\mE_i^{p-\ell}, H_i^{p-\ell} )
- \dsp \sum_{\ell=1}^{4} \; \sum_{p=0}^{k-1} \; \ptre^{p+\ell} \;
A_H^{(\ell)} (\mE_i^p, H_i^p ) \\[12pt]
& + & \dsp \ptre^{k} \; \partial_{\eta} H_i^p \times n - \sum_{\ell=1}^{4} \;
\ptre^{\ell + k}  \; A_H^{(\ell)} (0, H_i^k),
\end{array} \right.
$$
or equivalently
$$
\left| \begin{array}{lll}
s_{i}^{\delta,k} & = & \dsp \sum_{p=k}^{k+\ell-1} \; \ptre^{p} \; 
\dsp \sum_{\ell=1}^{4} A_H^{(\ell)} (\mE_i^{p-\ell}, H_i^{p-\ell} )
\\[12pt]
& + & \dsp \ptre^{k} \; \partial_{\eta} H_i^p \times n - \sum_{\ell=1}^{4} \;
\ptre^{\ell + k}  \; A_H^{(\ell)} (0, H_i^k)
\end{array} \right.
$$
This time, we see that we can write
$$
i \omega E^\ptre_{i,k} + \rt H^\ptre_{i,k} -
  \frac{1}{\ptre^2} E^\ptre_{i,k} = \frac{1}{\delta} \; \sum_{q=0}^3 \; \ptre^q \; 
h_{k,q} ^{\ptre} \quad \mbox{in } \supp \chi,
$$
where the expression of $h_{k,q}^{\ptre}$ is similar to the
$g_{k,q}$'s (see formula  (\ref{formgq})) and implies in particular that
\begin{equation} \label{prophq}
\|\chi \; h_{k,q}^{\ptre}\|_{L^2(\Oint)} \leq C_{k,q} \; \ptre^{\frac{1}{2}},
  \quad q = 0,1,2,3.
\end{equation} 
In summary,  taking the difference between  (\ref{formulecompose}) and
(\ref{trans_eq_2})  we have shown that
$$%\begin{equation}\label{eq-error}
\left\{ \begin{array}{lll} \dsp  \rt \,
\bfdHki + i \omega \bfdEki - \frac{1}{\omega \ptre^2} \; \bfdEki&= \dsp
\delta^{k-1} \; \chi \; \Big(\sum_{q=0}^3 \; \ptre^q \; h_{k,q}^{\ptre} \Big)
+ \nabla \chi \times H^\ptre_{i,k}, \\[12pt]
\dsp 
\rt \, \bfdEki - i \omega \bfdHki  &= \dsp \delta^k \; \chi \;
\Big(\sum_{q=0}^1 \; \ptre^q \; g_{k,q}^{\ptre} \Big) +
\nabla \chi \times E^\ptre_{i,k},
\end{array} \right.
$$%\end{equation}
where eliminating $\bfdHki$ we get 
\begin{equation}\label{electrique}
\rt \; \rt \, \bfdEki - \omega^2 \; \bfdEki +
\frac{i}{\delta^2} \; \bfdEki = f_k^{\ptre}
\end{equation}
with 
$$%\begin{equation}\label{second-membre}
\left| \begin{array}{lll} 
f_k^{\ptre} & := & \dsp \ptre^k \; \chi \; \Big(\sum_{q=0}^1 \; \ptre^q \; \rt \,
g_{k,q}^{\ptre} 
\Big) + \nabla \chi \times \Big(\sum_{q=0}^1 \; \ptre^q \; g_{k,q}^{\ptre}
\Big) + \nabla \chi \times \rt \, \big( \nabla \chi \times E^\ptre_{i,k} \big)
\\[12pt] & - & \dsp i \omega \; \delta^{k-1} \; \chi \; \Big(\sum_{q=0}^3 \;
\ptre^q \; h_{k,q}^{\ptre} \Big) \; - \;  i \omega \; \nabla \chi \times
H^\ptre_{i,k}. \end{array} \right.
$$%\end{equation}
Taking into account the form of
the functions $g_{k,q}^{\ptre}$ and the
exponential decay of the the fields $E_i^p$ and  $H_i^p$ (Theorem
\ref{lesEkHk}), and since the support of $\nabla \chi$
is separated from $\Gamma$, there exists a constant $\tau > 0$ such that:
\begin{equation*}%\label{estimnabla}
\left| \begin{array}{l}
\| \, \nabla \chi \times \rt \, \big( \nabla \chi \times E^\ptre_{i,k}
\big)\,\|_{L^2(\Oext)} \; \leq \; C_1(k) \; e^{- \tau \, \ptre},\\[12pt]
\dsp \|\,\nabla \chi \times \Big(\sum_{q=0}^1 \; \ptre^q \; g_{k,q}^{\ptre}
\Big)\,\|_{L^2(\Oext)}  \; \leq \; C_2(k) \; e^{- \tau \, \ptre},\\[18pt]
\|\,\nabla \chi \times
H^\ptre_{i,k} \,\|_{L^2(\Oext)}  \leq C_3(k) \; e^{- \tau \, \ptre}.
\end{array} \right.
\end{equation*}
Combining these inequalities with estimates (\ref{propgq}) and
(\ref{prophq}), we see that:
\begin{equation}\label{estimf}
\|f_k^{\ptre}\|_{L^2(\Oext)} \; \leq \; C_k \; \ptre^{k - \frac{1}{2}} \;
. 
\end{equation}
{\bf Error estimates.} We can now proceed with the final step of the
proof.
First, we multiply the equation  (\ref{errorOext}) by $\overline{\bfdEke}$ and
integrate over $\Oext$. Using the Stokes formula and the boundary condition in
(\ref{errorOext}), we get 
\begin{equation*}%\label{estExt}
\dsp \int_{\Omega_e} | \rt \, \bfdEke |^2 \;
dx - \omega ^2 \int_{\Omega _e} |\bfdEke|^2 \; dx \dsp - i \omega
\int_{\partial \Omega} | \bfdEke \times n |^2
\dsp  + \left<\rt \, \bfdEke \times n , \;
\big(\overline{\bfdEke}\big)_T \right>_{\Gamma}= 0.
\end{equation*}
Next, we multiply the equation  (\ref{electrique}) by $\bfdEki$ and integrate
over $\Oint$. We get 
\begin{equation*}%\label{estInt}
\left| \begin{array}{l}
\dsp \int_{\Omega_i} | \rt \bfdEki |^2 \;
dx - \omega ^2 \int_{\Omega_i} |\bfdEki|^2 \; dx  -
\frac{1}{\ptre^2} \int_{\Omega _i} |\bfdEki|^2
\dsp  - \left<\big(\rt \bfdEke \times n\big) \cdot
\big(\overline{\bfdEki}\big)_T \right>_{\Gamma}  \\[12pt]
\quad \quad \quad \quad \quad \quad
\quad \quad \quad \quad \quad \quad
\quad \quad \quad \quad \quad \quad
\quad \quad \quad \quad \quad 
\quad \quad \quad \quad \quad = \dsp
\int_{\Omega_i}  f_k^{\ptre} \cdot \overline{\bfdEki} \; dx .
\end{array} \right.
\end{equation*}
Adding the last two equalities and using the fact that $\bfdEk$ belongs to
$H(\rt ; \Omega)$, we get
\begin{equation} \label{presquefini}
\left| \begin{array}{l}
\dsp \int_{\Omega} | \rt \, \bfdEk |^2 - \omega ^2
\int_{\Omega} |\bfdEk|^2 - i \omega \Big( \int _{\partial \Omega }
| \dEke \times n |^2 + \frac{1}{\ptre^2} \int_{\Omega_i} |\dEke|^2 \Big)
\\[12pt]
\quad \quad \quad \quad \quad \quad \dsp
= \left<\rt \, \bfdEke \times n - \rt \, \bfdEki \times n, 
\big(\overline{\bfdEk}\big)_T \right>_{\Gamma} + \int_{\Omega_i}  f_k^{\ptre} \cdot \overline{\bfdEki} \; dx
\end{array} \right.
\end{equation}
It remains to compute the jump
$$
\rt \, \bfdEke \times n - \rt \, \bfdEki \times n \equiv \rt \, E^\ptre_{e,k} \times n - \rt \, E^\ptre_{i,k} \times n
$$
across $\Gamma$. Taking the trace on $\Gamma$ of equation (\ref{utile}), we
get, with $\rho_{i}^{\delta,k}(\xgg) = r_{i}^{\delta,k}(\xgg , 0)$,
$$ 
\rt E^\ptre_{i,k} \times n =  i \omega H^\ptre_{i,k} \times n +
\rho_{i}^{\delta,k} \quad \mbox{on } \Gamma.
$$
The function $\rho_{i}^{\delta,k}$ is not zero but small. In particular,
according to , we have 
\begin{equation} \label{estimrho}
\|\rho_{i}^{\delta,k}\|_{H^{\frac{1}{2}}(\Gamma)} \leq C_k \; \delta^k.
\end{equation}
On the other hand, taking the trace on $\Gamma$ of the first equation of
(\ref{errorOextpre1}) we get
$$ 
\rt E^\ptre_{e,k} \times n =  i \omega H^\ptre_{e,k} \times n, 
\quad \mbox{on } \Gamma
$$
The continuity conditions (\ref{cond_infini_limites_EH}) imply $H^\ptre_{i,k}
\times n = H^\ptre_{e,k} \times n \mbox{ on } \Gamma$ so that
\begin{equation} \label{expsaut}
\rt \, \bfdEke \times n - \rt \, \bfdEki \times n \equiv \rt \, E^\ptre_{e,k}
\times n - \rt \, E^\ptre_{i,k} \times n = \rho_{i}^{\delta,k}.
\end{equation}
Substituting (\ref{expsaut}) into (\ref{presquefini}) and using the estimates
(\ref{estimf}) and (\ref{estimrho}), we get 
\begin{equation*}% \label{presquepresquefini}
\left| \begin{array}{l}
\dsp \left| \int_{\Omega} | \rt \bfdEk |^2 - \omega ^2
\int_{\Omega} |\bfdEk|^2 - i \omega \Big( \int _{\partial \Omega }
| \dEke \times n |^2 + \frac{1}{\ptre^2} \int_{\Omega_i} |\dEke|^2 \Big) \right|
\\[12pt]
\quad \quad \quad \quad \quad \quad \quad 
\; \leq \; C_k \; \big( \; \delta^{k-\frac{1}{2}} \;
\| \bfdEk \times n \|_{L^2(\Oint)} + \;
\delta^k \; \| \bfdEk \times n \|_{H^{-\frac{1}{2}}(\Gamma)} \; \big) \; . 
\end{array} \right.
\end{equation*}
We can finally apply Lemma~\ref{stab-eps} with ${\bf E}^{\ptre} = \bfdEk$ and $s =
k - \frac{1}{2}$, which provides the desired estimates.  $\square$
%\proofend
\section{Analysis of the GIBCs} \label{Sec-Analysis}
\subsection{Well-posedness of the approximate problems}
\label{well-posedness} We shall prove in
this section that the approximate fields  
$(E^{\ptre,k}, H^{\ptre,k})$, solution of (\ref{eqinter}, \ref{GIBCcond}) 
for $k=$ 0, 1, 2, 3, are well defined. In fact, for $k \leq 2$  this result is an application (or an adaptation) of
classical results about Maxwell equations with an impedance boundary
condition of the form
$$%\begin{equation}\label{impedance}
E \times n + \omega \,Z \; H_T = 0 \mbox{
on } \Gamma,
$$%\end{equation}
where $Z$ is a function with positive real part (see for instance \cite{Monk}). To include the case $k=3$ it
is sufficient to extend these results to the cases where $Z$ is a continuous operator
form $L^2_t(\Gamma)$ into $L^2_t(\Gamma)$ with positive definite real
part. More precisely we shall assume that there exists two positive constants
$z_*$ and $z^*$ such that
\begin{equation} \label{prop-cd}
\left\{ \begin{array}{ll}
(i) & \|Z\varphi\|_{\Gamma} \leq z^* \; \|\varphi\|_{\Gamma}, \\[6pt]
(ii) & 
\Ree (Z\varphi, \varphi)_\Gamma \geq  z_* \; \|\varphi\|_{\Gamma}^2,
\end{array} \right.
\end{equation}
for all $\varphi \in L^2_t(\Gamma)$.
These properties are satisfied by the operators ${\cal D}^{\ptre, k}$, $k=1, 2,
3$ for $\ptre$ sufficiently small, and can be seen as a special consequence of Lemma~\ref{lembuy} (stated and
proved in next section) where the dependence of the constants $z_*$ and $z^*$
in terms of $\ptre$ is also given (which is important for error analysis). The functional space adapted
to this type of boundary conditions is the same as for constant impedances,
namely $\hcurlt{\Oext}$ (see (\ref{defhcurlt}) for the definition of this space).
% \begin{equation} \label{espaceX}
% X = \{ H \in H(\rt, \Oext) \; /\; {H_T}|_{\Gamma} \in L^2_t(\partial \Omega),
% {H_T}|_{\Gamma} \in L^2_t(\Gamma) \; \} \quad (\equiv {\cal V}_H^3), 
% \end{equation}
% which is an Hilbert space when equipped with the norm
% \begin{equation} \label{normX}
% \|H\|^2_X = \|H\|_{L^2(\Oext)}^2 + \|\, \rt \, H\|_{L^2(\Oext)}^2 + \|
% H_T\|_{\Gamma}^2 + \|H_T\|_{\partial \Omega}^2.
% \end{equation}

\begin{theorem}
\label{exist-unicite2*} Let $f\in L^2(\Oext)$ be compactly supported in
$\Oext$ and $g \in L^2_t(\partial \Oext)$. Then the boundary value problem

\begin{equation*}%\label{GIBC-equation2*}
\left \{\begin{array}{llll}
\rt \, H + i \omega E= 0,  &\mbox{ in } \Oext,\\[6pt]
\rt \, E - i \omega H  = 0,  &\mbox{ in } \Oext,\\[6pt]
E\times n  + H_T  = g, &\mbox{
on } \partial \Omega,\\[6pt]
E \times n + \omega Z H_T = g , &\mbox{
on } \Gamma.
\end{array} \right.
\end{equation*}
has a unique solution $(E,H)$ in $\hcurlt{\Oext}\times \hcurlt{\Oext}$.
\end{theorem}

\proof
The proof  uses basically the same
arguments as for classical impedance conditions (see  for instance \cite{Monk})
and is detailed here for the reader convenience.  The approach
is divided in three steps: 
\begin{enumerate}
\item One eliminates the electric field $E$ and formulate a boundary value problem
  in $H$ only, then writes the associated weak formulation  in an appropriate
  functional framework.
\item One shows that the Fredholm alternative can be applied for solving
  this problem.
\item One shows the uniqueness of the solution which also implies the
  existence.
\end{enumerate}
{\bf Step 1.} The problem to be solved for $H\in \hcurlt{\Oext}$ is
\begin{equation} \label{pbH}
\left\{ \begin{array}{ll}
\rt  \, \rt \, H - \omega^2 \; H = - \rt \, f & \mbox{ in } \; \Oext,
\\[6pt]
\rt \, H \times n - i \; \omega^2 Z H_T = -i\omega g, 
& \mbox{ on } \; \Gamma, \\[6pt]
\rt \, H \times n - i \; \omega \; H_T = -i\omega g, 
& \mbox{ on } \; \partial \Omega.
\end{array} \right.
\end{equation}
The equivalent weak formulation of the
boundary value problem (\ref{pbH}) can be written in the form:
\renewcommand{\wH}{H'}
\begin{equation} \label{weakform}
\left\{ \begin{array}{l}
\mbox{ Find } H \in \hcurlt{\Oext} \mbox{ such that}, \\[6pt]
\dsp \int_{\Oext} \big( \, \rt \, H \cdot \overline{\rt \, \wH} - \omega^2 H \cdot \overline{\wH} \,
\big) \; dx  
\; - \; i \; \omega (H_T, \wH_T)_{\partial \Omega} - \; i \; \omega^2 \; 
(ZH_T, \wH_T)_\Gamma\\[12pt] 
\quad \quad \quad \quad  \dsp = - \; i \; \omega (g, \wH_T)_{\partial \Oext} + \int_{\Oext} f \cdot \overline{\rt \, \wH} \; dx , \quad \forall
\; \wH \in \hcurlt{\Oext}. 
\end{array} \right.
\end{equation}
Next, we formulate a variational problem, equivalent to (\ref{weakform}),
posed in a subspace $\hcurltz{\Oext}$ of $\hcurlt{\Oext}$ having additional compactness properties, by
using the Helmholtz decomposition. Let us consider the
closed subspace of $H^1(\Omega)$, 
\begin{equation*}% \label{defS}
S := \{ v \in H^1(\Omega) \; /\; (u|_{\Gamma},u|_{\partial \Omega}) \in
P_0(\Gamma) \times P_0(\partial \Omega) \}  
\end{equation*}
where $P_0(\Gamma)$ is the space of constant
functions on $\Gamma$ (the same for $P_0(\partial \Omega$)). Then we define
$\hcurltz{\Oext}$ as
\begin{equation*}% \label{defX0}
\hcurltz{\Oext} := \{ H \in \hcurlt{\Oext} \; /\; \int_{\Oext} H \cdot \nabla v \; dx = 0, \forall \; v
\in S \},
\end{equation*}
which forms a closed subset of $\hcurlt{\Oext}$. One also has the
orthogonal  decomposition (in $L^2(\Oext)$)
\begin{equation} \label{decompX}
\hcurlt{\Oext} = \hcurltz{\Oext} \oplus \nabla S, \quad \mbox{where } \nabla S := \{ \nabla v \; /\; v
\in S \}.
\end{equation}
\begin{remark}
By Green's theorem, the reader will notice that $\hcurltz{\Oext}$ is nothing but the
subspace of $\hcurlt{\Oext}$ made of vector fields whose divergence vanishes and whose
normal trace on each connected component of $\partial \Oext$, namely $\Gamma$
and $\partial \Omega$, has zero mean value.
\end{remark}
We claim that the problem (\ref{weakform}) is equivalent to
\begin{equation} \label{weakformbis}
\left\{ \begin{array}{l}
\mbox{ Find } H \in \hcurltz{\Oext} \mbox{ such that}, \\[6pt]
\dsp \int_{\Oext} \big( \, \rt \, H \cdot \overline{\rt \, \wH} - \omega^2 H \cdot \overline{\wH} \,
\big) \; dx  
\; - \; i \; \omega (H_T, \wH_T)_{\partial \Omega}- \; i \; \omega^2 \; 
(ZH_T, \wH_T)_{\Gamma}\\[12pt] 
\quad \quad \quad \quad  \dsp  = - \; i \; \omega (g, \wH_T)_{\partial \Oext}+\int_{\Oext} f \cdot \overline{\rt \, \wH} \; dx , \quad \forall
\; \wH \in \hcurltz{\Oext}. 
\end{array} \right.
\end{equation}
Indeed, let $H \in \hcurlt{\Oext}$ be a solution of (\ref{weakform}). Since, $\hcurltz{\Oext} \subset
\hcurlt{\Oext}$, to prove that $H$ is solution of (\ref{weakformbis}), we only have to
prove that $H \in \hcurltz{\Oext}$. According to (\ref{decompX}), we have:
\begin{equation} \label{decompH}
H = H_0 + \nabla u, \quad H_0 \in \hcurltz{\Oext}, \quad u \in S,
\end{equation}
and we simply have to prove that $\nabla u$ = 0. If we choose $\wH = \nabla u$
with $v \in S$ in (\ref{weakform}), remarking that $\rt \,\wH = 0$ and $\wH_T
= 0$ on $\Gamma$ and $\partial \Omega$, we get, using (\ref{decompH}):
$$
\omega^2 \int_{\Oext} |\nabla u|^2\; dx = 0 \quad \Longrightarrow \quad \nabla u = 0.
$$
Reciprocally, one proves that any solution of (\ref{weakformbis}) is solution of
(\ref{weakform}) using the same type of argument (decompose the test function
instead of the solution).\\[12pt]
{\bf Step 2.} We now rewrite (\ref{weakformbis}) in the form
\begin{equation} \label{weakform2}
\left\{ \begin{array}{l}
\mbox{ Find } H \in \hcurltz{\Oext} \mbox{ such that}, \\[6pt]
a(H, \wH) = b(H, \wH) + L(\wH),  \quad \forall \; \wH \in \hcurltz{\Oext},
\end{array} \right.
\end{equation}
where we have defined:\\[4pt]
$$
\left. \begin{array}{lll}
\dsp a(H, \wH) & = & \dsp
\int_{\Oext} \big( \, \rt \, H \cdot \overline{\rt \, \wH} +H \cdot \overline{\wH }\,
\big) \; dx - \; i \; \omega (H_T, \wH_T)_{\partial \Omega}
 - \; i \; \omega^2 \; (ZH_T, \wH_T), \\[10pt]
\dsp b(H, \wH) & = & \dsp (1 + \omega^2) \; 
\int_{\Oext} H \cdot \overline{\wH} \; dx ,\\[10pt]
\dsp L(\wH) & = & - \; i \; \omega (g, \wH_T)_{\partial \Oext}+\dsp
\int_{\Oext} f \cdot \overline{\rt \, \wH }\; dx.  
\end{array} \right.
$$
According to (\ref{prop-cd}), one has, with $\alpha < 1$ denoting an arbitrary
positive constant 
$$
\begin{array}{lcl}
|a(H, H)| &\ge& \frac{\alpha}{2} |\Ree a(H, H)| +\frac{1}{2} |\Ime a(H, H)|
\\[6pt]
& \ge& \frac{\alpha}{2}(\|H\|_{H(\rt, \Oext)}^2 - z^* \|H_T\|_\Gamma^2) +
\frac{\omega^2 z_*}{2}\|H_T\|_\Gamma^2 + \frac{\omega}{2} \|H_T\|_{\partial \Omega}^2.
\end{array}
$$
Choosing $\alpha$ such that $\alpha z^*  < \omega^2 z_*$ proves that  the sesquilinear form
$a$ is coercive on $\hcurltz{\Oext}$. On the other hand the sesquilinear form
$b$ is continuous with respect to the $L^2(\Oext)$ norm, which compactly
contains $\hcurltz{\Oext}$. Therefore the Fredholm alternative can be applied to
(\ref{weakform2}): existence uniqueness and stability is equivalent  to the
uniqueness of solutions.
\\[12pt]
{\bf Step 3}. We shall now prove the uniqueness of the solution of
(\ref{weakform}). 
Let $H$ be a solution of (\ref{weakform}) with $f=0 $ and $g=0$. Let us take $\wH =
\overline{H}$ and condider the imaginary part of the resulting equality. We get
\begin{equation*}% \label{identite_im}
\omega \; |H|_{\partial \Omega}^2 + \omega^2 \; \Ree(ZH_T, H_T) = 0.
\end{equation*}
From (\ref{prop-cd}) one deduces that $H_T = 0$ on
$\partial \Omega$ and, using the boundary condition on $\partial \Omega$ that 
$\rt \, H \times n = 0$ on $\partial \Omega$ (this makes sense since, $\rt \;\rt 
\, H\; = \omega^2 \; H$,  and therefore $\rt \,H$ belongs to $H(\rt,
\Oext)$). One concludes using standard unique continuation theorems for
Maxwell's equations (see \cite{Monk}). \proofend 

\subsection{Error estimates for the GIBCs}\label{errorGIBC}
The error estimates rely of some key properties of the boundary operator $
{\cal D}^{\ptre, k}$ that we shall summarize in the following lemma. 
We recall that ${\cal D}^{\ptre, k} = 0$, $\ptre \sqrt{i}$ and $\ptre \sqrt{i}
+ \ptre^2 ({\cal H} -{\cal C})$ for $k= 0, 1 $ and $2$, respectively.
%\subsection{Study of the generalized impedance boundary value
%problems}
For $k=3$, we denote by $A^\ptre$ and $B^\ptre$ the two operators
\begin{equation*}% \label{defAB}
A^{\ptre} := \big(1   - {\ptre^2} \; \grag \divg\big)^{-1}, 
\quad B^{\ptre} = \big( 1 + {\ptre^2} \;  \Rotg \rotg \big)^{-1} \; .
\end{equation*}
By Lax-Milgram's Lemma these operators are well defined as continuous
operators from $L^2_t(\Gamma)$ to respectively $H(\divg, \Gamma)$ and
$H(\rotg, \Gamma)$.
Setting
$$
\begin{array}{l}
\alpha_\delta := \dsp\frac{ 1}{2\sqrt{2}}   +  \ptre ({\cal 
  H}-{\cal C})  + \frac{\ptre^2}{2 \sqrt{2}} \left({\cal C}^2 -  {\cal H}^2
  +\eps_r \omega^2\right),
\\[12pt]
\beta_\delta := \dsp\frac{ 1}{2\sqrt{2}}   - \frac{\ptre^2}{2 \sqrt{2}} \left({\cal C}^2 -  {\cal H}^2
  +\eps_r \omega^2\right),
\end{array}
$$
the expression (\ref{D3modified}) of ${\cal D}^{\ptre, 3}$ can be written in
the form \begin{equation*}% \label{defd3}
\left\{ \begin{array}{lcl}
{\cal D}^{\ptre, 3}\varphi &= & \dsp \ptre \, \alpha_\ptre \varphi +
\frac{\sqrt{2}}{4} \;\ptre\left(  A^\ptre \varphi + \ptre^2\Rotg\rotg B^\ptre\varphi  \right)
\\[12pt]
& + & \dsp i \ptre \,\beta_\ptre \varphi + i
\frac{\sqrt{2}}{4} \; \ptre\left( B^\ptre \varphi - \ptre^2\grag \divg A^\ptre\varphi\right).
\end{array} \right.
\end{equation*}
The fundamental properties of the operators ${\cal D}^{\ptre, k}$ are
summarized in the following lemma.
\begin{lemma} \label{lembuy} Let $k = 1,2 $ or $3$.
There exist a constant $\ptre_k >0$ and two constants $C_1 > 0$ and $C_2 >0$,  independent of $\ptre$, such that
\begin{equation}
\label{errorbuy}
\begin{array}{ll}
(i) & \|{\cal D}^{\ptre, k} \varphi\|_\Gamma \le C_1 \; \ptre \;
\|\varphi\|_\Gamma, 
\\ [6pt]
(ii) & \Ree ({\cal D}^{\ptre, k} \varphi, \varphi)_\Gamma \ge C_2 \; \ptre
\; \|\varphi\|_\Gamma^2,
%\\[6pt]
%(iii) & \Im ({\cal D}^{\ptre, k} \varphi, \varphi) \ge 0
\end{array}
\end{equation}  
for all $\varphi \in L^2_t(\Gamma)$ and $\ptre \le \ptre_k$. 
\end{lemma}

\proof
These properties are straightforward for $k = 1 $ and $2$. We shall
concentrate on the case $k=3$. We first observe that 
$\alpha_\ptre$ and $\beta_\ptre$ are bounded functions on $\Gamma$,  and
if $\varphi \in L^2_t(\Gamma)$ then
$\ptre^2 \grag \divg A^\ptre\varphi = (A^\ptre\varphi - \varphi) \in
L^2_t(\Gamma)$  and $ \ptre^2 \Rotg\rotg B^\ptre\varphi = (-B^\ptre\varphi + \varphi)\in
L^2_t(\Gamma)$. Therefore ${\cal D}^{\ptre, 3}\varphi \in L^2_t(\Gamma)$ and
one has
\begin{equation} \label{defd33}
\left\{ \begin{array}{lcl}
({\cal D}^{\ptre, 3}\varphi, \psi)_\Gamma &= & \dsp \ptre \; (\alpha_\ptre \varphi, \psi)_\Gamma +
\frac{\sqrt{2}}{4} \; \ptre\left( (A^\ptre \varphi, \psi)_\Gamma + \ptre^2 \; (\Rotg\rotg B^\ptre\varphi ,\psi)_\Gamma \right)
\\[12pt]
& + & \dsp i \ptre \; (\beta_\ptre \varphi, \psi)_\Gamma + i
\frac{\sqrt{2}}{4} \; \ptre\left( (B^\ptre \varphi, \psi)_\Gamma - \ptre^2\; (\grag \divg A^\ptre\varphi,\psi)_\Gamma \right)
\end{array} \right.
\end{equation}
 for all $\varphi, \;\psi \in L^2_t(\Gamma)$. For $\ptre$ sufficiently small, the functions $\alpha_\ptre$ and $\beta_\ptre$
satisfy
\begin{equation}\label{cort}
0< \alpha_*<|\alpha_\ptre|<\alpha^* \quad \mbox{and} \quad 0<
\beta_*<|\beta_\ptre|<\beta^*
\end{equation}
for some positive constants  $\alpha_*$, $\alpha^*$, $\beta_*$ and $\beta^*$
independent of $\ptre$. 
On the other hand, from the identities
\begin{equation*}%\label{adbd}
(1 - \ptre^2 \grag \divg ) A^\ptre \varphi = \varphi  \quad \mbox{and} \quad(1
+ \ptre^2 \Rotg \rotg ) B^\ptre \varphi = \varphi
\end{equation*}
one respectively deduces 
\begin{equation}\label{adel}
\left\{\begin{array}{ccl}
(A^\ptre \varphi, \varphi)_\Gamma &= & \|A^\ptre\varphi\|_\Gamma^2 + \ptre^2 \|\divg A^\ptre\varphi\|_\Gamma^2
\\[12pt]
-(\grag \divg A^\ptre \varphi, \varphi)_\Gamma &= & \|\divg A^\ptre\varphi \|_\Gamma^2+
\ptre^2 \|\grag\divg A^\ptre\varphi\|_\Gamma^2
\end{array}\right.
\end{equation}
and
\begin{equation}\label{bdel}
\left\{\begin{array}{ccl}
(B^\ptre \varphi, \varphi)_\Gamma &= & \|B^\ptre\varphi\|_\Gamma^2 + \ptre^2 \|\rotg B^\ptre\varphi\|_\Gamma^2
\\[12pt]
(\Rotg\rotg B^\ptre \varphi, \varphi)_\Gamma &= & \|\rotg B^\ptre\varphi \|_\Gamma^2+
\ptre^2 \|\Rotg\rotg B^\ptre\varphi\|_\Gamma^2
\end{array}\right.
\end{equation}
Property $(ii)$ is obtained as an immediate consequence of (\ref{adel}),
(\ref{bdel}) and (\ref{cort}) when applied to (\ref{defd33}) with $\psi =
\varphi$. %Property $(iv)$ is obtained   
 Identities (\ref{adel}) and
(\ref{bdel}) also respectively imply, 
\\[6pt]
\centerline{$\;$\quad $
\|A^\ptre\varphi\|_\Gamma \le \|\varphi\|_\Gamma$, \quad$ \ptre^2 \|\grag\divg
A^\ptre\varphi\|_\Gamma \le \|\varphi\|_\Gamma$,}
\\[6pt]
\centerline{$\;$\quad $ \|B^\ptre\varphi\|_\Gamma \le
\|\varphi\|_\Gamma $  and $ \ptre^2 \|\Rotg\rotg B^\ptre\varphi\|_\Gamma \le \|\varphi\|_\Gamma.$}
\\[6pt]
Property $(i)$ is then easily obtained from (\ref{defd33}) with $\psi =
{\cal D}^{\ptre, 3}\varphi$ and using these estimates, as well as
(\ref{cort}).
\proofend
 
\subsection{Error estimates for the GIBCs}\label{errorGIBC-section}
We shall set for $k = 0,1,2,3$,
\begin{equation} \label{error2}
\left\{ \begin{array}{l}
\tdEke = E^{\ptre,k}_e-\sum_{\ell=0}^k E_e^\ell, \\[12pt]
\tdHke = H^{\ptre,k}_e-\sum_{\ell=0}^k H_e^\ell.
\end{array} \right.
\end{equation}
Using (\ref{bdry-gibc-1}), together with (\ref{errorD}) and (\ref{boundR})
when $k=3$,  we see that $(\tdEke, \tdHke)\in
{\cal V}_E^k \times {\cal V}_H^k$ is solution  of the boundary value problem:
\begin{equation} \label{bof1}
\left\{
\begin{array}{ll}
\rt \; \tdHke  +i \omega \; \tdEke = 0,  &\;\;\mbox{ in } \; \Oext,\\[12pt]
\rt \; \tdEke -i \omega \; \tdHke  = 0,  &\;\; \mbox{ in } \; \Oext,\\[12pt]
\big(\tdEke\big)_T - \tdHke  \times n  = 0, &\;\;\mbox{ on } \;
\partial \Omega,\\[12pt]
\tdEke \times n + \omega \, {\cal D}^{\ptre, k}(
\tdHke)_T= \ptre^{k+1} \, \varphi^\ptre_k &\;\;
\mbox{ on } \Gamma, 
\end{array}
\right.
\end{equation}
where the tangential vector fields $\varphi^\ptre_k$ remain bounded with
respect to $\ptre$ in all 
spaces $H^s_t(\Gamma)^3$.\\[12pt]
Eliminating $\tdEke$, we see that $\tdHke\in {\cal V}_H^k$ satisfies
$$%\begin{equation} \label{pbH}
\left\{ \begin{array}{ll}
\rt \big( \, \rt \, \tdHke) - \omega^2 \; \tdHke = 0 & \mbox{ in } \; \Oext,
\\[12pt]
\rt \, \tdHke \times n - i \; 
\omega^2 \;{\cal D}^{\ptre, k}(\tdHke)_T= \ptre^{k+1} \, \varphi^\ptre_k, 
& \mbox{ on } \; \Gamma, \\[12pt]
\rt \, \tdHke \times n - i \; \omega \; (\tdHke)_T = 0
, & \mbox{ on } \; \partial \Omega.
\end{array} \right.
$$%\end{equation}
The proof of error estimates is based on some key a priori estimates that we
shall give hereafter.   We multiply by $\overline{\tdHke}$ the equation satisfied by ${\tdHke}$ in 
$\Oext$, integrate over $\Oext$ and use Green's formula to obtain, after 
having used the boundary conditions on $\partial \Omega$ and $\Gamma$: 
\begin{equation} \label{esti-enrerg}
\hspace*{-0.3cm} \left| \begin{array}{rl}
\dsp \int_{\Oext} \big( |\rt \, \tdHke|^2 - \omega^2 |\tdHke|^2 \big) \; dx -
i \; \int_{\partial \Omega} |(\tdHke)_T|^2 \; d \sigma  & \\[12pt]
\dsp - \; i \; \omega^2 \; \big({\cal D}^{\ptre, k}(\tdHke)_T, (\tdHke)_T\big)_\Gamma & =
\dsp \delta^{k+1} \left< \varphi_{k}^\ptre,  (\tdHke)_T\right>_\Gamma
, 
\end{array} \right.
\end{equation}
where $\left<\; , \;\right>_\Gamma$ here denotes a duality pairing between
  $H^{-1/2}(\mathrm{div}, \Gamma)$ and $H^{-1/2}(\rt, \Gamma)$. Considering the imaginary part of (\ref{esti-enrerg}) and using
(\ref{errorbuy})-(ii) together with trace theorems in $H(\rt, \Oext)$,  one obtains the existence of two non negative constants
$C_1$ and $C_2$ independent of $\ptre$ such that
\begin{equation}\label{bri2}
 C_1 \; \ptre \; \| (\tdHke)_T \|_{\Gamma}^2 + \| (\tdHke)_T \|_{\partial \Omega}^2 \le C_2 \; \ptre^{k+1} \; \| \tdHke \|_{H(\rt, \Oext)}.
\end{equation}
More precisely we have $C_1 =0$ for $k=0$ and $C_1 >0$ for $k \neq 0$. Using (\ref{errorbuy})-(i) and (\ref{bri2}) one also deduces that 
\begin{equation} \label{estimI}
|\Ime \big({\cal D}^{\ptre, k}(\tdHke)_T, (\tdHke)_T\big)_\Gamma| \le  C_3 \; 
\ptre^{k+1} \; \| \tdHke \|_{H(\rt, \Oext)},
\end{equation}
for some constant $C_3$ independent of $\ptre$. 
Now considering  the real part of (\ref{esti-enrerg})
and using (\ref{estimI}) as well as the trace theorem in $H(\rt, \Oext)$,
one gets the existence of two positive constants
$C_4$ and $C_5$ independent of $\ptre$ such that
\begin{equation} \label{bri1}
 \| \tdHke \|_{H(\rt, \Oext)}^2 \le C_4 \; \ptre^{k+1} \; \| \tdHke \|_{H(\rt, \Oext)} + C_5 \; \| \tdHke \|_{L^2(\Oext)}^2.
\end{equation}
Based on these a priori estimates we are in position to prove the following result.
\begin{lemma} \label{lemerreor}
For $k = 0, 1, 2 $ or 3, there exist a constant $C$ independent of $\ptre$ and $\delta_0 >0$ such
that
$$
\| \tdEke \|_{H(\rt, \Oext)} + \| \tdHke \|_{H(\rt, \Oext)} \le C_k \;
\ptre^{k+1} 
$$
for all $\ptre \le \delta_0$.
\end{lemma}
\proof
According to estimate (\ref{bri1}) and first two equations of (\ref{bof1}) it
is sufficient to prove the existence of a constant $C$ independent of $\ptre$
such that 
\begin{equation}\label{suff2}
\| \tdHke \|_{L^2(\Oext)} \le C \; \ptre ^{k+1}.
\end{equation}
Let us assume that (\ref{suff2}) does not hold, i.e. $\lambda_\ptre := \ptre
^{k+1}/ \| \tdHke \|_{L^2(\Oext)}$ goes to $0$ as $\ptre \to 0$,   and consider the scaled fields
$$h^\ptre =  \tdHke / \| \tdHke \|_{L^2(\Oext)} \quad \mbox{ and } \quad
e^\ptre = \tdEke / \| \tdHke \|_{L^2(\Oext)}.$$
Dividing (\ref{bri1}) by $\| \tdHke \|_{L^2(\Oext)}^2$ implies in particular
that $(h^\ptre)$  is a bounded sequence in ${H(\rt, \Oext)}$. Dividing (\ref{bri2}) by the
same quantity and using the latter result shows that
\begin{equation}  \label{toiuy}
C_1 \; \ptre \; \| h^\ptre_T \|_{\Gamma}^2 + \| h^\ptre_T \|_{\partial
  \Omega}^2 \rightarrow 0 \mbox{ as } \ptre \rightarrow 0.
\end{equation}
The last boundary condition in (\ref{bof1}) combined
with property (\ref{errorbuy})-(i) shows in particular that
$$
\| e^\ptre_T \|_{\Gamma} \le C_6 \; \delta \; \| h^\ptre_T \|_{\Gamma} + \lambda_\ptre \, \|\varphi^\ptre_k\|_{\Gamma} %\ptre ^{k+1}/\| \tdHke \|_{L^2(\Oext)}
$$
(with the alternative, $C_6$ and $C_1 >0$ or $C_6=C_1 =0$). We therefore
conclude that  $\| e^\ptre_T \|_{\Gamma}$ goes to $0$ as $\ptre \rightarrow
0$. The first three equations of (\ref{bof1}) implies that $e^\ptre$ is a
bounded sequence  in $\hcurltz{\Oext}$ (see definition in the proof of Theorem~\ref{exist-unicite2*}). Therefore, up to extracted subsequence, one can
assume that   $e^\ptre$ converges strongly in $L^2(\Oext)$ and weakly in $\hcurltz{\Oext}$
to some $e\in \hcurltz{\Oext}$. Passing to the limit as $\ptre \to 0$ in (\ref{bof1}), we
observe that   $e\in \hcurltz{\Oext}$ is solution of  
$$%\begin{equation} \label{bof1}
\left\{
\begin{array}{ll}
\rt \; \rt \;e  -\omega^2 \; e = 0,  &\;\;\mbox{ in } \; \Oext,\\[6pt]
\rt \; e\times n  - i\omega e_T  \times n  = 0, &\;\;\mbox{ on } \;
\partial \Omega,\\[6pt]
e\times n =0 &\;\;
\mbox{ on } \Gamma, 
\end{array}
\right.
$$%\end{equation}
and therefore $e=0$. We then deduce that $\rt \, h^\ptre$ strongly converges to $0$ in
$L^2(\Oext)$. Coming back to identity (\ref{esti-enrerg}) and considering the
real part, one deduces after division by  $\| \tdHke \|_{L^2(\Oext)}^2$ that
\begin{equation} \label{bri3}
 \| h^\ptre \|_{H(\rt, \Oext)}^2 \le \lambda_\ptre \, \tilde C_0  \; \|
 h^\ptre \|_{H(\rt, \Oext)} + \tilde C_1 \; \| \rt h^\ptre
 \|_{L^2(\Oext)}^2 + \tilde C_2; |\Ime ({\cal D}^{\ptre, k}h^\ptre_T,
 h^\ptre_T)|. 
\end{equation}
Property (\ref{errorbuy})-(i) shows that
$$
|\Ime ({\cal D}^{\ptre, k}h^\ptre_T,  h^\ptre_T) |\le \tilde C_3 \ptre
|h^\ptre_T|_\Gamma^2 \quad \to \quad 0
$$
according to (\ref{toiuy}). Therefore, considering the limit as $\ptre \to 0 $ in (\ref{bri3}) implies that  $h^\ptre$
strongly converges to $0$ in $H(\rt, \Oext)$, which contradicts
$\|h^\ptre\|_{L^2(\Oext)} = 1$.
\proofend

%%%%%%%%%%%%%%%%%%%%%%%%%%%%%%%%%%%%%%%%%%%%%%%%%%%%%%%%%%%%%%%%%%%%%%%%

\begin{appendix}
\section{Technical Lemmas} \label{AppA}

The first Lemma is a slight variation of the classical trace lemma in $\hcurl{O}$.

\begin{lemma}\label{trace-lemma}
Let $O$ be a bounded open subset of $\R^3$ of  class
$C^2$ (and which is locally from one side of its normal). Then there exists a constant $C$ depending only on $O$ such
that
\begin{equation} \label{trace-ineq}
\dsp \|u\times n\|^2_{H^{-\frac{1}{2}}(\partial O) } \leq C \; 
\|u\|_{L^2(O)} \; \left( \|\rt\, u\|_{L^2(O)} + \|u\|_{L^2(O)} \right)
\quad \forall u \in \hcurl{O}.
\end{equation}
\end{lemma}
\begin{remark} \label{HS}
When $O = \R^3_+ := \big \{ y \in \R^3 \; /\;  y_3 >0 \big\}$, one can easily check,
using partial Fourier transform in the variables $(y_1, y_2)$, that
$$%\begin{equation} \label{trace-ineq}
\dsp \forall \; u \in \hcurl{\R^3_+}, 
\quad \|u\times n\|^2_{H^{-\frac{1}{2}}(\partial \R^3_+) } 
\leq \; 2 \; \|u\|_{L^2(\R^3_+)} \; ( \|\rt\, u\|_{L^2(\R^3_+)} +
\|u\|_{L^2(\R^3_+)} ). 
$$%\end{equation}
\end{remark}
\proof The idea is to see how the proof for the half-space (cf. previous remark)
  is modified when the boundary of $O$ is not flat.\\[12pt]
Let us consider first the cases where there exists $a>0,~ b>0,
~\delta>0$ and $h \in C^2(\R^2) \cap W^{2, \infty}(\R^2)$ such that
\begin{eqnarray*}
\left | \begin{array}{llll} \dsp \partial O \cap \mbox{ supp } u
\subset  \Sigma :=\dsp \Big\{ \; \varphi(y_1, y_2, 0) \; ;\; (y_1, y_2) \in
]-a,a\,[\times ]-b,b\,[ \; \Big\},
\\[12pt]
\dsp \varphi\Big(]-a,a[
\times ]-b,b[ \times  ]0,\delta[ \Big)\subset O,
\\[12pt]
\dsp u \mbox{ is compactly supported in } \varphi\Big( \; ]-a,a \, [
\times ]-b,b \, [ \times  ]0,\delta\,[ \; \Big),
\end{array} \right.
\end{eqnarray*}
where
$$
\varphi(y_1,y_2,y_3) = \big(y_1,y_2,y_3 + h(y_1,y_2)\big).
$$
Setting $ \tilde{u} := u \circ \varphi$ and $ \tilde{n} := u \circ \varphi$ ,
one has for $y_3=0$,
\begin{equation}\label{uTimesn}
\tilde{u} \times \tilde{n} = f \; 
\Big( \tu_2 + \tu_3
\frac{\partial h}{\partial y_2} \, , \; \tu_3 \frac{\partial h}{\partial
y_1} - \tu_1, \; - \; \tu_1 \frac{\partial h}{\partial y_2} + \tu_2
\frac{\partial}{\partial y_1} h \Big), \quad f := \big( 1+|Dh|^2 \big)^{-\frac{1}{2}}
\end{equation}
% on $\Sigma$, where $\dsp n:= \frac{1}{\sqrt{1+|Dh|^2}} \; \Big (-
% \frac{\partial h}{\partial y_1}, -\frac{\partial h}{\partial y_2},
% 1 \Big)$. 
On the other hand, $(e_1, e_2, e_3)$ denoting the canonical basis in the
$y$-space, 
\begin{equation}\label{RefLem8}
\left|\begin{array}{llll}  (i) & (\rt \, u \circ \varphi) \cdot e_1  &=& \dsp
\frac{\partial\tu_3}{\partial y_2}  - \frac{\partial\tu_3}{\partial
y_3}
 \frac{\partial h}{\partial y_2} - \frac{\partial \tu_2}{\partial y_3},
\\[12pt]
 (ii) & (\rt \, u \circ \varphi)\cdot e_2  &=&  \dsp \frac{\partial\tu_1}{\partial y_3}
 - \frac{\partial \tu_3}{\partial y_1}
 + \frac{\partial\tu_3}{\partial y_3}  \frac{\partial h}{\partial y_1}, 
\\[12pt]
 (iii) & (\rt \, u \circ \varphi)\cdot e_3  &=&  \dsp \frac{\partial\tu_2}{\partial
y_1}  - \frac{\partial \tu_1}{\partial y_2} +
\frac{\partial \tu_3}{\partial y_3} (\frac{\partial h}{\partial y_2}
 - \frac{\partial h}{\partial y_1}).
\end{array} \right.
\end{equation}
Let us set $\dsp u_2^*= \tu_2 + \tu_3 \frac{\partial h}{\partial y_2} \;$, 
since, 
$$
\frac{\partial u_2^*}{\partial y_3}  = \frac{\partial\tu_2}{\partial y_3}
 + \frac{\partial\tu_3}{\partial y_3}
\frac{\partial h}{\partial y_2}
+ \tu_3 \; \frac{\partial^2 h }{\partial y_3 \partial y_2}  \\
$$
using the formula (\ref{RefLem8})-$(i)$ one gets, setting $r:= (\rt \, u \circ
\varphi) = (r_1, r_2, r_3)$
\begin{equation} \label{f1}
\frac{\partial u_2^*}{\partial y_3} = - r_1 +
\frac{\partial\tu_3}{\partial y_2} 
+ \tu_3 \frac{\partial^2 h}{\partial y_3 \partial y_2}  \\
\end{equation}
In what follows, we use the Fourier transform ${\cal F}$ in the variables
$(y_1, y_2)$ and denote by $(\xi_1, \xi_2)$ the dual variable, by $\hu_i$ the 
Fourier transform of $\wu_i$ and $\hu_2^*$ the Fourier transform of $u_2^*$. 
By definition of the norm in $H^{-\frac{1}{2}}(\R^2)$,
$$
\|u_2^*(\cdot, \cdot,0)\|^2_{H^{-\frac{1}{2}}} = 
\int_{\R^2} (1 + |\xi|^2)^{-\frac{1}{2}} \; |\hu_2^*(\xi_1,\xi_2,0)|^2 \; d\xi_1 d\xi_2 
$$
Since
$$
|\hu_2^*(\xi_1,\xi_2,0)|^2 
=\dsp - 2 \Ree {\dsp \int_0^\delta} \Big[ \; \frac{\partial
  \hu_2^*}{\partial y_3} \; \overline {\hu_2^*} \; \Big] (\xi_1,\xi_2,y_3)\;
dy_3
$$
using (\ref{f1}), we have, $\widehat{r}_1$ being the Fourier transform of $r_1$,
\begin{equation*}
\left| \begin{array}{lcl}
|\hu_2^*(\xi_1,\xi_2,0)|^2 
&=& \dsp 2 \; \Ree {\dsp\int_0^\delta } \Big[ \; \widehat{r}_1 \; \overline
{\hu_2^*} \; \Big] (\xi_1,\xi_2,y_3) \; dy_3 
\\[12pt] 
&-& 2 \; \Ree {\dsp\int_0^\delta } i \xi_2
\big[ \wu_3  \; \overline {\hu_2^*} \, \big] (\xi_1,\xi_2,y_3) \; dy_3
\\[12pt]  
&-& \dsp 2 \; \Ree {\dsp\int_0^\delta } \Big[ \; {\cal F}\big(\tu_3
\frac{\partial^2 h}{\partial y_3 \partial y_1} \big) \; \overline {\wu_2^*}\;
\Big] (\xi_1,\xi_2,y_3) \; dy_3.
\end{array}\right.
\end{equation*}
We divide the above equality by $(1 + |\xi|^2)^{-\frac{1}{2}}$ and integrate
over $\xi$. Next we use $(1 + |\xi|^2)^{-\frac{1}{2}} \leq 1$, $ |\xi_2| \; (1
+ |\xi|^2)^{-\frac{1}{2}} \leq 1$ and Plancherel's theorem to obtain, since $h
\in W^{2,\infty}(\R^2)$, 
\begin{equation*}
\|u_2^*(\cdot, \cdot,0)\|^2_{H^{-\frac{1}{2}}}
\leq \dsp 2 {\dsp\int_0^\delta }\int_{\R^2}  |{r}_1| \; |u_2^*| \; dy
+ C{\dsp\int_0^\delta }\int_{\R^2} |u_3|  \; |u_2^*| \; dy.
\end{equation*}
Coming back to the variable $x$ through the change of variable $x =
\varphi(y)$, we easily get, since $|u_2^*| \leq |u_2| + C |u_3|$
\begin{equation*}
\|u_2^*(\cdot, \cdot,0)\|^2_{H^{-\frac{1}{2}}} \leq C \left( \|u\|^2_{L^2(O)} + \|u\|_{L^2(O)} \|\rt
u\|_{L^2(O)} \right)
\end{equation*}
where the constant $C$ only depends on $h$.\\[12pt]
Finally, using the lemma  \ref{lemPro} (notice that
$f := \big( 1+|Dh|^2 \big)^{-\frac{1}{2}}$ belongs to $W^{1, \infty}$), we
get 
\begin{equation} \label{estrace1}
\|f \big( \tu_2 + \tu_3 \frac{\partial h}{\partial y_2} \big) \|^2_{H^{-\frac{1}{2}}} \leq C \left( \|u\|^2_{L^2(O)} + \|u\|_{L^2(O)} \|\rt 
u\|_{L^2(O)} \right) 
\end{equation}
In the same way, one obtains 
\begin{equation}\label{estrace2}
\left\{\begin{array}{lllll}
&\dsp \|f \big(\tu_3 \frac{\partial h}{\partial y_1} -
\tu_1\big) \|_{H^{-\frac{1}{2}}(\Gamma)} &\le & C \dsp \left(
\|u\|^2_{L^2(O)} + \|u\|_{L^2(O)} \|\rt u\|_{L^2(O)} \right) \\[12pt]
& \dsp \| f \big(- \tu_1 \frac{\partial h}{\partial y_2} + \tu_2
\frac{\partial h}{\partial y_1} \big)\|_{H^{-\frac{1}{2}}(\Gamma)} &\le
& \dsp C \left( \|u\|^2_{L^2(O)} + \|u\|_{L^2(O)} \|\rt
u\|_{L^2(O)} \right)
\end{array}\right.
\end{equation}
Thanks to (\ref{uTimesn}) and by definition of the norm in
$H^{-\frac{1}{2}}(\Gamma)$, estimates (\ref{estrace1}) and (\ref{estrace2}) lead to the
desired inequality.\\[12pt] 
Obtaining the same  inequality in the general case can be deduced by the using
a partition of unity $(\varphi_i)_{i=1, ..., N}$ of $O$ and
noticing that 
$$
\|\rt \varphi_i u\|_{L^2(O)} = \|\varphi_i ~ \rt u + \nabla \varphi_i \times
u\|_{L^2(O)} \le \|\varphi_i\|_\infty \; \| \rt u \|_{L^2(O)}+ \|\nabla
\varphi_i\|_\infty \;  
\|u\|_{L^2(O)}.
$$
\proofend
\begin{lemma}\label{lemPro} Let $ f \in W^{1,\infty} (R^n) $ and
$ g \in H^{-\frac{1}{2}}(R^n) $  then $f\,g \in H^{-\frac{1}{2}}
(R^n)$ and one has, 
\begin{eqnarray*}
\|f \, g\|_{H^{-\frac{1}{2}}(R^n)} \leq 3^{\frac{1}{4}} \; 
\|f\|_{W^{1,\infty}(R^n)} \; \|g\|_{H^{-\frac{1}{2}}(R^n)}.
\end{eqnarray*}
\end{lemma}
\proof We first remark that if $\psi \in H^{\frac{1}{2}}(R^n)$, then
$$
\|f \, \psi\|_{H^{\frac{1}{2}}(R^n)}\leq 3^{\frac{1}{4}} \; 
\|f\|_{W^{1,\infty}(R^n)} \; \|\psi\|_{H^{\frac{1}{2}}(R^n)}
$$
which is deduced by interpolation from the (obvious) inequalities:
$$
\left\{ \begin{array}{ll}
\forall \; \psi \in L^2(\R^n), \quad & \|f \,
\psi\|_{H^{\frac{1}{2}}(R^n)}\leq   
\|f\|_{L^{\infty}(R^n)} \; \|g\|_{L^{{2}}(R^n)},
\\[12pt]
\forall \; \psi \in H^1(\R^n), \quad & \|f \, \psi\|_{H^{1}(R^n)}\leq
\Big( 2 \; \|f\|_{L^{\infty}(R^n)}^2 + \|f\|_{W^{1,\infty}(R^n)}^2 \Big)^
{\frac{1}{2}} \; \|\psi\|_{H^{1}(R^n)}.
\end{array} \right.
$$
Next, if $g \in H^{-\frac{1}{2}}(R^n)$, we have, for any $\psi \in
H^{\frac{1}{2}}(R^n)$, 
$$
|\left< fg, \psi \right>| = |\left< g, f \psi \right>| \leq
\|g\|_{H^{-\frac{1}{2}}(R^n)} \; \|f \, \psi\|_{H^{\frac{1}{2}}(R^n)} \leq 
3^{\frac{1}{4}} \; \|g\|_{H^{-\frac{1}{2}}(R^n)} \;
\|f\|_{W^{1,\infty}(R^n)} \; \|\psi\|_{H^{-\frac{1}{2}}(R^n)},
$$
from which one easily concludes.
\proofend

\begin{lemma}\label{lemMaxpre2}
Let  $O \subset R^3$ be a bounded open set with a $C^2$ boundary $\Gamma$. There exists a constant $C$ that depends only on
$\Gamma$ such that
\begin{eqnarray*}
\nodt{u}{H^{\frac{1}{2}}(\Gamma)} \leq C (\nodt{\nabla_{\Gamma} u}{H^{-\frac{1}{2}}(\Gamma)} + \nod{u}{L^2(\Gamma)}) \qquad \forall\, u \in H^{\frac{1}{2}}(\Gamma).
\end{eqnarray*}
\end{lemma}
\proof In the case $O=\{(x_1,x_2,x_3)\in \R^3 \; /\;  x_3 \ge 0\}$ one can check by using Fourier transform in the plane $(x_1, x_2)$ that
$$
\nodt{u}{H^{\frac{1}{2}}(\Gamma)}^2 = \nodt{\nabla_{\Gamma} u}{H^{-\frac{1}{2}}(\Gamma)}^2 + \nod{u}{L^2(\Gamma)}^2.
$$
The inequality is therefore trivially verified in this case. The general case
can be easily deduced by using local parameterizations of the boundary
$\Gamma$. This is where the $C^2$-regularity of $\Gamma$ is taken into
account. 
\proofend

\noindent Our next lemma is a sharper version of classical compact-embedding
theorem for spaces of  $L^2$ functions with bounded divergence and curl into
$L^2$. We set 
$$H(\rt, \dv, O) := \{ u \in L^2(O)^3 \; / \; \rt\,u \in L^2(O)^3 \mbox{ and } \dv\, u \in L^2(O)\}$$
equipped with the norm 
$$
\nodt{u}{H(\rt, \dv, O)}^2 = \nodt{u}{L^2(O)}^2 + \nodt{\rt\,u}{L^2(O)}^2 + \nodt{\dv\,u}{L^2(O)}^2.
$$

\begin{lemma}\label{compact-lemma}
Let $O \subset R^3$ be a bounded simply connected open set with 
$C^2$ boundary $\Gamma$. Then  every bounded  sequence $(u_k)_{k \in \N }$ 
of  $H(\rt, \dv, O)$ such that
\begin{equation}\label{Hypcompact-lemma}
 ({u_k}_{|\Gamma}  \times n)_{k \in \N} \textrm{ is  convergent in }
H^{-\frac{1}{2}}_t(\Gamma)
\end{equation}
 has a convergent  subsequence $ (u_{k'}) $
 in $ L^2(O)^3$.
\end{lemma}
\proof Our proof is an adaptation of the proof given by Costabel in the case
where, instead of (\ref{Hypcompact-lemma}), one has an $L^2$ control of the
boundary term of the sequence (see Theorem~2 of \cite{Costabel}). 
\\[12pt]
The idea is to make a Helmholtz decomposition of $u_k$ of the form:
\begin{equation} \label{decompuk}
u_k = w_k + \nabla p_k, \quad (w_k, p_k) \in H^1(O) \times L^2(O), 
\quad \dv \, w_k = 0, 
\end{equation}
constructed in such a way that:
\begin{itemize}
\item[(i)] $w_k$ is bounded in $H^1(O)$ (and thus admits a converging 
subsequence in $L^2(O)^3$) : this uses the fact that $\rt \, u_k$ 
is bounded in $L^2(O)^3$,
\item[(i)] $\nabla p_k$  admits a converging 
subsequence in $L^2(O)^3$ : this uses the fact that $\dv \, u_k$ 
is bounded in $L^2(O)^3$ and that $({u_k}  \times n)_{|\Gamma}$ converges according to (\ref{Hypcompact-lemma}).
\end{itemize}
In order to construct $w_k$ from $\rt \; u_k$ we first construct an extension 
of $\rt \; u_k$ in $\R^3$ which has a compact support (independent of $k$) 
and is divergence free. For this, we choose a ball $B$ containing 
$\overline{O}$ in its interior and will constitute the support of 
the extension of $\rt \; u_k$.\\[12pt]
First notice that since $\rt \, u_k \in H(\dv;O)$, the trace 
$\rt \, u_k \cdot n_{|\Gamma}$ is well defined in $H^{-\frac{1}{2}}(\Gamma)$ 
and satisfies, since $\dv \, \rt \, u_k = 0$:
\begin{equation} \label{esti-rotuk}
\|\rt \, u_k \cdot n\|_{H^{-\frac{1}{2}}(\Gamma)} \leq C \; 
\| \rt \, u_k\|_{L^2(O)}, 
\quad \left< \rt \, u_k \cdot n, 1 \right>_{\Gamma} =0.
\end{equation}
Therefore, the following Neumann problem in $B \setminus O$:
\begin{equation}\label{proPre}
\left\{\begin{array}{ll}
\Delta \varphi_k = 0 & \textrm{ in } B \setminus \overline{O}, \\[12pt]
\partial_n \varphi_k = (\rt \, u_k)\cdot  n &\textrm{ on } \partial O,  \\[12pt]
\partial_n \varphi_k = 0 &\textrm{ on } \partial B,
\end{array}\right.
\end{equation}
admits a solution $\varphi_k \in H^1(B \setminus O)$, which 
is unique if we impose in addition
\begin{equation*}% \label{int0}
\dsp \int_{B\setminus O}\varphi_k \; dx=0.
\end{equation*}
Moreover, using the Poincare-Wirtinger inequality and (\ref{esti-rotuk})
$$
\int_{B \setminus \overline{O}} |\nabla \varphi_k|^2 \; dx = 
\left< \,  (\rt \, u_k)\cdot  n , \overline{\varphi_k} \; \right>_{\Gamma} 
\leq C \; \|\rt \, u_k\|_{L^2(O)} \;
\| \nabla \varphi_k\|_{L^2(B \setminus \overline{O})} 
$$
from which we deduce
\begin{equation} \label{estiphik}
\| \nabla \varphi_k\|_{L^2(B \setminus \overline{O})} \leq
C \; \|\rt \, u_k\|_{L^2(O)}.
\end{equation}
~\\[0pt]
We then introduce the extension $\chi_k$ of $\rt \; u_k$ as
\begin{equation*}
\chi_k = \begin{cases}
\rt \; u_k &\qquad \textrm{ in } O,  \\
\nabla \varphi_k& \qquad\textrm{ in }  B \setminus \overline{O}, \\
0 &\qquad\textrm{ in } R^3 \setminus \overline{B}.\\
\end{cases}
\end{equation*}
Of course, $\chi_k \in L^2(R^3)$, is compactly supported in $\overline{B}$ and satisfies, thanks to (\ref{estiphik})
\begin{equation*} %\label{esti-chik}
\|\chi_k\|_{L^2(\R^3)} \leq C \; \|\rt \, u_k\|_{L^2(O)} \; .
\end{equation*}
Moreover, the two boundary conditions in
(\ref{proPre}) have been chosen in order to enforce the continuity of
 $\chi_k \cdot n$ across $\partial O$ and $\partial B$ so that:
$$\dv \, \chi_k=0 \quad \mbox{in } \R^3. $$
Next, we introduce $\Psi_k$ as the unique solution of he following Laplace 
problem in $\R^3$
\begin{equation} \label{defpsik}
- \; \Delta \Psi_k = \chi_k \, \quad \mbox{in } \R^3, \quad 
\Psi_k \in H^2_{loc}(\R^3)^3, \quad \nabla \Psi_k \in L^2(\R^3)^3.
\end{equation}
It is well known that $\Psi_k$ is given by:
\begin{equation*} %\label{exppsik}
\Psi_k = G * \chi_k, \quad G = \frac{1}{4 \pi |x|}, \quad \mbox{(the 
fundamental solution of the Laplace operator)}
\end{equation*}
and satisfies in particular
\begin{equation} \label{estipsik}
\|\Psi_k\|_{H^1(O)} \leq C \; \| \chi_k \|_{L^2(\R^3)} \; .
\end{equation}
Moreover,
\begin{equation} \label{divpsik}
\dv \, \chi_k =0 \quad \Longrightarrow \quad \dv \, \Psi_k =0,
\end{equation}
Next we define
\begin{equation} \label{defwk}
w_k = \rt\, \Psi_k, \quad \in L^2(\R^3)
\end{equation}
whose restriction to $O$ is the good candidate for (\ref{decompuk}). Indeed
$$
\rt \; w_k = \rt \, (\rt \, \Psi_k) = \grad (\dv \, \Psi_k) - \Delta
\Psi_k = \chi_k  \quad \mbox{in } \R^3,
\quad \mbox{(see (\ref{defpsik}) and (\ref{divpsik}))}
$$
which implies in particular
$$
\rt \; w_k = \rt \; u_k \quad \mbox{in }  O.
$$
Before constructing $p_k$, we first check property (i). The fact that $w_k$ 
is bounded in $L^2(O)^3$ results directly from (\ref{estipsik}). Next, we show 
that  $w_k$ is bounded in $H^1(\R^3)$, which implies in particular (i).
Indeed using the Fourier 
transform in $\R^3$ we deduce from (\ref{defwk}) and (\ref{defpsik}) that
 ($\xi$ denotes the dual variable of $x$ and $\widehat{u}$ the Fourier 
transform of $u$):
$$
|\xi|^2 \; |\widehat{w}_k(\xi)|^2 = 
\frac{| \, \xi \times \widehat{\chi}_k(\xi) \, |^2}
{|\xi|^2} \leq \; |\widehat{\chi}_k(\xi)|^2
$$
which yields, by Plancherel's theorem
\begin{equation*}% \label{esti-gradwk}
\| \nabla w_k \|^2_{L^2(\R^3)} \; \leq \; \| \chi_k \|^2_{L^2(\R^3)} 
\leq C \; \|\rt \; u_k\|^2_{L^2(O)}.
\end{equation*}
From now on, we can therefore assume that (up to extracted subsequence) $w_k$ converges in $L^2(O)$.\\[12pt]
Since $\rt (u_k -w_k) = 0$ and $O$ is simply connected, one can construct $p_k$
(unique up to an additive constant) such that $\nabla p_k =u_k -w_k$ (use for
instance Theorem 2.9 of \cite{GR-book86}). Fixing
$p_k$ by imposing that $\int_O \, p_k dx = 0$ gives raise to a bounded sequence
$p_k$ in $H^1(O)$ by the Poincar\'e-Wirtinger inequality.  
Since we further have that $(\dv \,u_k)$ is bounded in $L^2(O)$ and $(w_k)$ is
bounded in $H^1(O)$, then, up to
extracted subsequence, one can assume that $\dv u_k$ is convergent in
$H^{-1}(O)$, $w_k|_{\partial O}$ is convergent in $H^{-\frac{1}{2}}(\partial
O)$ and $p_k|_{\partial O}$ is convergent in $L^2(\partial
O)$.
We shall deduce that $p_k$ is strongly convergent in  $H^1(O)$.
We
first observe that $p_k$ satisfies 
\begin{equation*}% \label{Pbpk}
\left\{ \begin{array}{ll}
- \Delta p_k = \dv \, u_k, & \mbox{in } O,\\[12pt]
\nabla p_k \times n = u_k \times n - w_k \times n, & \mbox{on } \partial O.
\end{array} \right.
\end{equation*}
Let $m$ and $ k$ be two indexes. From 
\begin{equation}\label{RefLem1}
\left\{ \begin{array}{ll}
\Delta (p_k - p_m) = \dv (u_k - u_m) & \mbox{ in } O, \\[10pt]
\nabla (p_k - p_m) \times n = (u_k - u_m) \times n - (w_k
-w_m) \times n & \mbox{ on } \partial O.
\end{array} \right.
\end{equation}
and using the classical  theory for elliptic equations one gets the existence
of a constant $C_1$ such that
\begin{equation}\label{RefLem2}
\|\nabla p_k - \nabla p_m\|_{L^2(\Omega)} \leq C_1 \; \Big( \, \|\dv (u_k -
u_m)\|_{H^{-1}(O)} + \|p_k - p_m\|_{H^\frac{1}{2}(\partial O)} \; \Big).
\end{equation}
On the other hand, using Lemma~\ref{lemMaxpre2} one has
\begin{equation}\label{RefLem3}
\|p_k - p_m\|_{H^\frac{1}{2}(\partial O)} \leq C_2 \left (\| \nabla
(p_k - p_m) \times n \|_{H^{-\frac{1}{2}}(\partial O)} + \| p_k -
p_m\|_{L^2(\partial O)} \right).
\end{equation}
From the second equation of (\ref{RefLem1}), (\ref{RefLem2}) and
(\ref{RefLem3}) it is easily seen that
\begin{equation*}%\label{RefLem4}
\begin{array}{lll}
\|\nabla p_k - \nabla p_m\|_{L^2(\Omega)} \leq   & C_3\left( \|\dv
(u_k - u_m)\|_{H^{-1}(O)} + \|(u_k - u_m) \times
n\|_{H^{-\frac{1}{2}}(\partial O)} \right.\\[10pt]
& \left. \; + \; \quad \|w_k - w_m\|_{H^{-\frac{1}{2}}(\partial
O)} + \| p_k - p_m\|_{L^2(\partial O)} \right).
\end{array}
\end{equation*}
Using assumption (\ref{Hypcompact-lemma}) one concludes $\nabla p_k$ is a Cauchy sequence in
$L^2(\partial O)$.  The result of the lemma is then proved since $u_k = w_k
+ \nabla p_k$.
\proofend

\noindent Lemma \ref{compact-lemma} also applies to domains $\Oint$ that are
not  simply connected. This is proved in the following lemma.
\begin{lemma}\label{NewCompact-lemma}
The result of Lemma~\ref{compact-lemma} applies to  bounded open domains $O
\in R^3$ of class $C^2$.
\end{lemma}
\proof Let $x$ be an arbitrary point in  $\overline{O}$. If $x \in
O$, one defines $U_x$ as a ball centered at $x$ such that
$\overline{U}_x \subset O$. If not, one defines $U_x$ as a
neighborhood of $x$ such that there exits a bijective map $\phi_x : Q \mapsto
U_x$ such that
$$
\phi_x \in C^1(\overline{Q}), ~ \phi^{-1} \in C^1(\overline{U}_x),
~ \phi(Q_+) = U_x \cap O, \mbox{ et } \phi (Q_0) = U \cap
\partial O,
$$
where $Q$ denotes the unit cube of $R^3$, $Q_+ := \{ x \in Q ~|~ x_3 > 0
\}$, and $Q_0 = \{ x \in Q ~|~ x_3 = 0\}$.

With this definition one observes that   $U_x \cap
O$ is a simply connected domain for all $x \in
\overline{O}$. By the compactness of $\overline{O}$ one can extract a finite
covering of  $\overline{O}$ from $\{U_x; x \in \overline O\}$. Let us denote by
$\{ U_i, \; i \in I\}$ this finite covering and consider a partition of unity
$(\theta_i)_{i \in I} \subset C^\infty(R^3)$ subordinated to this covering, i.e.
$$
\supp \theta_i \subset U_i ~ ,~ \sum_{i \in I} \theta_i =1 \mbox{
on } \overline{O}.
$$
Then define $u_n^i := \theta_i \, u_n$ for all $i \in I$. It is easy to see
that for every $i$, the sequence $(u_n^i)$ satisfies the hypotheses of
Lemma~\ref{compact-lemma} with $O$ replaced by  $U_i$. Using a finite diagonal process, one can therefore
assume that 
there exists a subsequence $n_k$ such that 
$$
u_{n_k}^i \mbox{ converges in } L^2(U_i) \mbox{ for all } i \in I.
$$
Consequently, the sequence  $u_{n_k} = \sum_{i \in I} u_{n_k}^i$ is convergent in
$L^2(O)$. \proofend
%\section{The case of the operator $\tilde {\cal D}^{\ptre, 3}$}
\end{appendix}

\tableofcontents

\end{document}